%% file: Versione_4.tex
\DeclareSymbolFont{rsfs}{U}{rsfs}{m}{n}

\documentclass[10pt,a4paper,oneside,leqno,noamsfonts]{amsart}
           \makeatletter
           \newcommand{\mylabel}[2]{#2\def\@currentlabel{#2}\label{#1}}
           \renewcommand\@biblabel[1]{#1.}
           \makeatother
\usepackage[english]{babel}
\usepackage[dvipsnames]{xcolor}
\usepackage{graphicx,pifont}
\usepackage[a4paper,top=2.8cm,bottom=2.8cm,left=2.7cm,
           right=2.7cm,bindingoffset=5mm,marginparwidth=60pt]{geometry}
\usepackage[utf8]{inputenc}
\usepackage{braket,caption,comment,mathtools,stmaryrd}
\usepackage[usestackEOL]{stackengine}
\usepackage{multirow,booktabs,microtype,relsize}
\usepackage[colorlinks,bookmarks]{hyperref} 
      \hypersetup{colorlinks,%
            citecolor=britishracinggreen,%
            filecolor=black,%
            linkcolor=cobalt,%
            urlcolor=black}
      \numberwithin{equation}{section}
\input{macros.tex}

\usepackage{youngtab} 
\newcommand{\minusone}{$-1$}
\newcommand{\ijbox}{$ij$}
\usepackage{ytableau} 
\title[Double  nested Hilbert schemes and the local stable pairs theory of curves]{Double  nested Hilbert schemes and the local stable pairs theory of curves}

\author{Sergej Monavari}
\address{Mathematical Institute, Utrecht University, P.O.~Box 80010 3508 TA Utrecht, The Netherlands}
\email{s.monavari@uu.nl}

\begin{document}
\maketitle

\begin{abstract}
We propose a variation of the classical Hilbert scheme of points --- the \emph{double nested Hilbert scheme of points} --- which parametrizes flags of zero-dimensional subschemes whose nesting is dictated by a Young diagram. Over a smooth quasi-projective curve,  we compute the generating series of topological Euler characteristic of these spaces, by exploiting the combinatorics of reversed plane partitions. Moreover, we realize this moduli space as the zero locus of a section of a vector bundle over a smooth ambient space, which therefore admits a virtual fundamental class.\\
We apply this construction to the  stable pair theory of a local curve, that is the total space of the direct sum of two line bundles over a curve. We show that the invariants localize to virtual intersection numbers on double nested Hilbert scheme of points on the curve, and that the localized contributions to the invariants are controlled by three universal series for every Young diagram, which can be  explicitly determined after the anti-diagonal restriction of the equivariant parameters. Under the anti-diagonal restriction, the invariants are matched with the Gromov-Witten invariants of local curves of Bryan-Pandharipande, as predicted by the MNOP correspondence.\\
Finally, we discuss  $K$-theoretic refinements à la Nekrasov-Okounkov.
\end{abstract}
\section{Introduction}
\subsection{Double nested Hilbert scheme of points}
Let $X$ be a  quasi-projective scheme over $\BC$. We denote by $X^{[n]}$ the \emph{Hilbert scheme of $n$ points} on $X$, which parametrizes 0-dimensional closed subschemes $Z\subset X$ of length $n$. Given a tuple of non-decreasing integers $\mathbf{n}=(n_0 \leq \dots\leq n_d)$, the \emph{nested Hilbert scheme of points} $X^{[\mathbf{n}]}$ parametrizes flags of zero-dimensional subschemes $(Z_0\subset\dots \subset Z_d)$ of $X$, where each $Z_i$ has length $n_i$. The scheme structure of these moduli spaces has been intensively studied in the literature, see for example \cite{Che_cellular_decomposition}.\\
We propose a variation of this moduli space, by parametrizing flags of subschemes nesting in two directions. Let $\lambda$ be a Young diagram and $\mathbf{n}_\lambda=(n_\Box)_{\Box\in \lambda}$ a reversed plane partition, that is a labelling of $\lambda$ by non-negative integers non-decreasing in rows and columns. We denote by $X^{[\mathbf{n}_\lambda]}$ the \emph{double nested Hilbert scheme of points}, the moduli space parametrizing flags of 0-dimensional closed subschemes $(Z_\Box)_{\Box\in \lambda}\subset X$ 
\begin{equation*}
  \begin{tikzcd}
    Z_{00}\arrow[r, phantom, "\subset"]\arrow[d, phantom, "\cap"] &Z_{01}\arrow[d, phantom, "\cap"]\arrow[r, phantom, "\subset"]&Z_{02}\arrow[d, phantom, "\cap"]\arrow[r, phantom, "\subset"]\arrow[d, phantom, "\cap"]&Z_{03}\arrow[r, phantom, "\subset"]\arrow[d, phantom, "\cap"]&\dots\\
      Z_{10}\arrow[r, phantom, "\subset"]\arrow[d, phantom, "\cap"]&Z_{11}\arrow[d, phantom, "\cap"]\arrow[r, phantom, "\subset"]&Z_{12}\arrow[d, phantom, "\cap"]\arrow[r, phantom, "\subset"]&Z_{13}\arrow[d, phantom, "\cap"]\arrow[r, phantom, "\subset"]&\dots\\
        Z_{20}\arrow[d, phantom, "\cap"]\arrow[r, phantom, "\subset"]&Z_{21}\arrow[d, phantom, "\cap"]\arrow[r, phantom, "\subset"]&Z_{22}\arrow[d, phantom, "\cap"]\arrow[r, phantom, "\subset"]&\dots &\\
          \dots&\dots&\dots& &\\
  \end{tikzcd}  
\end{equation*}
where each $Z_\Box$ has length $n_\Box$. If $\lambda$ is a horizontal or vertical Young diagram, the nesting is linear and we recover the usual nested Hilbert scheme of points.\\
The scheme structure
of these moduli spaces is interesting already in dimension one, for a smooth  curve $C$. Cheah proved \cite{Che_cellular_decomposition} that the nested Hilbert scheme $C^{[\mathbf{n}]}$ is smooth,  being isomorphic to a  product of symmetric powers of $C$ via a Hilbert-Chow type morphism. However, as soon as we allow double nestings, $C^{[\mathbf{n}_\lambda]}$ can have several  irreducible components (see Example \ref{exam: singularities}), therefore failing to be smooth.\\
Our first result is a closed formula for the generating series of  topological Euler characteristic of $C^{[\mathbf{n}_\lambda]}$ in terms of the hook-lengths $h(\Box)$ of $\lambda$.
\begin{theorem}[Theorem \ref{thm: euler char double nested}]
Let $C$ be a smooth quasi-projective curve  and $\lambda$ a Young diagram. Then 
\begin{align*}
    \sum_{\mathbf{n}_\lambda} e(C^{[\mathbf{n}_\lambda]})q^{|\mathbf{n}_\lambda|}=\prod_{\Box\in \lambda}(1-q^{h(\Box)})^{-e(C)}.
\end{align*}
\end{theorem}
This is achieved by exploiting  the power structure on the Grothendieck ring of varieties $K_0(\Var_{\BC})$, by which we reduce   to the combinatorial problem of counting the number of reversed plane partitions of a given Young diagram, which was solved by Stanley and Hillman-Grassl \cite{Stanley_theory_application_I_II,HG_reversed_plane_partitions}. Motivic analogues of this formula are studied  in \cite{MR_nested_Quot}.
\subsection{Virtual fundamental class}
The double nested Hilbert scheme $C^{[\mathbf{n}_\lambda]}$ is  in general  singular, making it hard to perform intersection theory. To remedy this, we show that $C^{[\mathbf{n}_\lambda]}$ admits a \emph{perfect obstruction theory} in the sense of Behrend-Fantechi and Li-Tian \cite{BF_normal_cone, LT_virtual_cycle}. In fact, we can (globally!) realize $C^{[\mathbf{n}_\lambda]}$ as the zero locus of a section of a vector bundle over a smooth ambient space.
\begin{theorem}[Theorem \ref{thm: zero locus}]\label{thm: zero locus intro}
Let $C$ be an irreducible smooth quasi-projective curve. There exists a section $s$ of a vector bundle $\CE$ over a smooth scheme $A_{C,\mathbf{n}_\lambda}$ such that 
\[
 \begin{tikzcd}
    &\CE  \arrow[d] \\
   C^{[ \mathbf{n}_\lambda]}\cong Z(s)\arrow[r, hook]
&A_{C,\mathbf{n}_\lambda}. \arrow[u, bend right, swap,  "s"]\end{tikzcd}
\]
 \end{theorem}
 By this construction $C^{[ \mathbf{n}_\lambda]}$ naturally admits a  perfect obstruction theory (see Example \ref{exa:toy model pot}) and in particular carries a virtual fundamental class $[C^{[ \mathbf{n}_\lambda]}]^{\vir}$, which recovers the usual fundamental class in the case where the nesting is linear. We pause a moment to explain this construction in the easiest interesting example, that is for the reversed plane partition 
 \[
 \ytableausetup{centertableaux}
\begin{ytableau}
 n_{00} & n_{01} \\
 n_{10} & n_{11}
\end{ytableau}
\]
 The embedding in the smooth ambient space is
 \begin{align*}
   C^{[ \mathbf{n}_\lambda]}&\hookrightarrow  A_{C,\mathbf{n}_{\lambda}}:=C^{[n_{00}]}\times C^{[n_{10}-n_{00}]}\times C^{[n_{11}-n_{01}]}\times C^{[n_{01}-n_{00}]}\times C^{[n_{11}-n_{10}]},\\
   (Z_{00}, Z_{01}, Z_{10}, Z_{11})&\mapsto (Z_{00}, Z_{10}-Z_{00},Z_{11}-Z_{01},Z_{01}-Z_{00},Z_{11}-Z_{10}).
 \end{align*}
In other words, $A_{C,\mathbf{n}_{\lambda}}$ records the subscheme  in position $(0,0)$ and all possible vertical and horizontal differences of subschemes, where  sum and difference are well-defined by seeing the closed subschemes $Z_{ij}$ as divisors on $C$. At the level of closed points, the image of the embedding is given by all $(Z_{00}, X_1, X_2,Y_1 , Y_2)\in A_{C,\mathbf{n}_{\lambda}}$ such that $X_1+Y_2=Y_1+X_2$ --- again, as divisors. Notice that $X_1+Y_2$ and $Y_1+X_2$ are effective divisors of the same degree, therefore they are equal if and only if one is contained into the other, say $ X_1+Y_2\subset Y_1+X_2$.\\
This  relation is encoded into a section of a vector bundle $\CE$, as we now explain. Denote by $\CX_1, \CX_2, \CY_1, \CY_2$ the universal divisors on $ A_{C,\mathbf{n}_{\lambda}}\times C$ and set 
\begin{align*}
    \Gamma^1&=\CY_{1}+\CX_{2},\\
    \Gamma^2&=\CX_{1}+\CY_{2}.
\end{align*}
The vector bundle $\CE$ is defined as
\begin{align*}
     \CE= \pi_*\oO_{\Gamma^2  }(\Gamma^1),
\end{align*}
where $\pi:A_{C,\mathbf{n}_{\lambda}}\times C\to A_{C,\mathbf{n}_{\lambda}}$ is the projection. The section $s$ of $\CE$ is the one induced --- via $\pi_*$ --- by the section of $\oO_{A_{C,\mathbf{n}_{\lambda}}\times C}(\Gamma^1)$ which vanishes on $\Gamma^1$ and then restricted to $\Gamma^2$.
\subsection{Stable pair invariants of local curves}
Let $C$ be a smooth projective curve and $L_1, L_2$ two line bundles over $C$. We denote by \emph{local curve} the total space $X=\Tot_C(L_1\oplus L_2)$ with its natural $\TT=(\BC^*)^2$-action on the fibers. \\
For $d> 0$ and $n\in \BZ$, we denote by $P_X=P_n(X, d[C])$ the moduli space of stable pairs $[\oO_X\xrightarrow{s} F]\in \Db (X)$ with curve class $d[C]$ and $\chi(F)=n$. The moduli space $P_X$ has a perfect obstruction theory \cite{PT_curve_counting_derived}, but is in general non-proper. Still, the $\TT$-action on $X$ induces one on $P_X$ with proper $\TT$-fixed locus $P_X^\TT$, therefore we can define invariants via Graber-Pandharipande virtual localization \cite{GP_virtual_localization}
\begin{align*}
    \PT_{d,n}(X):=\int_{[ P_X^{\TT}]^{\vir}}\frac{1}{e^{\TT}(N^{\vir})}\in \BQ(s_1,s_2),
\end{align*}
where $s_1,s_2$ are the generators of the $\TT$-equivariant cohomology and $N^{\vir}$ is the virtual normal bundle. We denote its generating series by
\begin{align*}
    \PT_d(X;q):=\sum_{n\in \BZ}q^n\cdot \PT_{d,n}(X)\in \BQ(s_1,s_2)(\!(q)\!).
\end{align*}
Pandharipande-Pixton extensively studied stable pair theory on local curves \cite{PP_descendants_local_curves, PP:local_curves_stationary} using degeneration techniques and relative invariants, focusing on the rationality of the generating series, including the case of descendent insertions. The novelty of this paper is the different approach which only relies on the Graber-Pandharipande localization --- without degenerating the curve $C$ --- and the virtual structure constructed on the double nested Hilbert schemes $C^{[\mathbf{n}_\lambda]}$. This is in particular useful to address the $K$-theoretic generalizations of stable pair invariants (cf. Section \ref{sec: intro K}).

Our main result is  that the generating series $\PT_d(X;q)$ of such invariants is controlled by some universal series   and determine them under the anti-diagonal restriction $s_1+s_2=0$. 
\begin{theorem}[Theorems \ref{thm: full PT invariants}, \ref{thm: generating series with antidiagonal restriction}]\label{thm: main result intro}
There are universal series $A_{\lambda}(q), B_\lambda(q), C_\lambda(q)\in \BQ(s_1,s_2)\llbracket q\rrbracket $ such that
\begin{align*}
\PT_d(X;q)=\sum_{\lambda\vdash d}\left(q^{-|\lambda|}A_{\lambda}(q)\right)^{g-1}\cdot\left(q^{-n(\lambda)} B_{\lambda}(q)\right)^{\deg L_1}\cdot\left(q^{-n(\overline{\lambda})} C_{\lambda}(q)\right)^{\deg L_2},
\end{align*}
where  $\overline{\lambda}$ is the conjugate partition of $\lambda$, $n(\lambda)=\sum_{i=0}^{l(\lambda)}i\cdot\lambda_{i}$ and $g=g(C)$.
Moreover, under the anti-diagonal restriction $s_1+s_2=0$
\begin{align*}
   A_{\lambda}(q,s_1,-s_1)&=(-s_1^2)^{|\lambda|}\cdot\prod_{\Box\in \lambda}h(\Box)^2,\\
    B_{\lambda}(-q,s_1,-s_1)&=(-1)^{n(\lambda)}\cdot s_1^{-|\lambda|}\cdot\prod_{\Box\in \lambda}h(\Box)^{-1}\cdot\prod_{\Box\in \lambda}(1-q^{h(\Box)}),\\
    C_{\lambda}(-q,s_1,-s_1)&=(-1)^{n(\overline{\lambda})}\cdot (-s_1)^{-|\lambda|}\cdot\prod_{\Box\in \lambda}h(\Box)^{-1}\cdot\prod_{\Box\in \lambda}(1-q^{h(\Box)}).
\end{align*}
\end{theorem}
We sketch now the main steps required in proving Theorem \ref{thm: main result intro}.
\subsection{Proof of the main theorem}
The connected components of the $\TT$-fixed locus $P_n(X,d[C])^\TT$  are double nested Hilbert schemes of points $C^{[\mathbf{n}_\lambda]}$, for suitable reversed plane partitions $\mathbf{n}_\lambda$ and Young diagram $\lambda$. In fact,  pushing forward via $X\to C$ a $\TT$-fixed stable pair $[\oO_X\xrightarrow{s}F]$, corresponds a decomposition $\bigoplus_{(i,j)\in \BZ^2}[\oO_C\xrightarrow{s_{ij}} F_{ij}]$ on $C$, where every $F_{ij}$ is a line bundle with section $s_{ij}$. These data produce divisors $Z_{ij}\subset C$ satisfying the nesting conditions dictated  by $\lambda$, in other words an element of $ C^{[\mathbf{n}_\lambda]}$. \\

On each connected component, there is an induced virtual fundamental class $[C^{[\mathbf{n}_\lambda]}]^{\vir}_{\PT}$, coming from the deformation of stable pairs. This virtual cycle coincides with the one constructed by the zero-locus construction of Theorem \ref{thm: zero locus intro}.
By determining the class in $K$-theory of the virtual normal bundle,  stable pair invariants on $X$ are reduced to ($\TT$-equivariant) virtual intersection numbers on $ C^{[\mathbf{n}_\lambda]}$, namely
\begin{align}\label{eqn: intro inv double nested}
    \int_{[C^{[\mathbf{n}_\lambda}]^{\vir}}e^{\TT}(-N^{\vir}_{C,L_1,L_2})\in \BQ(s_1, s_2).
\end{align}
The generating series of these invariants, for every fixed Young diagram $\lambda$, is controlled by three universal series (Theorem \ref{thm: universal series})
\begin{align*}
    \sum_{\mathbf{n}_\lambda}q^{|\mathbf{n}_\lambda|} \int_{[C^{[\mathbf{n}_\lambda]}]^{\vir}}e^{\TT}(-N_{C,L_1,L_2}^{\vir})=A_{\lambda}^{g-1}\cdot B_{\lambda}^{\deg L_1}\cdot C_{\lambda}^{\deg L_2}\in \BQ(s_1,s_2)\llbracket q \rrbracket.
\end{align*}
This universal structure is proven by following the strategy of \cite{EGL_cobordism}. In fact, these invariants are \emph{multiplicative} on triples of the form $(C, L_1, L_2)=(C'\sqcup C'', L'_1\oplus L''_1, L'_2\oplus L''_2)$ and are polynomial in   the Chern numbers of $(C, L_1, L_2)$. The latter is obtained by pushing the virtual intersection number to $C^{[\mathbf{n}_\lambda]}$ on the smooth ambient space $A_{C, \mathbf{n}_\lambda}$ --- a product of symmetric powers of $C$ ---
and later to a product of Jacobians $\Pic^{n_i}(C)$, where the integrand is a polynomial on well-behaved cohomology classes.\\

By the universal structure  any computation is reduced  to a basis of the three-dimensional $\BQ$-vector space of Chern numbers  of triples $(C, L_1, L_2)$. A simple basis consists of the Chern numbers of  $(\BP^1,\oO, \oO) $ and any two $(\BP^1, L_1, L_2)$ with $L_1\otimes L_2 =K_{\BP^1}$. In both cases, the invariants are explictly determined under the anti-diagonal restriction $s_1+s_2=0$ by further applying the virtual localization formula.
\subsection{Toric computations}
The $\BC^*$-action on $\BP^1$ canonically lifts to the double nested Hilbert scheme ${\BP^1}^{[\mathbf{n}_\lambda]}$, with only finitely many $\BC^*$-fixed points, therefore we can further $\BC^*$-localize the invariants \eqref{eqn: intro inv double nested} to obtain
\begin{align*}
    \int_{\left[{\BP^{1}}^{[\mathbf{n}_\lambda]}\right]^{\vir}}
    e^\TT(-N_{\BP^{1},L_1,L_2}^{\vir})
    &=\left.\left(\sum_{\underline{Z}\in {\BP^{1}}^{[\mathbf{n}_\lambda],\BC^*}}e^{\TT\times \BC^*}(-T^{\vir}_{\underline{Z}}-N^{\vir}_{\BP^{1},L_1,L_2, \underline{Z} }) \right)\right|_{s_3=0},
\end{align*}
where $s_3$ is the generator of the $\BC^*$-equivariant cohomology and $T^{\vir}_{\underline{Z}}$ is the virtual tangent bundle of $ {\BP^{1}}^{[\mathbf{n}_\lambda]}$ at the fixed point $\underline{Z}$.\\
Under the anti-diagonal restriction $s_1+s_2=0$, this translates the computation of the invariants into a purely combinatorial problem, which we explictly solve in the trivial vector bundle case $L_1=L_2=\oO_{\BP^1}$ and in the Calabi-Yau case $L_1\otimes L_2=K_{\BP^1}$. A few  remarks are in order. In the trivial vector bundle case, the solution is equivalent to the  vanishing
\begin{align*}
    \left.\int_{[{\BP^1}^{[\mathbf{n}_\lambda]}]^{\vir}}e^{\TT}(-N_{\BP^1,\oO,\oO}^{\vir})\right|_{s_1+s_2=0}=0,
\end{align*}
for every reversed plane partition of positive size $|\mathbf{n}_\lambda|>0$. This relies on the vanishing $e^{\TT\times \BC^*}(-T^{\vir}_{\underline{Z}} - N^{\vir}_{\BP^{1},L_1,L_2, \underline{Z} })=0 $, which comes from a simple vanishing property of the topological vertex in stable pair theory proved in \cite{MPT_curves_K3}.\\
In the Calabi-Yau case, the invariants turn out to be  \emph{topological}, under the anti-diagonal restriction.
\begin{theorem}[Theorem \ref{thm: g=0 case}]\label{thm: intro g=0}
Let $X$ be Calabi-Yau. Then the generating series of the invariants \eqref{eqn: intro inv double nested} coincides, up to a sign,  with the generating series of the topological Euler characteristic
\begin{align*}
   \sum_{\mathbf{n}_\lambda}q^{|\mathbf{n}_\lambda|}\cdot \left.\left( \int_{\left[{\BP^{1}}^{[\mathbf{n}_\lambda]}\right]^{\vir}}e^{\TT}(-N^{\vir}_{\BP^1, L_1, L_2})\right)\right|_{s_1+s_2=0}=(-1)^{\deg L_1(c_\lambda+|\lambda|)+ |\lambda|}\cdot \sum_{\mathbf{n}_\lambda}(-q)^{|\mathbf{n}_\lambda|}e\left({\BP^1}^{[\mathbf{n}_\lambda]}\right),
\end{align*}
where $c_\lambda=\sum_{(i,j)\in \lambda}(j-i)$.
\end{theorem}
This happens as, under the anti-diagonal restriction, each  $\BC^*$-fixed point $\underline{Z}$ contributes with a sign
\begin{align}\label{eqn:intro sign}
   \left. e^{\TT\times \BC^*}(-T^{\vir}_{\underline{Z}}-N^{\vir}_{\BP^{1},L_1,L_2, \underline{Z} })\right|_{s_1+s_2=0}=(-1)^{\deg L_1(c_\lambda+|\lambda|)+ |\lambda|+|\mathbf{n}_\lambda|},
\end{align}
which is independent of $\underline{Z}$ and the invariants amount to a (signed) count of the $\BC^*$-fixed points. It is not a priori clear how to obtain the the same sign through  the vertex formalism for stable pairs developed by Pandharipande-Thomas \cite{PT_vertex}.\\
Nevertheless, the topological nature of the invariants in the Calabi-Yau case is not surprising also for a non-toric curve $C$. If $X$ is Calabi-Yau and $P_n(X, d[C])$ is proper --- which happens  only in rare cases --- the anti-diagonal restriction would compute its virtual Euler characteristic and Behrend's  weighted Euler characteristic, which is a purely topological invariant of a scheme with a symmetric perfect obstruction theory \cite{Beh_DT_via_microlocal}. 
\subsection{Gromov-Witten/stable pairs correspondence}
In the seminal work \cite{MNOP_1}, a conjectural correspondence - known as the MNOP conjecture - between Gromov-Witten invariants and Donaldson-Thomas invariants of projective threefolds is formulated, proven for toric varieties in \cite{MNOP_1, MNOP_2, MOOP_GW/DT_toric} for primary insertions. By defining the GW/DT invariants via equivariant residues, the conjecture has been extended to local curves in \cite{BP_local_GW_curves} and proven by combining the results of \cite{BP_local_GW_curves, OP_local_theory_curves}. \\
Stable pair invariants were later introduced by Pandharipande-Thomas \cite{PT_curve_counting_derived} to give a more natural geometric interpretation of the MNOP conjecture through the DT/PT correspondence proved by Toda and Bridgeland in \cite{Toda_curve_counting_DT_PT, Bri_Hall_algebras_curve_counting} using wall-crossing and Hall algebra techniques. The Gromov-Witten/stable pairs correspondence has been subsequently extended to include descendent insertions and to quasi-projective varieties whenever invariants can be defined through virtual localization. The correspondence had been confirmed by Pandharipande-Pixton for  Calabi-Yau and Fano complete intersections in product of projective spaces and toric varieties \cite{PP_GW/PT_quintic, PP_GW/PT_toric} and had been recently addressed in \cite{OOP_GW/PT_descendents}. See \cite{Pandha_descendents_stable_pairs_3folds} for a complete survey on the subject.
\subsection{The local GW theory of curves}
For $X=\Tot_C(L_1\oplus L_2)$  a local curve, let $\overline{M}^\bullet_h(X, d[C])$ denote the moduli space of stable maps (with possibly disconnected domain) of genus $h$ and degree $d[C]$. Define the partition function of Gromov-Witten invariants of $X$ (with a shifted exponent)
\begin{align*}
    \GW_d(g|\deg L_1, \deg L_2;u)=u^{2-2g+\deg L_1+\deg L_2}\sum_{h\in \BZ}u^{2h-2}\int_{[ \overline{M}^\bullet_h(X, d[C])^\TT]^{\vir}}\frac{1}{e^{\TT}(N^{\vir})}\in \BQ(s_1, s_2)(\!( u )\!),
\end{align*}
where the dependence is only on the genus $g=g(C)$, the degrees of the line bundles and the degree $d$. The  Gromov-Witten theory of local curves had been solved by Bryan-Pandharipande \cite[Theorem 7.1]{BP_local_GW_curves} using a TQFT approach. Moreover  they deduced an  explicit closed formula for the partition function under the anti-diagonal restriction $s_1+s_2=0$.
\begin{theorem}[Bryan-Pandharipande]\label{thm: GW series}
The partition function of Gromov-Witten invariants satisfies
\begin{multline*}
    \left.\GW_{d}(g|\deg L_1,\deg L_2;u)\right|_{s_1+s_2=0}=(-1)^{d(g-1-\deg L_2)}s_1^{d(2g-2-\deg L_1-\deg L_2)}\cdot\\\sum_{\lambda\vdash d} Q^{\frac{1}{2}c_\lambda(\deg L_1-\deg L_2)}\prod_{\Box\in \lambda}h(\Box)^{2g-2-\deg L_1-\deg L_2}\cdot i^{-\deg L_1-\deg L_2}\left(Q^{\frac{h(\Box)}{2}}-Q^{-\frac{h(\Box)}{2}}\right)^{\deg L_1+\deg L_2},
\end{multline*}
where we set $Q=e^{i u}$ and $i=\sqrt{-1}$.
\end{theorem}
With this explicit expression it is immediate to check the Gromov-Witten/stable pairs correspondence under the anti-diagonal restriction.
\begin{corollary}[Corollary \ref{cor: GW/PT correspondence}]
Let $X$ be a local curve. Under the anti-diagonal restriction $s_1+s_2=0$ the GW/stable pair correspondence holds
\begin{align*}
    (-i)^{d(2-2g+\deg L_1+\deg L_2)}\cdot \GW_d(g|\deg L_1, \deg L_2;u)=(-q)^{-\frac{1}{2}\cdot d(2-2g+\deg L_1+\deg L_2)}\PT_d(X, q),
\end{align*}
after the change of variable $q=-e^{iu}$.
\end{corollary}
\subsection{$K$-theoretic refinement}\label{sec: intro K}
$K$-theoretic refinement of Donaldson-Thomas theory and stable pair theory attracted much attention recently, both in Mathematics and String Theory: see for example \cite{Tho_equivariant_K_theory, Afg_refinement_KT, Arb_K-theo_surface, FMR_higher_rank} for Calabi-Yau threefolds,  \cite{Nek_magnificient_4, NP_colors, CKM_K_theoretic} for Calabi-Yau fourfolds and \cite{NO_membranes_and_sheaves, Okounk_Lectures_K_theory,KOO_2_legDT, Oko_Takagi} for local curves.\\
A scheme $X$ with a perfect obstruction theory is endowed not only with a virtual fundamental class, but also with a \emph{virtual structure sheaf} $\oO^{\vir}_X\in K_0(X)$. If $X$ is proper, $K$-theoretic invariants are simply of the form 
\[
\chi(X, \oO^{\vir}_X\otimes V)\in \BZ,
\]
where $V\in K_0(X)$. If $X$ is a local curve the moduli space of stable pairs $P_X$ is in general not proper and $K$-theoretic stable pair invariants  are defined by virtual localization \cite{FG_riemann_roch} on the proper $\TT$-fixed locus $P_X^\TT$, that is one set
\[
\chi(P_X, \oO^{\vir}_{P_X}\otimes V):= \chi\left(P_X^{\TT}, \frac{\oO^{\vir}_{P_X^{\TT}}\otimes V|_{P_X^\TT}}{\Lambda^\bullet N^{\vir, *}} \right)\in \BQ(\tf_1, \tf_2).
\]
In Section \ref{sec: K-th} we show that, also in the $K$-theoretic setting, the invariants are controlled by  universal series.\\
The naive generalization of cohomological invariants is for $V=\oO_X$, that is no insertions. However, we learn from Nekrasov-Okounkov \cite{NO_membranes_and_sheaves} that it is more natural to consider the \emph{twisted virtual structure sheaf}
\[  
\widehat{\oO}^{\vir}=\oO^{\vir}\otimes K_{\vir}^{1/2},
\]
where $ K_{\vir}^{1/2}$ is a square root\footnote{This square root may not exist as a line bundle, but it does exists as a class in $K$-theory after inverting 2.} of the \emph{virtual canonical bundle}. Denote by $\PT^{\widehat{K}}_d(X;q)$ the generating series of $K$-theoretic invariants with $V=K_{\vir}^{1/2}$.
\begin{theorem}[Corollary \ref{cor: K-theoretic NO universal}]
There exist universal series  $ A_{\widehat{K},\lambda}(q), B_{\widehat{K},\lambda}(q), C_{\widehat{K},\lambda}(q)\in \BQ(\tf_1^{1/2}, \tf_2^{1/2})\llbracket q \rrbracket$ such that
\begin{align*}
\PT^{\widehat{K}}_d(X;q)=\sum_{\lambda\vdash d}\left(q^{-|\lambda|}A_{\widehat{K},\lambda}(q)\right)^{g-1}\cdot\left(q^{-n(\lambda)} B_{\widehat{K},\lambda}(q)\right)^{\deg L_1}\cdot\left(q^{-n(\overline{\lambda})} C_{\widehat{K},\lambda}(q)\right)^{\deg L_2}.
\end{align*}
Moreover, the universal series are explicitly computed under $\tf_1\tf_2=1$.
\end{theorem}
We are  not aware of a $K$-theoretic Gromov-Witten refinement for which a \emph{refined} GW/stable pairs correspondence holds.
\subsection*{Acknowledgements}
I am grateful to Martijn Kool for suggesting this problem and for many helpful discussions.  I wish to thank Jan-Willem van Ittersum and Andrea Ricolfi for useful discussions. S.M. is supported by NWO grant TOP2.17.004.
\subsection*{Conventions}
All \emph{schemes} are  separated  and of finite type over $\BC$.  We let $K^0(X)$ be the K-group of vector bundles on $X$. When $X$ carries an action by an algebraic torus $\TT$, we let $K^0_\TT(X)$ be the K-group of $\TT$-\emph{equivariant} vector bundles on $X$. Similarly, we let $K_0(X)$ denote the K-group of \emph{coherent sheaves} on $X$, and we let $K_0^{\TT}(X)$ be the K-group of (the abelian category of) $\TT$-\emph{equivariant} coherent sheaves on $X$. When $X$ is smooth, the natural $\BZ$-linear map $K^0(X) \to K_0(X)$, resp.~$K^0_{\TT}(\pt)$-linear map $K^0_{\TT}(X) \to K_0^{\TT}(X)$, is an isomorphism.
Chow groups $A^\ast (X)$ and cohomology groups $H^\ast(X)$ are taken with rational coefficients.
We use $(\cdot)^\vee$ for the derived dual of complexes and $(\cdot)^*$ for the (underived) dual of coherent sheaves. For clarity of exposition, we suppress various pullback maps, whenever they are clear from the context.
\section{Double nested Hilbert schemes}
\subsection{Young diagrams}\label{sec: notation partitions}
By definition, a \emph{partition} $\lambda$ of $d\in \BZ_{\geq 0}$ is a finite sequence of positive integers 
\[
\lambda=(\lambda_0\geq \lambda_1\geq \lambda_2\geq \dots)
\]
where 
\[
|\lambda|=\sum_{i}\lambda_i=d.
\]
The number of parts of $\lambda$ is called the \emph{length} of $\lambda$ and is denoted by $l(\lambda)$. A partition $\lambda$ can be equivalently described by its associated \emph{Young diagram}, which is the collection of $d$ boxes  in $\BZ^2$ located at $(i,j)$ where $0\leq j< \lambda_{i}$.\footnote{This notation was borrowed by \cite[Sec. 3.1]{BP_local_GW_curves}; however, in our conventions, $(i,j)$ labels the box's corner closest to the origin.}\\
Given a partition $\lambda$, a \emph{reversed plane partition} $\mathbf{n}_\lambda=(n_{\Box})_{\Box\in \lambda}\in \mathbb{Z}_{\geq 0}$ is a collection of non-negative integers such that $n_{\Box}\leq n_{\Box'}$ for any $\Box,\Box'\in \lambda$ such that  $\Box\leq \Box'$.
In other words, a reversed plane partition is a Young diagram labelled with non-negative integers which are non-decreasing in rows and columns. The \emph{size} of a reversed plane partition is
\begin{align*}
    |\mathbf{n}_\lambda|=\sum_{\Box\in \lambda}n_{\Box}.
\end{align*}
\begin{figure}[!h]
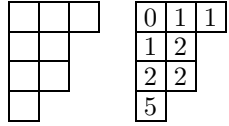

    \centering
   \yng(3,2,2,1) \quad \young(011,12,22,5)
    \caption{On the left, a Young diagram of size 8. On the right, a reversed plane partition of size 14.}
    \label{fig:my_label}
\end{figure}
The \emph{conjugate partition} $\overline{\lambda}$ is obtained by reflecting the Young diagram of $\lambda$ about the $i=j$ line.

In the paper we will require the following standard quantities. Given a box in the Young diagram $ \box\in \lambda$, define the \emph{content} $c(\Box)=j-i$ and the \emph{hooklength} $h(\Box)=\lambda_{i}+\overline{\lambda}_{j}-i-j-1$. The total content
\[ c_\lambda=\sum_{\Box \in \lambda}c(\Box) \]
satisfies the following identities (cf \cite[pag. 11]{MacD_symmetric_functions_Hall}):
\begin{align}\label{eqn: combinatorial identities}
    \sum_{\Box \in \lambda}h(\Box)=n(\lambda)+n(\overline{\lambda})+|\lambda|, \quad c_\lambda=n(\overline{\lambda})-n(\lambda),
\end{align}
where 
\[
n(\lambda)=\sum_{i=0}^{l(\lambda)}i\cdot\lambda_{i}.
\]
For any Young diagram $\lambda$ there is an associated  graph, where  any box of $\lambda$ corresponds to a vertex and  any face common to two boxes correspond to an edge connecting the corresponding vertices. A \emph{square} of this graph is a circuit made of four different edges.
\begin{figure}[!h]
    \centering
   \[\yng(3,2,2,1) \quad \begin{tikzcd}
     \bullet \arrow[r, dash] \arrow[d, dash] & \bullet \arrow[r, dash]  \arrow[d, dash]&\bullet\\
     \bullet \arrow[r, dash] \arrow[d, dash] & \bullet  \arrow[d, dash]& \\
     \bullet \arrow[r, dash] \arrow[d, dash]& \bullet & \\
     \bullet &  & 
   \end{tikzcd} \]
    \caption{A Young diagram and its associated graph, with 8 vertices, 9 edges and 2 squares.}
    \label{fig: graph}
\end{figure}
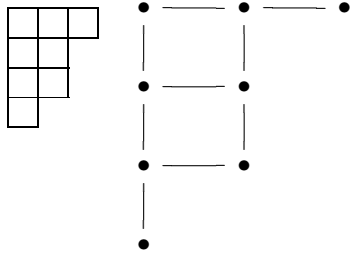
\begin{lemma}\label{lemma: vertices and edges of graph}
Let $\lambda$ be a Young diagram and denote by $V,E,Q$ respectively the number of  vertices, edges and squares of the associated graph. Then 
  \[V-E+Q-1=0. \]
\end{lemma}
 \begin{proof}
 We prove the claim by induction on the size of $\lambda$. If $|\lambda|=1$, this is clear. Suppose it holds for all $\lambda$ with $|\lambda|\leq n-1$. Then we construct $\lambda$ of size $n$ by adding a box with  lattice coordinates $(i,j)$ to a Young diagram $ \Tilde{\lambda}$ of size $n-1$. There are two possibilities: either one of $i,j$ is zero, so  we added one vertex and one edge, or both $i,j$ are non-zero, so we  added one vertex, one square and two edges. In both cases the claim is proved.
 \end{proof}
\subsection{Double nested Hilbert schemes}
Let $X$ be a projective scheme and $\oO(1) $ a fixed ample line bundle.   The \emph{Hilbert polynomial} of a closed subscheme $Y\subset X$ is defined by 
\[m\mapsto \chi(\oO_Y\otimes \oO(m)).\]
Given a polynomial $p(m)$, the \emph{Hilbert scheme} is the moduli space parametrizing closed subschemes  $Y\subset X$ with Hilbert polynomial $p(m)$, which is representable by a projective scheme (e.g. by  \cite{Gro_existence_Hilbert}). We consider here a more general situation, where we replace  closed subschemes by flags of closed subschemes, satisfying certain nesting conditions dictated by Young diagrams.\\

Let $\lambda$ be a Young diagram and   $\mathbf{p}_\lambda=(p_{\Box})_{\Box\in \lambda}\in \mathbb{Z}[x]$ be a collection of polynomials indexed by $\lambda$. If all $p_{\Box}$ are  non-negative integers which are non-decreasing in rows and columns, $\mathbf{p}_\lambda=\mathbf{n}_\lambda$ is a reversed plane partition.
\begin{definition}\label{def: double nested}
Let $X$ be a projective scheme and $\mathbf{p}_\lambda$ as above. The \emph{double nested Hilbert functor} of $X$ of type $\mathbf{p}_\lambda$ is the moduli functor
\[
 \underline{\Hilb}^{\mathbf{p}_{\lambda}}(X): \Sch^{\op}\to \Sets,
\]
\[  T\mapsto \left\{ \begin{tabular}{l|l}$
(\CZ_{\Box})_{\Box\in \lambda}\subset X\times T$ & $\CZ_{\Box}\mbox{ a } T\mbox{-flat closed subscheme with Hilbert polynomial } p_{\Box},$\\ 
    &  $\mbox{such that } \CZ_{\Box} \subset \CZ_{\Box'} \mbox{ for } \Box\leq \Box'.$
    \end{tabular}
    \right\}. \]
\end{definition}
\begin{remark}
If $|\lambda|=1$ we recover the classical  Hilbert scheme, while if $\lambda$ is a horizontal (or vertical) Young diagram we recover the \emph{nested Hilbert scheme}, already widely studied in the literature.
\end{remark}
\begin{prop}\label{prop: representability}
Let $X$ be a projective scheme and $\mathbf{p}_\lambda$ as above. Then $\underline{\Hilb}^{\mathbf{p}_{\lambda}}(X)$ is representable by a projective scheme $\Hilb^{\mathbf{p}_{\lambda}}(X) $, which we call the \emph{double nested Hilbert scheme}.
\end{prop}
\begin{proof}
We prove our claim in the case $\mathbf{p}_{\lambda}$ is
\[
 \ytableausetup{centertableaux}
\begin{ytableau}
 p_{00} & p_{01} \\
 p_{10} & p_{11}
\end{ytableau}
\]
as the general case will follow by an analogous reasoning. There are forgetful maps between nested Hilbert functors
\begin{align*}
   &\underline{\Hilb}^{[p_{00}, p_{01},p_{11}]}(X)\to \underline{\Hilb}^{[p_{00},p_{11}]}(X)\\
    &\underline{\Hilb}^{[p_{00}, p_{10},p_{11}]}(X)\to \underline{\Hilb}^{[p_{00},p_{11}]}(X),
\end{align*}
which forget the second subscheme of the corresponding flag. Consider their fiber product
\[
\begin{tikzcd}
  \underline{\Hilb}^{[p_{00}, p_{01},p_{11}]}(X)\times_{\underline{\Hilb}^{[p_{00},p_{11}]}(X)} \underline{\Hilb}^{[p_{00}, p_{10},p_{11}]}(X)\arrow[d]\arrow[r]&\underline{\Hilb}^{[p_{00}, p_{01},p_{11}]}(X)\arrow[d]\\
  \underline{\Hilb}^{[p_{00}, p_{01},p_{11}]}(X)\arrow[r]& \underline{\Hilb}^{[p_{00},p_{11}]}(X).
\end{tikzcd}
\]
There is an obvious morphism of functors
\[
\underline{\Hilb}^{\mathbf{p}_{\lambda}}(X)\to \underline{\Hilb}^{[p_{00}, p_{01},p_{11}]}(X)\times_{\underline{\Hilb}^{[p_{00},p_{11}]}(X)} \underline{\Hilb}^{[p_{00}, p_{10},p_{11}]}(X),
\]
which is easily checked to be an isomorphism by comparing each flat family of flags over every scheme $T$. We conclude by the fact that the nested Hilbert functors (and their fiber products) are representable by a projective scheme by \cite[Thm. 4.5.1]{Ser_deformation}.
\end{proof}
Thanks to representability, double nested Hilbert schemes are equipped with universal subschemes, for any $\Box\in \lambda $,
\begin{align*}
    \CZ_{\Box}\subset X\times \Hilb^{\mathbf{p}_{\lambda}}(X),
\end{align*}
such that the fiber over a point $\underline{Z}=(Z_\Box)_{\Box\in \lambda}\in \Hilb^{\mathbf{p}_{\lambda}}(X)$ is
\[
{\CZ_\Box}|_{\underline{Z}}=Z_\Box\subset X.
\]
\begin{remark}
If  $ \mathbf{p}_{\lambda}=\mathbf{n}_{\lambda}$, Definition \ref{def: double nested}  generalizes to $X$ quasi-projective. In fact, let $X\subset \overline{X} $ be any compactification of $X$. We define the \emph{double nested Hilbert scheme points} as the open subscheme
\begin{align*}
X^{[\mathbf{n}_\lambda]}:=\Hilb^{\mathbf{n}_{\lambda}}(X)\subset \Hilb^{\mathbf{n}_{\lambda}}(\overline{X})
\end{align*}
consisting of the 0-dimensional subschemes supported  on $X\subset \overline{X}$.
\end{remark}
Double nested Hilbert schemes of points are rarely smooth varieties. Some smooth examples consist of
\begin{itemize}
    \item $|\lambda|=1$, $X$ a smooth quasi-projective curve or surface (see e.g.  \cite{Nak_lectures_Hilb_schemes}),
    \item $\lambda$ a vertical/horizontal Young diagram, $X$ a smooth quasi-projective curve (see e.g.   \cite{Che_cellular_decomposition}).
\end{itemize}
In general,  $X^{[\mathbf{n}_\lambda]}$ is singular even for $X$ a smooth quasi-projective curve.
\begin{example}\label{exam: singularities}
Let $C$ be a smooth curve and consider the reversed plane partition $\mathbf{n}_\lambda$
\begin{figure}[!h]
    \centering
   \young(01,12)
\end{figure}

There are two types of flags of divisors, of the form
\[
\begin{tikzcd}
 \varnothing\arrow[r, phantom, "\subset"]\arrow[d, phantom, "\cap"] & P\arrow[d, phantom, "\cap"]\\
  Q\arrow[r, phantom, "\subset"]&P+Q,
\end{tikzcd}
\qquad 
\begin{tikzcd}
 \varnothing\arrow[r, phantom, "\subset"]\arrow[d, phantom, "\cap"] & P\arrow[d, phantom, "\cap"]\\
  P\arrow[r, phantom, "\subset"]&P+Q,
\end{tikzcd}
\]
where $P,Q\in C$.
Therefore its reduced scheme structure consists of two irreducible components $C\times C\cup C\times C$,  intersecting at the diagonals of $C\times C$.
\end{example}
Singularities make it hard to perform intersection theory on $X^{[\mathbf{n}_\lambda]}$. To remedy this we construct, in special cases, \emph{virtual fundamental classes} in $A_*(X^{[\mathbf{n}_\lambda]})$. We  briefly recall the language of perfect obstruction theories of Behrend-Fantechi  and Li-Tian \cite{BF_normal_cone,LT_virtual_cycle}.
\subsection{Perfect obstruction theories}\label{sec: pot}
A \emph{perfect obstruction theory} on a scheme $X$ is the datum of a morphism 
\[
\phi\colon \BE\to\BL_X
\]
in $\derived^{[-1,0]}(X)$, where $\BE$ is a perfect complex of perfect amplitude contained in $[-1,0]$, such that $h^0(\phi)$ is an isomorphism and $h^{-1}(\phi)$ is surjective. Here, $\BL_X = \tau_{\geq -1}L_X^\bullet$ is the cut-off at $-1$ of the full cotangent complex $L_X^\bullet \in \derived^{[-\infty,0]}(X)$. A perfect obstruction theory is called \emph{symmetric} (see \cite{BF_symm_pot}) if there exists an isomorphism $\theta\colon \BE \xrightarrow{\sim} \BE^\vee[1]$ such that $\theta = \theta^\vee[1]$.
The \emph{virtual dimension} of $X$ with respect to $(\BE,\phi)$ is the integer $\vd = \rk \BE$. This is just $\rk E^0 - \rk E^{-1}$ if one can write $\BE = [E^{-1}\to E^0]$. \\
A perfect obstruction theory determines a cone
\[
\mathfrak C \into E_1 = (E^{-1})^*.
\]
Letting $\iota\colon X \into E_1$ be the zero section of the vector bundle $E_1$, the induced \emph{virtual fundamental class} on $X$ is the refined intersection
\[
[X]^{\vir} = \iota^![\mathfrak C]\, \in\, A_{\vd}(X). 
\]
By a result of Siebert \cite[Thm.~4.6]{Sie_virtual_fund_classes}, the virtual fundamental class depends only on the K-theory class of $\BE$.
\begin{example}\label{exa:toy model pot}
Let $\iota: Z=Z(s)\into A$ be the zero locus of a section $s\in \Gamma(A,\CE)$, where $\CE$ is a vector bundle over a smooth quasi-projective variety $A$. Then there exists a induced perfect obstruction theory on $Z$
\begin{equation*}\label{symmetric_POT_quot}
\begin{tikzpicture}[baseline=(current  bounding  box.center)]
\node (N1) at (-1.98,0.95) {$\mathbb E$};
\node (N2) at (-1.37,0.94) {$=$};
\node (N3) at (-1.97,-0.88) {$\BL_{Z}$};
\node (N4) at (-1.38,-0.88) {$=$};
\node (O1) at (-0.1,0.93) {$\big[\CE^*|_{Z}$};
\node (O2) at (2.99,0.93) {$\Omega_{A}|_{Z}\big]$};
\node (O3) at (-0.1,-0.88) {$\big[ \CI/\CI^2$};
\node (O4) at (2.99,-0.88) {$\Omega_{A}|_{Z}\big]$};
\path[commutative diagrams/.cd, every arrow, every label]
(N1) edge node[swap] {$\phi$} (N3)
(O1) edge node {$\dd s^*$} (O2)
(O1) edge node[swap] {$s^*$} (O3)
(O3) edge node {$\dd$} (O4)
(O2) edge node {$\mathrm{id}$} (O4);
\end{tikzpicture}
\end{equation*}
in $\derived^{[-1,0]}(Z)$, where we represented the truncated cotangent complex by means of the exterior derivative $\dd$ constructed out of the ideal sheaf $\CI \subset \oO_{Z}$ of the inclusion $Z \into A$. Moreover
\begin{align*}
    \iota_*[Z]^{\vir}=e(\CE)\cap [A]\in A_*(A),
\end{align*}
where $e(\cdot)$ denotes the Euler class.
\end{example}
\subsection{Points on Curves}\label{sec: pot on points on curves}
Let $C$ be an irreducible smooth quasi-projective curve and $ \mathbf{n}_{\lambda}$ a reversed plane partition.  In this section we show that $C^{[\mathbf{n}_{\lambda}]}$ is the zero locus of a section of a vector bundle over a smooth ambient space, and therefore admits a perfect obstruction theory as in Example \ref{exa:toy model pot}. \\
We define 
\begin{align*}
A_{C,\mathbf{n}_{\lambda}}=C^{[n_{00}]}\times \prod_{\substack{(i,j)\in \lambda\\ i\geq 1}} C^{[n_{ij}-n_{i-1,j}]}\times \prod_{\substack{(l,k)\in \lambda\\ k\geq 1}} C^{[n_{lk}-n_{l,k-1}]}.
\end{align*}
As $ C^{[n]}\cong C^{(n)}$ is a symmetric product via the Hilbert-Chow morphism,  $A_{C,\mathbf{n}_{\lambda}}$ is a smooth quasi-projective variety of dimension \[\dim(A_{C,\mathbf{n}_{\lambda}})=n_{00}+\sum_{\substack{(i,j)\in \lambda\\ i\geq 1}}(n_{ij}-n_{i-1,j})+\sum_{\substack{(l,k)\in \lambda\\ k\geq 1}}(n_{lk}-n_{l,k-1}).\]
To ease the notation, we denote its elements  by $\underline{Z}=((Z_{00}, X_{ij},Y_{lk}))_{ij,lk}\in A_{C,\mathbf{n}_{\lambda}}$, where
$Z_{00}\subset C$ is a divisor of length $n_{00}$ and $X_{ij}\subset C$ (resp. $Y_{lk}\subset C$) is a divisor  of length $ n_{ij}-n_{i-1,j}$ (resp. $n_{lk}-n_{l,k-1}$).\\
$A_{C,\mathbf{n}_{\lambda}}$ comes equipped with universal divisors, which we denote by
\begin{align*}
    \CZ_{00}, \CX_{ij},\CY_{lk}\subset C\times A_{C,\mathbf{n}_{\lambda}},
\end{align*}
with fibers are
\begin{align*}
    \CZ_{00}|_{\underline{Z}}&=Z_{00},\\
    \CX_{ij}|_{\underline{Z}}&=X_{ij},\quad i\geq 1,\\
    \CY_{lk}|_{\underline{Z}}&=Y_{lk},\quad k\geq 1.
\end{align*}
For every $(i,j)\in \lambda$ with $i,j\geq 1$ define the  universal effective divisors
\begin{align*}
    \Gamma^1_{ij}&=\CX_{i,j}+\CY_{i-1,j},\\
    \Gamma^2_{ij}&=\CY_{i,j}+\CX_{i,j-1}.
\end{align*}
\begin{theorem}\label{thm: zero locus}
Let $C$ be an irreducible smooth quasi-projective curve, $\pi:C\times A_{C,\mathbf{n}_\lambda}\to A_{C,\mathbf{n}_\lambda}$ be the natural projection and define the vector bundle
 \begin{align*}
     \CE= \bigoplus_{\substack{(i,j)\in \lambda\\ i,j\geq 1}}\pi_*\oO_{\Gamma^2_{ij}  }(\Gamma^1_{ij}).
 \end{align*}
Then there exists a section $s$ of $\CE$ whose zero set is isomorphic to $C^{[ \mathbf{n}_\lambda]} $
\[
 \begin{tikzcd}
    &\CE  \arrow[d] \\
   C^{[ \mathbf{n}_\lambda]}\cong Z(s)\arrow[r, hook]
&A_{C,\mathbf{n}_\lambda}. \arrow[u, bend right, swap,  "s"]\end{tikzcd}
\]
 \end{theorem}
 \begin{proof}
Notice that $\CE$ is a vector bundle, as by cohomology and base change all higher direct images vanish 
\[\RR^k\pi_*\oO_{\Gamma^2_{ij}  }(\Gamma^1_{ij})=0, \quad k>0.\]
 There is a closed immersion 
\[C^{[\mathbf{n}_\lambda] }\hookrightarrow A_{C,\mathbf{n}_\lambda},\] 
given on closed points $(Z_{ij})_{(i,j)\in \lambda}\in C^{[\mathbf{n}_\lambda]}$ as 
 \[
(Z_{ij})_{(i,j)\in \lambda}\mapsto (Z_{00}, (Z_{ij}-Z_{i-1,j}),(Z_{ij}-Z_{i,j-1}))_{(i,j\in \lambda)}\in A_{C,\mathbf{n}_\lambda}.
\]
Define sections $\Tilde{s}_{ij}\in H^0(C\times A_{C,\mathbf{n}_\lambda}, \oO_{\Gamma^2_{ij}  }(\Gamma^1_{ij}) )$ as the composition 
 \begin{align*}
     \Tilde{s}_{ij}:  \oO_{C\times A_{C,\mathbf{n}_\lambda}}\xrightarrow{s'_{ij}}  \oO_{C\times A_{C,\mathbf{n}_\lambda}}(\Gamma^1_{ij})\to\oO_{\Gamma^2_{ij}}(\Gamma^1_{ij}),
 \end{align*}
 where $s'_{ij}$ is the section vanishing on $ \Gamma^1_{ij}$ while the second morphism is the restriction  along $j:\Gamma^2_{ij}\to  A_{C,\mathbf{n}_\lambda}\times C$; in other words, $\Tilde{s}_{ij}=j_*j^* s'_{ij} $.  The sections $\Tilde{s}_{ij}$ induce sections $s_{ij}=\pi_* \Tilde{s}_{ij}$ of $\pi_* \oO_{\Gamma^2_{ij}  }(\Gamma^1_{ij})$ and we set $s=(s_{ij})_{ij}\in H^0(A_{C,\mathbf{n}_\lambda},\CE)$. We claim that 
\begin{align*}
    C^{[\mathbf{n}_\lambda]}\cong Z(s).
\end{align*}
 To prove it, we follow the strategy of \cite[Prop. 2.4]{CK1}. For $(i,j)\in \lambda$ with $i,j\geq 1$, consider the universal divisors 
 \[
 \begin{tikzcd}
    &\Gamma^2_{ij}  \arrow[ld]\arrow[rd,"\Tilde{\pi}"]\arrow[r, hook, "j"]&A_{C, \mathbf{n}_\lambda}\times C\arrow[d, "\pi"] \\
C & & A_{C,\mathbf{n}_\lambda}. \end{tikzcd}
\]
Let  $\underline{Z}_T=(Z^T_{00}, X^T_{ij},Y^T_{lk})_{ij,lk}$  be any $T$-flat family with corresponding classifying morphism $f:T\to A_{C, \mathbf{n}_\lambda} $, where $Z^T_{00}, X^T_{ij},Y^T_{lk}\subset T\times C $ have zero-dimensional fibers of appropriate length. Consider the commutative diagram
 \[
 \begin{tikzcd}
    &\Gamma^{2,T}_{ij}  \arrow[ld]\arrow[rd,  "\Tilde{\pi}_T"]\arrow[r, hook, "j_T"]&T\times C \arrow[d, "\pi_T"]\\
 C& & T, \end{tikzcd}
\]
where $\Gamma^{1,T}_{ij},\Gamma^{2,T}_{ij}$ are the pullbacks of the universal divisors $ \Gamma^{1}_{ij},\Gamma^{2}_{ij}$ along $f\times \id_C$. To prove the claim it suffices to show that $\underline{Z}_T$ is a $T$-point of  $C^{[\mathbf{n}_\lambda]}$ if and only if $f$ factors through $Z(s)$.\\
Now  $\underline{Z}_T$ is a $T$-point of  $C^{[\mathbf{n}_\lambda]}$ if and only if $ \Gamma^{2,T}_{ij}  =\Gamma^{1,T}_{ij}$ for all $(i,j)\in \lambda $ such that $(i,j)\geq 1$. Notice that  the inclusion  $ \Gamma^{2,T}_{ij}  \subset\Gamma^{1,T}_{ij}$ is enough to have the equality, as all fibers are divisors in $C$ of the same degree.\\
On the other hand, we have that $f$ factors through $Z(s)$ if and only if $f^*s $ is the zero section of $f^*\CE$, i.e. if $f^*s_{ij}=0$ for all $(i,j)\in \lambda $ such that $(i,j)\geq 1$. Repeatedly applying flat base change, we obtain
\begin{align*}
    f^*s_{ij}=\Tilde{\pi}_* j_T^*(f\times \id_C)^* s'_{ij}.
\end{align*}
Therefore $ f^*s_{ij}$ is the zero section if and only if $ \Gamma^{2,T}_{ij}  \subset\Gamma^{1,T}_{ij}$ as required.
 \end{proof}
Thanks to Theorem \ref{thm: zero locus}, $C^{[\mathbf{n}_\lambda]}$ falls in the situation of Example \ref{exa:toy model pot} and we obtain a virtual fundamental class.
\begin{corollary}\label{cor: pot on double nested}
 Let $C$ be a smooth quasi-projective curve and $\mathbf{n}_\lambda$ a reversed plane partition. Then $C^{[\mathbf{n}_\lambda]}$ has a perfect obstruction theory
 \begin{align}\label{eqn: pot coming from zero locus double nested}
     [\CE^*|_{C^{[\mathbf{n}_\lambda]}}\to {\Omega^1_{A_{C,\mathbf{n}_\lambda}}}|_{C^{[\mathbf{n}_\lambda]}} ]\to \BL_{ C^{[\mathbf{n}_\lambda]}}.
 \end{align}
In particular there exists a virtual fundamental class
 \begin{align*}
     [C^{[\mathbf{n}_\lambda]}]^{\vir}\in A_{*}(C^{[\mathbf{n}_\lambda]} ).
     \end{align*}
 \end{corollary}
\subsection{Topological Euler characteristic}
Recall that we can view Euler characteristic weighted by a constructible function as a Lebesgue integral, where the measurable sets are constructible sets, measurable functions are constructible functions and the measure of a set is given by its Euler characteristic (cf. \cite[Sec. 2]{BK_elliptic_surface_DT}). In this language we have
\begin{align*}
    e(X)=\int_X 1\cdot de, 
\end{align*}
for any constructible set  $X$. The following lemma is reminiscent of  the existence of a power structure on the Grothendieck ring of varieties.
\begin{lemma}[{\cite[Lemma 32]{BK_elliptic_surface_DT}}]\label{lemma:Bryan_Kool}
Let $B$ be a scheme of finite type over $\BC$ and $e(B)$ its topological Euler characteristic. Let $g:\BZ_{\geq 0}\to \BZ(\!(p)\!)$ be any function with $g(0)=1$. Let $G:\Sym^n B\to \BZ(\!(p)\!) $ be the constructible function defined by 
\begin{align*}
    G(\mathbf{ax})=\prod_ig(a_i),
\end{align*}
for all $\mathbf{ax}=\sum_i a_ix_i\in \Sym^n B$, where $x_i\in B$ are distinct closed points. Then 
\begin{align*}
    \sum_{n=0}^{\infty}q^n \int_{\Sym^n B}G\cdot de=\left( \sum_{a=0}^\infty g(a)q^a\right)^{e(B)}.
\end{align*}
\end{lemma}
Using Lemma \ref{lemma:Bryan_Kool} we compute the topological Euler characteristic of double nested Hilbert schemes of points of any quasi-projective smooth curve.
\begin{theorem}\label{thm: euler char double nested}
Let $C$ be a smooth quasi-projective curve and $\lambda$ a Young diagram. Then 
\begin{align*}
    \sum_{\mathbf{n}_\lambda} e(C^{[\mathbf{n}_\lambda]})q^{|\mathbf{n}_\lambda|}=\prod_{\Box\in \lambda}(1-q^{h(\Box)})^{-e(C)}.
\end{align*}
\end{theorem}
\begin{proof}
Consider the constructible map 
\begin{align*}
    \rho_n:\bigsqcup_{|\mathbf{n}_\lambda|=n} C^{[\mathbf{n}_\lambda]}\to \Sym^{n} C,
\end{align*}
defined, for $\underline{Z}=(Z_{\Box})_{\Box\in \lambda}\in C^{[\mathbf{n}_\lambda]} $, by 
\begin{align*}
    \rho_n(\underline{Z})=\sum_{\Box\in \lambda}Z_{\Box}\in \Sym^{n} C.
\end{align*}
In other words, $\rho$ just forgets the distribution and the nesting of the divisor $\sum_{\Box\in \lambda}Z_{\Box}$ among all $\Box\in \lambda$.\\
Let $\mathbf{ax}=\sum_{i}a_ix_i\in \Sym^n C $, with $x_i$ different to each other. The fiber $\rho_n^{-1}(\mathbf{ax}) $ is clearly 0-dimensional and satisfies 
\begin{align}\label{eqn: 1}
    \rho_n^{-1}(\mathbf{ax})\cong \prod_i\rho_{a_i}^{-1}(a_ix_i).
\end{align}
In particular the Euler characteristic of the fiber $\rho_n(nx)$ does not depend on the point $x\in C$ and counts the number of reversed plane partition of size $n$ and underlying Young diagram $\lambda$
\begin{align}\label{eqn: 2}
    e(\rho_n^{-1}(nx))=\sum_{|\mathbf{n}_\lambda|=n}1.
\end{align}
Consider now
\begin{align*}
    \int_{\bigsqcup_{|\mathbf{n}_\lambda|=n} C^{[\mathbf{n}_\lambda]} }1\cdot de=\int_{\Sym^n C}{\rho_n}_*1\cdot de,
\end{align*}
where for any $ \mathbf{ax}\in \Sym^n C$ with $x_i$ different to each other, using \eqref{eqn: 1} and \eqref{eqn: 2}
\begin{align*}
    {\rho_n}_* 1(\mathbf{ax})&=e(\rho_n^{-1}(\mathbf{ax}))\\
    &= \prod_i   \sum_{|\mathbf{n}_\lambda|=a_i}1.
\end{align*}
Now,  $g(a)=\sum_{|\mathbf{n}_\lambda|=a}1$ and $G(\mathbf{ax} )={\rho_n}_* 1(\mathbf{ax})$  satisfy the hypotheses of Lemma \ref{lemma:Bryan_Kool} and therefore
\begin{align*}
     \sum_{n=0}^\infty\sum_{|\mathbf{n}_\lambda|=n} e(C^{[\mathbf{n}_\lambda]})q^n&= \sum_{n=0}^\infty q^n \int_{\bigsqcup_{|\mathbf{n}_\lambda|=n} C^{[\mathbf{n}_\lambda]} }1\cdot de\\
     &=  \sum_{n=0}^\infty q^n \int_{\Sym^n C}{\rho_n}_*1\cdot de\\
     &= \left(\sum_{n=0}^\infty \sum_{|\mathbf{n}_{\lambda}|=n}q^{|\mathbf{n}_\lambda|}  \right)^{e(C)}.
\end{align*}
A closed formula for the generating series of reversed plane partitions was given by Stanley \cite[Prop. 18.3]{Stanley_theory_application_I_II} and by Hillman-Grassl \cite[Thm. 1]{HG_reversed_plane_partitions} using hook-lengths
\begin{align*}
    \sum_{n=0}^\infty \sum_{|\mathbf{n}_{\lambda}|=n}q^{|\mathbf{n}_\lambda|} =\prod_{\Box\in \lambda}(1-q^{h(\Box)})^{-1},
\end{align*}
by which we conclude the proof.
\end{proof}

\subsection{Double nesting of divisors}
We conclude this section with a generalization of the zero-locus construction of Theorem \ref{thm: zero locus}.\\

Let $X$ be a smooth projective variety of dimension $d$ and $\bm{\beta}_\lambda=(\beta_\Box)_{\Box\in \lambda}$ be a collection of homology classes $\beta_\Box\in H_{n-2}(X,\BZ)$. Denote by $H_{\bm{\beta}_\lambda}$ the double nested Hilbert scheme of effective divisors on $X$, which parametrizes flags of divisors $(Z_{\Box})_{\Box\in \lambda}\subset X$ satisfying the nesting condition dictated by $\bm{\beta}_\lambda$. Denote by
\begin{align*}
A_{X,\bm{\beta}_{\lambda}}=H_{\beta_{00}}\times \prod_{\substack{(i,j)\in \lambda\\ i\geq 1}}H_{\beta_{ij}-\beta_{i-1,j}}\times \prod_{\substack{(l,k)\in \lambda\\ k\geq 1}}H_{\beta_{lk}-\beta_{l,k-1}},
\end{align*}
where $H_\beta$ is the usual Hilbert scheme of divisors on $X$ of class $\beta$. Analogously to Section \ref{sec: pot on points on curves}, $A_{X,\bm{\beta}_{\lambda}}$ comes equipped with universal (Cartier) divisors $   \CZ_{00}, \CX_{ij},\CY_{lk}\subset X\times A_{X,\bm{\beta}_{\lambda}}$ and for every $(i,j)\in \lambda$ with $i,j\geq 1$ we define the  universal effective divisors
\begin{align*}
    \Gamma^1_{ij}&=\CX_{i,j}+\CY_{i-1,j},\\
    \Gamma^2_{ij}&=\CY_{i,j}+\CX_{i,j-1}.
\end{align*}
Define the coherent sheaf 
 \begin{align*}
     \CE= \bigoplus_{\substack{(i,j)\in \lambda\\ i,j\geq 1}}\pi_*\oO_{\Gamma^2_{ij}  }(\Gamma^1_{ij}),
 \end{align*}
 where $\pi:X\times A_{X,\bm{\beta}_{\lambda}}\to A_{X,\bm{\beta}_{\lambda}}$ is the natural projection. Under some extra assumptions on $X$ and $ \bm{\beta}_{\lambda}$ Theorem \ref{thm: zero locus} generalizes.
\begin{prop}
Assume that $ A_{X,\bm{\beta}_{\lambda}}$ is smooth and $\CE$ is a vector bundle. Then there exists a section $s$ of $\CE$ such that 
\[
 \begin{tikzcd}
    &\CE  \arrow[d] \\
  H_{\bm{\beta}_\lambda}\cong Z(s)\arrow[r, hook]
&A_{X,\bm{\beta}_{\lambda}}. \arrow[u, bend right, swap,  "s"]\end{tikzcd}
\]
In particular, $  H_{\bm{\beta}_\lambda}$ has a perfect obstruction theory.
 \end{prop}
 \begin{corollary}
Let $X=\BP^{n_1}\times \cdots\times  \BP^{n_s}$. Then there exists a virtual fundamental class  $  [H_{\bm{\beta}_\lambda}]^{\vir}\in A_*( H_{\bm{\beta}_\lambda})$.
 \end{corollary}
 \begin{proof}
 The smoothness of $A_{X,\bm{\beta}_{\lambda}} $ follows by the smoothness of  $H_{\beta}\cong \BP(H^0(X,\oO_X( \beta)))$ for every $\beta\in H_{n-2}(X,\BZ)$. Let $D_1, D_2$ be two effective divisors on $X$ such that $[D_1]=[D_2]\in H_{n-2}(X,\BZ) $; in particular, $D_1,D_2$ are linearly equivalent. Combining the long exact sequence in cohomology of the short exact sequence
\begin{align*}
    0\to \oO(D_1-D_2)\to \oO(D_1)\to O_{D_2}(D_1)\to 0
\end{align*}
and the vanishings
\begin{align*}
    &H^k(X, \oO(D_1))=0, \quad k=1, \dots, \dim X -1,\\
    &H^k(X, \oO(D_1-D_2))=0, \quad k=2,\dots \dim X,
\end{align*}
 yields  that $H^k(X, \oO_{D_2 }(D_1))=0 $ for $k\geq 1$; cohomology and base  implies that $\RR^k\pi_*\oO_{\Gamma^2_{ij}  }(\Gamma^1_{ij})=0$ for $k\geq 1$. Finally, We have that, if $\dim X\geq 2$,
 \begin{align*}
     \dim H^0(X, \oO_{D_2}(D_1))&=\dim
H^0(X, \oO_{X}(D_1))-1,
\end{align*}
 which depends only on the degree $[D_1]=[D_2]\in H_2(X, \BZ)$ and
 implies that the dimension of the fibers of  $\CE$ is constant, thus $\CE$ is a vector bundle.
 \end{proof}
\section{Moduli space of stable pairs}
\subsection{Moduli space of stable pairs}
Moduli spaces of stable pairs were introduced by Pandharipande-Thomas \cite{PT_curve_counting_derived} in order to give a geometric interpretation of the MNOP conjectures \cite{MNOP_1}, through the DT/PT correspondence proved by Toda (for Euler characteristic) and Bridgeland in \cite{Toda_curve_counting_DT_PT, Bri_Hall_algebras_curve_counting} using wall-crossing and Hall algebra techniques.\\
For a smooth quasi-projective threefold $X$, a curve class $\beta\in H_2(X, \BZ)$ and $n\in \BZ$, we define $P_{n}(X,\beta)$ to be the moduli space of pairs
\[I^\bullet=[\oO_X\xrightarrow{s} F]\in  \derived^{\mathrm{b}}(X)\]
in the derived category of $X$ where $F $ is a pure 1-dimensional sheaf with proper support $[\mathrm{supp}( F)]=\beta$ with $\chi(F)=n$ and $s$ is a section with 0-dimensional cokernel.\\
By the work of   Huybrechts-Thomas \cite{HT_obstruction_theory},  the Atiyah class gives a \emph{perfect obstruction theory} on $P_n(X,\beta)$
\begin{equation}\label{eqn: obstruction theory}
     \BE=\RR\hom_\pi(\BI, \BI)^\vee_0[-1]\to \BL_{P_n(X,\beta)},
\end{equation}
where 
$(\cdot)_0$ denotes the trace-free part,
 $\pi:X\times P_n(X,\beta)\to P_n(X,\beta)$ is the canonical projection and $\BI^\bullet=[\oO\to \BF]$ is the universal stable pair on $X\times P_n(X,\beta) $.\\
If $X$ is projective, the perfect obstruction theory induces a virtual fundamental class $[ P_n(X,\beta)]^{\vir}\in A_*(P_n(X,\beta)) $ and one defines \emph{stable pair} (or \emph{PT}) \emph{invariants} by integrating cohomology classes $\gamma\in H^*(P_{n}(X,\beta), \BZ)$ against the virtual fundamental class
\begin{align}\label{eqn: PT invariants compact}
    \PT_{\beta,n}(X,\gamma)=\int_{[ P_n(X,\beta)]^{\vir}}\gamma\in \BZ.
\end{align}
We focus here in the case of $X$ a \emph{local curve}, i.e.  $X=\Tot_C(L_1\oplus L_2)$ the total space of the direct sum of two line bundles $L_1, L_2$ on a smooth  projective curve $C$ and $\beta=d[C]\in H_2(X,\BZ)\cong H_2(C,\BZ)$  a multiple of the zero section of $X\to C $.\\
$X$ is a smooth \emph{quasi-projective} threefold, therefore the moduli space of stable pairs $ P_n(X,\beta)$ is hardly ever a proper scheme and one cannot define invariants as in \eqref{eqn: PT invariants compact}. Nevertheless, the algebraic torus $\TT=(\BC^*)^2$ acts on $X$ by scaling the fibers and the action naturally lifts to $P_n(X,d[C]) $, making the perfect obstruction theory naturally $\TT$-equivariant by \cite[Example 4.6]{Ric_equivariant_Atiyah}. Moreover, the $\TT$-fixed locus $P_n(X,d[C])^\TT$ is proper (cf. Prop. \ref{prop: iso of schemes fixed locus}), therefore by Graber-Pandharipande \cite{GP_virtual_localization} there is naturally an induced perfect obstruction theory on  $P_n(X,d[C])^\TT$ and a virtual fundamental class $[P_n(X,d[C])^\TT]^{\vir}\in A_*(P_n(X,d[C])^\TT) $. We define \emph{$\TT$-equivariant stable pair invariants} by  Graber-Pandharipande virtual localization formula
\begin{align}\label{eqn: localized PT invariants}
    \PT_{d,n}(X)=\int_{[ P_n(X,d[C])^{\TT}]^{\vir}}\frac{1}{e^{\TT}(N^{\vir})}\in \BQ(s_1,s_2),
\end{align}
where $s_1, s_2$ are the generators of $\TT$-equivariant cohomology and the virtual normal bundle
\begin{align}\label{eqn: virtual normal bundle}
    N^{\vir}=(\BE|^\vee_{ P_n(X,d[C])^{\TT}})^{\mov}\in K^0_\TT(P_n(X,d[C])^{\TT})
\end{align}
is the $\TT$-moving part of the restriction of the dual of the perfect obstruction theory. Stable pair invariants with descendent insertions on local curves have been studied in \cite{PP:local_curves_stationary, PP_descendants_local_curves, Oblomkov_EGL_type}.
\subsection{Torus representations and their weights}
Let $\TT = (\BC^*)^g$ be an algebraic torus, with character lattice $\widehat{\TT} = \Hom(\TT,\BC^\ast) \cong \BZ^g$. Let $K^0_{\TT}(\pt)$ be the $K$-group of the category of $\TT$-representations. Any finite dimensional $\TT$-representation $V$ splits as a sum of $1$-dimensional representations called the \emph{weights} of $V$. Each weight corresponds to a character $\mu \in \widehat{\TT}$, and in turn each character corresponds to a monomial $\tf^\mu = \tf_1^{\mu_1}\cdots \tf_g^{\mu_g}$ in the coordinates of $\TT$. The map
\begin{equation*}
\tr\colon K^0_{\TT}(\pt) \to \BZ \left[\tf^\mu \mid \mu \in \widehat{\TT}\right],\quad V\mapsto \tr_V,
\end{equation*}
sending the class of a $\TT$-representation to its decomposition into weight spaces is a ring isomorphism, where tensor product on the left corresponds to the natural multiplication on the right. We will therefore sometimes identify a (virtual) $\TT$-representation with its character. \\
If $X$ is a scheme with a trivial $\TT$-action, every $\TT$-equivariant coherent sheaf on $X$ decomposes as $F=\bigoplus_{\mu\in \widehat{\TT}}F_\mu\otimes \tf^\mu $.
\subsection{The fixed locus}
In this section we prove that the $\TT$-fixed locus $P_n(X,d[C])^\TT$ is a disjoint union of double nested Hilbert schemes of points  $C^{[\mathbf{n}_\lambda]}$, for suitable reversed plane partitions $\mathbf{n}_\lambda$, where $\lambda $ are  Young diagram of size $|\lambda|=d$. Our strategy is similar to Kool-Thomas \cite[Sec. 4]{KT_vertical} for local surfaces.\\
Given a $\TT$-equivariant coherent sheaf on $X$, its pushdown along $p:X\to C$ decomposes into weight spaces (e.g. by \cite[Ex.  II.5.17, II.5.18]{Hart_Algebraic_Geometry})
\begin{align*}
    p_*F=\bigoplus_{(i,j)\in \BZ^2}F_{ij}\otimes \tf_1^{i}\tf_2^{j},
\end{align*}
where $F_{ij}$ is a coherent sheaf on $C$. For example 
\begin{align*}
    p_*\oO_X=\bigoplus_{i,j\geq 0}L_1^{-i}\otimes L_2^{-j}\otimes \tf_1^{-i}\tf_2^{-j}.
\end{align*}
Since $p$ is affine, the pushdown does not lose any information, and we recover the $\oO_X$-module structure of $F$ by the $   p_*\oO_X$-action that $ p_*F$ carries. This is generated by the action of the $-1$ pieces $L_1^{-1}\otimes\tf_1^{-1},L_2^{-1}\otimes\tf_2^{-1}$, so we find that the $\oO_X$-module structure is determined by the maps
\begin{align}\label{eqn: eqn maps module structure}
\begin{split}
      &\left(\bigoplus_{i,j}F_{ij}\otimes \tf_1^{i}\tf_2^{j} \right)\otimes L_1^{-1}\otimes \tf_1^{-1}\to \bigoplus_{i,j}F_{ij}\otimes \tf_1^{i}\tf_2^{j},\\
    &\left(\bigoplus_{i,j}F_{ij}\otimes \tf_1^{i}\tf_2^{j} \right)\otimes L_2^{-1}\otimes \tf_2^{-1}\to \bigoplus_{i,j}F_{ij}\otimes \tf_1^{i}\tf_2^{j},
\end{split}
\end{align}
which commute with both the actions of $\oO_C$ and $\TT$. In other words, \eqref{eqn: eqn maps module structure} are $\TT$-equivariant maps of $\oO_C$-modules. By $\TT$-equivariance, they are sums of maps
\begin{align}\label{eqn: maps weight -1}
\begin{split}
     &F_{ij}\otimes L_1^{-1}\to F_{i-1,j},\\
 &F_{ij}\otimes L_2^{-1}\to F_{i,j-1}. 
\end{split}
\end{align}
Let now $(F,s)\in P_n(X, d[C])^\TT$ be a $\TT$-fixed stable pair. Then  $s$ is a $\TT$-equivariant section of a  $\TT$-equivariant coherent sheaf  $F$  on $X$. Applying $p_*$ to $\oO_X\xrightarrow{s} F $ gives a graded map
which commutes with the maps \eqref{eqn: maps weight -1}.
 Writing 
\begin{align*}
    G_{ij}=F_{-i,-j}\otimes L_1^{i}\otimes L_2^{j},
\end{align*}
we find that the $\TT$-fixed stable pair $(F,s)$ on $X$ is equivalent to the following data of sheaves and commuting maps on $C$
\begin{equation}\label{eqn: diagram with G_i}
\begin{tikzcd}
& &&& \oO_C\arrow[rr, equal]\arrow[ld, equal]\arrow[ddd]& & \oO_C\arrow[rr, equal]\arrow[ld, equal]\arrow[ddd]& & \oO_C\arrow[r, equal] \arrow[ld, equal]\arrow[ddd]&\dots\\
 &  &&\oO_C\arrow[rr, crossing over, equal]\arrow[ld, equal]& &  \oO_C\arrow[rr, crossing over, equal]\arrow[ld, equal]& & \oO_C\arrow[r, crossing over, equal]\arrow[ld]&\dots&\\
  &  &  \dots& &  \dots& & \dots&  & &\\
  &\dots \arrow[r] & G_{0,-1}\arrow[rr]\arrow[ld] && G_{00}\arrow[rr]\arrow[ld]& & G_{01}\arrow[rr]\arrow[ld]& & G_{02}\arrow[r] \arrow[ld]&\dots\\
 \dots \arrow[r] & G_{1,-1}\arrow[rr]&&G_{10}\arrow[rr]\arrow[ld]\arrow[from=uuu, crossing over]& &  G_{11}\arrow[rr]\arrow[ld]\arrow[from=uuu, crossing over]& & G_{12}\arrow[r]\arrow[ld]\arrow[from=uuu, crossing over]&\dots&\\
  &  &  \dots& &  \dots& & \dots& & &\\
\end{tikzcd} 
\end{equation}
By the  purity of $F$, each $G_{ij}$ is either zero or a pure 1-dimensional coherent sheaf on $C$, and the "vertical" maps are generically isomorphisms. In particular, for every $(i,j)$ such that either $i<0$ or $j<0$, it follows that   $G_{ij}$ is zero-dimensional and therefore vanishes by the purity assumption. Moreover, if $G_{ij}$ is non-zero, it is a rank 1 torsion-free sheaf on a smooth curve, that is a line bundle on $C$ (with a section). Finally, any $\TT$-equivariant stable pair on $X$ is set-theoretically supported on $C$, thus is properly supported on $X$ and only finitely many $G_{ij}$ can be non-zero. This results in a diagram of the following shape
\begin{equation}
    \begin{tikzcd}\label{eqn: finitely many line bundles with divisors}
 && \oO_C\arrow[rr, equal]\arrow[ld, equal]\arrow[ddd]& & \oO_C\arrow[rr, equal]\arrow[ld, equal]\arrow[ddd]& & \oO_C\arrow[r, equal] \arrow[ld, equal]\arrow[ddd]&\dots\\
   &\oO_C\arrow[rr, crossing over, equal]\arrow[ld, equal]& &  \oO_C\arrow[rr, crossing over, equal]\arrow[ld,equal]& & \oO_C\arrow[r, crossing over, equal]\arrow[ld, equal]&\dots&\\
      \dots& &  \dots& & \dots&  & &\\
    && \oO_C(Z_{00})\arrow[rr, hook]\arrow[ld, hook]& & \oO_C(Z_{01})\arrow[rr, hook]\arrow[ld, hook]& & \oO_C(Z_{02})\arrow[r, hook] \arrow[ld, hook]&\dots\\
   &\oO_C(Z_{10})\arrow[rr, hook]\arrow[ld, hook]\arrow[from=uuu, crossing over]& &  \oO_C(Z_{11})\arrow[rr, hook]\arrow[ld, hook]\arrow[from=uuu, crossing over]& &\oO_C(Z_{12})\arrow[r, hook]\arrow[ld, hook]\arrow[from=uuu, crossing over]&\dots&\\
      \dots& &  \dots& & \dots& & &\\
\end{tikzcd}
\end{equation}
where $Z_{ij}$ are divisors on $C$ and all "horizontal" maps are injections of line bundles. Therefore, a $\TT$-fixed stable pair $(F,s)$ is equivalent to a nesting of divisors 
\begin{equation}\label{eqn: nesting divisors}
  \begin{tikzcd}
    Z_{00}\arrow[r, phantom, "\subset"]\arrow[d, phantom, "\cap"] &Z_{01}\arrow[d, phantom, "\cap"]\arrow[r, phantom, "\subset"]&Z_{02}\arrow[d, phantom, "\cap"]\arrow[r, phantom, "\subset"]\arrow[d, phantom, "\cap"]&Z_{03}\arrow[r, phantom, "\subset"]\arrow[d, phantom, "\cap"]&\dots\\
      Z_{10}\arrow[r, phantom, "\subset"]\arrow[d, phantom, "\cap"]&Z_{11}\arrow[d, phantom, "\cap"]\arrow[r, phantom, "\subset"]&Z_{12}\arrow[d, phantom, "\cap"]\arrow[r, phantom, "\subset"]&Z_{13}\arrow[d, phantom, "\cap"]\arrow[r, phantom, "\subset"]&\dots\\
        Z_{20}\arrow[d, phantom, "\cap"]\arrow[d, phantom, "\cap"]\arrow[r, phantom, "\subset"]&Z_{21}\arrow[d, phantom, "\cap"]\arrow[r, phantom, "\subset"]&Z_{22}\arrow[d, phantom, "\cap"]\arrow[r, phantom, "\subset"]&\dots &\\
          \dots&\dots&\dots& &\\
  \end{tikzcd}  
\end{equation}
where the nesting is dictated by a Young diagram $\lambda$. This results into a point of the double nested Hilbert scheme $C^{[\mathbf{n}_\lambda]}$, where by Riemann-Roch
\begin{align*}
    \chi(F)=|\mathbf{n}_\lambda|+ \mathbf{f}_{\lambda,g}( \deg L_1, \deg L_2),
\end{align*}
where for   a Young diagram $\lambda$  and $g,k_1, k_2\in \BZ$ we  define
\begin{align}\label{eqn: f lambda g}
\begin{split}
    \mathbf{f}_{\lambda,g}(k_1,k_2)&=\sum_{(i,j)\in \lambda}(1-g-i\cdot k_1-j\cdot k_2 )\\
    &=|\lambda|(1-g)-k_1\cdot n(\lambda)-k_2\cdot n(\overline{\lambda}).
\end{split}
\end{align}
 Conversely, any  nesting of divisors as in \eqref{eqn: nesting divisors} corresponds to a diagram of sheaves as in \eqref{eqn: finitely many line bundles with divisors}, which corresponds to a $\TT$-fixed stable pair on $X$. Therefore we have a bijection of set:
 \begin{align}\label{eqn: bijection of sets}
     P_n(X,d[C])^\TT= \bigsqcup_{\lambda\vdash d}\bigsqcup_{\mathbf{n}_\lambda} C^{[\mathbf{n}_\lambda]},
 \end{align}
 where the disjoint union is over all Young diagrams $\lambda$ of size $d$ and all reversed plane partitions $\mathbf{n}_\lambda $ satisfying   $n=|\mathbf{n}_\lambda|+ \mathbf{f}_{\lambda,g}( \deg L_1, \deg L_2)$. We mimic \cite[Prop. 4.1]{KT_vertical} to prove that the above bijection on sets is an isomorphism of schemes.
 \begin{prop}\label{prop: iso of schemes fixed locus}
There exists an isomorphism of schemes
\begin{align*}
    P_n(X,d[C])^\TT= \bigsqcup_{\lambda\vdash d}\bigsqcup_{\mathbf{n}_{\lambda}} C^{[\mathbf{n}_\lambda]},
\end{align*}
where the disjoint union is over all Young diagrams $\lambda$ of size $d$ and all reversed plane partitions $\mathbf{n}_\lambda $ satisfying
\begin{align*}
    n=|\mathbf{n}_\lambda|+ \mathbf{f}_{\lambda,g}( \deg L_1, \deg L_2).
\end{align*}
In particular, $ P_n(X,d[C])^\TT$ is proper.
 \end{prop}
 \begin{proof}
Let $B$ be any (connected) scheme over $\BC$. We  need  to adapt the construction of this section to a $\TT$-fixed stable pair on $X\times B$, flat over $B$. Pushing down by the affine map $p:X\times B\to C\times B$ gives a graded sheaf $\bigoplus_{i,j}F_{ij}$ on $C\times B$, flat over $B$ (therefore so are all its weight spaces $F_{ij}$). The original sheaf $F$ on $X\times B$ can be reconstructed from the maps \eqref{eqn: maps weight -1}. Therefore a $\TT$-fixed  pair $(F,s)$ on $X\times B$, flat over $B$, is equivalent to the data \eqref{eqn: diagram with G_i}, with each $G_{ij}$ on $C\times B$,  flat over B.\\
If $(F,s)$ is a \emph{stable pair}, over each closed fiber $C\times \{b\}$, where $b\in B$, we showed that each (non-zero) $G_{ij}$ is a line bundle. By \cite[Lemma 2.1.7]{HL_moduli_space}, this shows that each (non-zero) $G_{ij}$ is a line bundle on $C\times B$. Together with its non-zero section, this defines divisors $Z_{ij}\subset C\times B$, flat over $B$, satisfying the nesting condition of \eqref{eqn: nesting divisors}, which yields a $B$-point $B\to \bigsqcup_{\mathbf{n}_\lambda} C^{[\mathbf{n}_\lambda]}$. Conversely, any $B$-point $B\to \bigsqcup_{\mathbf{n}_\lambda} C^{[\mathbf{n}_\lambda]}$ defines a diagram \eqref{eqn: finitely many line bundles with divisors}, equivalent to a $\TT$-fixed stable pair $(F,s)$ on $X\times B$, flat over $B$.
 \end{proof}
 As a corollary, we compute the generating series of the topological Euler characteristic of the moduli space of stable pairs on a local curve.
 \begin{corollary}
 Let $C$ be a smooth projective curve of genus $g$, $L_1,L_2$ line bundle on $C$  and set $X=\Tot_C(L_1\oplus L_2)$. Then for any $d> 0$ we have
 \begin{align*}
     \sum_{n\in \BZ}e(P_n(X,d[C]))\cdot q^n= \sum_{|\lambda|=d}q^{|\lambda|(1-g)-\deg L_1 n(\lambda)-\deg L_2 n(\overline{\lambda})}\prod_{\Box\in \lambda}(1-q^{h(\Box)})^{2g-2}.
 \end{align*}
 \end{corollary}
 \begin{proof}
 The topological Euler characteristic of a $\TT$-scheme is the same of its $\TT$-fixed locus, therefore
 \begin{align*}
     \sum_{n\in \BZ}e(P_n(X,d[C]))\cdot q^n&=\sum_{n\in \BZ}e\left(P_n(X,d[C])^{\TT}\right)\cdot q^n\\
     &= \sum_{n\in \BZ}\sum_{|\lambda|=d}\sum_{\substack{|\mathbf{n}_\lambda|=n-\\\mathbf{f}_{\lambda,g}( \deg L_1, \deg L_2)}}q^n\cdot e(C^{[\mathbf{n}_\lambda]})\\
     &= \sum_{|\lambda|=d}q^{\mathbf{f}_{\lambda,g}( \deg L_1, \deg L_2)} \sum_{\mathbf{n}_\lambda} q^{|\mathbf{n}_\lambda|}\cdot e(C^{[\mathbf{n}_\lambda]})\\
     &= \sum_{|\lambda|=d}q^{\mathbf{f}_{\lambda,g}( \deg L_1, \deg L_2)}\prod_{\Box\in \lambda}(1-q^{h(\Box)})^{2g-2},
 \end{align*}
 where in the second line we applied Proposition \ref{prop: iso of schemes fixed locus} and in the last line Theorem \ref{thm: euler char double nested}.
 \end{proof}
 \section{K-theory class of  the perfect obstruction theory}
 Let $C$ be an irreducible smooth projective curve. On each connected component $C^{[\mathbf{n}_\lambda]}\subset P_{n}(X,\beta)^{\TT} $ of the $\TT$-fixed locus there exists and induced  virtual fundamental class $[C^{[\mathbf{n}_\lambda]}]_{\mathrm{PT}}^{\vir}$ coming from the perfect obstruction theory \eqref{eqn: obstruction theory}. We show in this section that $[C^{[\mathbf{n}_\lambda]}]_{\mathrm{PT}}^{\vir}$ agrees with the virtual fundamental class constructed in Corollary \ref{cor: pot on double nested} and compute the class in $K$-theory of the virtual normal bundle $N^{\vir}$.\\

We start by describing the class in $K$-theory of (the restriction of) the perfect obstruction theory $\BE\in K^{\TT}_0(C^{[\mathbf{n}_\lambda]})$. To ease  readability we will omit various pullbacks  whenever they are clear from the context. Recall the following identities in $K$-theory
 \begin{align*}
  &\BE=-\RR\pi_*\RR\hom(\BI, \BI )^{\vee}_0\in K_0^\TT(P_{n}(X,\beta)^\TT),\\
     &\BF=\sum_{(i,j)\in \lambda}i_*\oO_{C\times C^{[\mathbf{n}_\lambda]}}(\CZ_{ij}) \otimes L_1^{-i}\otimes L_2^{-j}\otimes \tf_1^{-i}\tf_2^{-j} \in K_0^\TT(X\times C^{[\mathbf{n}_\lambda]} ),\\
     & \BI=\oO_{X\times P_{n}(X,\beta)^{\TT}}-\BF\in K_0^\TT(X\times P_{n}(X,\beta)^\TT),
 \end{align*}
 where the various maps fit in the diagram
 \[
 \begin{tikzcd}
   C\arrow[r,"i", bend left] &X \arrow[l, "p"]  & X\times C^{[\mathbf{n}_\lambda]}\arrow[l]\arrow[r, "\pi"] \arrow[d, "p"]& C^{[\mathbf{n}_\lambda]} \\
  & & C\times C^{[ \mathbf{n}_\lambda]} \arrow[llu]\arrow[u, "i", bend left]\arrow[ur, "\pi "]& 
\end{tikzcd}
\]
 Again, to ease  notation, we keep denoting   $i\times \id_{C^{[\mathbf{n}_\lambda]}},p\times \id_{C^{[\mathbf{n}_\lambda]}} $ by $i,p$ and  by   $\pi$ the composition $\pi\circ i$. We compute
  \begin{align*}
     \BE^\vee&=\RR\pi_*\RR\hom( \oO,\BF)+\RR\pi_*\RR\hom(\BF, \oO)-\RR\pi_*\RR\hom( \BF,\BF)\\
     &= \RR\pi_*\BF-\left(\RR\pi_*(\BF\otimes K_X)\right)^\vee\otimes \tf_1\tf_2 -\RR\pi_*\RR\hom( \BF,\BF),
 \end{align*}
where in the second equality we applied Grothendieck duality\footnote{Even though $X$ is not proper, we can pass to a compactification $X\hookrightarrow \overline{X}$ and use Grothendieck duality exploiting the fact that the sheaves involved have  proper support, see \cite[footnote 4]{CK1}.}  and the projection formula and $K_X=K_C\otimes L_1^{-1}\otimes L_2^{-1}$. Define $\Lambda^\bullet(V)=\sum_{i=0}^{\rk V}(-1)^{i}\Lambda^i V$ for a locally free sheaf $V$ and extend it by linearity to any class in  $K^0_\TT(C)$. By \cite[Lemma 5.4.9]{CG_representation_theory} 
, for every $\TT$-equivariant coherent sheaf $\CF\in K^{\TT}_0(C)$, we have
 \begin{align*}
   \mathbf{L}i^*i_*\CF=\Lambda^\bullet N_{C/X}^*\otimes \CF\in K^{\TT}_0(C),  
 \end{align*}
where $ N_{C/X}=L_1\otimes \mathfrak{t}_1\oplus L_2\otimes \mathfrak{t}_2$ is the $\TT$-equivariant normal bundle of $i:C\to X$; an analogous statement holds for $i:C\times C^{[\mathbf{n}_\lambda]}\to X\times C^{[\mathbf{n}_\lambda]}$. Therefore
\begin{align*}
     \RR\pi_*\RR\hom( \BF,\BF)&=\sum_{(i,j),(l,k)\in \lambda}\RR\pi_*\RR\hom_X( i_*\oO_{C\times C^{[\mathbf{n}_\lambda]}}(\CZ_{ij}),i_*\oO_{C\times C^{[\mathbf{n}_\lambda]}}(\CZ_{lk})\otimes L_1^{i-l}L_2^{j-k} )\tf_1^{i-l}\tf_2^{j-k}\\
     &=\sum_{(i,j),(l,k)\in \lambda}\RR\pi_*\RR\hom_C( \mathbf{L}i^*i_*\oO_{C\times C^{[\mathbf{n}_\lambda]}}(\CZ_{ij}),\oO_{C\times C^{[\mathbf{n}_\lambda]}}(\CZ_{lk})\otimes L_1^{i-l}L_2^{j-k} )\tf_1^{i-l}\tf_2^{j-k}\\
     &=\sum_{(i,j),(l,k)\in \lambda}\RR\pi_*\left((\oO-L_1\tf_1-L_2\tf_2+ L_1 L_2\tf_1\tf_2)(\CZ_{lk}-\CZ_{ij})\otimes L_1^{i-l}L_2^{j-k} )\tf_1^{i-l}\tf_2^{j-k} \right),
 \end{align*}
where in the second line we used adjunction in the derived category. To simplify the notation, for any $(i,j), (l,k)\in \lambda$ set $\Delta_{ij;lk}=\CZ_{lk}-\CZ_{ij} $, which is an effective divisor if $(i,j)\leq (l,k) $. Putting all together we have the following identity in $K_0^\TT(C^{[\mathbf{n}_\lambda]})$
 \begin{multline}\label{eqn: dual pot}
     \BE^\vee=\sum_{(i,j)\in \lambda}\RR\pi_*(\oO_{C\times C^{[\mathbf{n}_\lambda]}}(\CZ_{ij}) \otimes L_1^{-i} L_2^{-j}) \tf_1^{-i}\tf_2^{-j}\\
     -\sum_{(i,j)\in \lambda}\left(\RR\pi_*(\oO_{C\times C^{[\mathbf{n}_\lambda]}}(\CZ_{ij}) \otimes K_X\otimes L_1^{-i} L_2^{-j} \tf_1^{-i}\tf_2^{-j})\right)^\vee \otimes \tf_1\tf_2\\
     -\sum_{(i,j),(l,k)\in \lambda}\RR\pi_*\left((\oO-L_1\tf_1-L_2\tf_2+ L_1L_2\tf_1\tf_2)(\Delta_{ij;lk})\otimes L_1^{i-l}L_2^{j-k} )\tf_1^{i-l}\tf_2^{j-k} \right).
 \end{multline}
 \begin{theorem}\label{thm: equality virtual classes}
 There is an identity of virtual fundamental classes
 \begin{align*}
     [C^{[\mathbf{n}_\lambda]}]^{\vir}_{\mathrm{PT}}=[C^{[\mathbf{n}_\lambda]}]^{\vir}\in A_*(C^{[\mathbf{n}_\lambda]}),
 \end{align*}
 where  the class on the left-hand-side is induced by \eqref{eqn: obstruction theory} by Graber-Pandharipande localization and the one on the right-hand-side is constructed in Corollary \ref{cor: pot on double nested}. 
 \end{theorem}
 \begin{proof}
By a result of Siebert \cite[Thm. 4.6]{Sie_virtual_fund_classes} any two virtual fundamental classes coincide if the classes in $K$-theory of their perfect obstruction theory agree. The class in $K$-theory of the dual of the induced PT perfect obstruction theory is the $\TT$-fixed part of $\BE^\vee $ by \cite[Prop. 1]{GP_virtual_localization}
\begin{multline*}
    (\BE^\vee)^{\fix}=\RR\pi_*\oO_{C\times C^{[\mathbf{n}_\lambda]}}(\CZ_{00})-\sum_{(i,j)\in \lambda}\RR\pi_*\oO_{C\times C^{[\mathbf{n}_\lambda]}}+\sum_{\substack{(i,j)\in \lambda\\ i\geq 1}}\RR\pi_*\oO_{C\times C^{[\mathbf{n}_\lambda]}}(\Delta_{i-1,j;ij})\\
    +\sum_{\substack{(i,j)\in \lambda\\ j\geq 1}}\RR\pi_*\oO_{C\times C^{[\mathbf{n}_\lambda]}}(\Delta_{i,j-1;ij}) -\sum_{\substack{(i,j)\in \lambda\\ i,j\geq 1}}\RR\pi_*\oO_{C\times C^{[\mathbf{n}_\lambda]}}(\Delta_{i-1,j-1;ij}).
\end{multline*}
We explained in Section \ref{sec: notation partitions} how to associate a graph to any Young diagram $\lambda$. Notice
 that boxes $(i,j)\in \lambda$ are in bijection with the vertices $V$,  boxes $(i,j)\in \lambda$, such that $ i\geq 1 $ (resp. $j\geq 1$) are in bijection with vertical (resp. horizontal) edges $E$ and boxes $(i,j)\in \lambda$, such that $ i,j\geq 1 $ are in bijection with squares $Q$ of the associated graph. By Lemma \ref{lemma: vertices and edges of graph}
\begin{align}\label{eqn:betti numbers of graph}
    V-E+Q-1=0. 
\end{align}
Combining this identity with the universal exact sequences
\begin{align}\label{eqn: universal exact sequences}
\begin{split}
     &0\to \oO\to \oO(\CZ_{00} )\to \oO_{\CZ_{00}}(\CZ_{00} )\to 0,\\
     & 0\to \oO\to \oO(\Delta_{ij;lk} )\to \oO_{\Delta_{ij;lk}}(\Delta_{ij;lk} )\to 0,
\end{split}
\end{align}
whenever $ \Delta_{ij,lk}$ is an effective divisor, one gets
 \begin{multline*}
    (\BE^\vee)^{\fix}=\RR\pi_*\oO_{\CZ_{00}}(\CZ_{00})+\sum_{\substack{(i,j)\in \lambda\\ i\geq 1}}\RR\pi_*\oO_{\Delta_{i-1,j;ij}}(\Delta_{i-1,j;ij})\\
    +\sum_{\substack{(i,j)\in \lambda\\ j\geq 1}}\RR\pi_*\oO_{\Delta_{i,j-1;ij}}(\Delta_{i,j-1;ij}) -\sum_{\substack{(i,j)\in \lambda\\ i,j\geq 1}}\RR\pi_*\oO_{\Delta_{i-1,j-1;ij}}(\Delta_{i-1,j-1;ij}).
\end{multline*}
Moreover, in the expression above, all higher direct images $\RR^{k}\pi_*$ vanish for $k> 0$ by cohomology and base change, therefore
\begin{multline*}
    (\BE^\vee)^{\fix}=\pi_*\oO_{\CZ_{00}}(\CZ_{00})+\sum_{\substack{(i,j)\in \lambda\\ i\geq 1}}\pi_*\oO_{\Delta_{i-1,j;ij}}(\Delta_{i-1,j;ij})\\
    +\sum_{\substack{(i,j)\in \lambda\\ j\geq 1}}\pi_*\oO_{\Delta_{i,j-1;ij}}(\Delta_{i,j-1;ij}) -\sum_{\substack{(i,j)\in \lambda\\ i,j\geq 1}}\pi_*\oO_{\Delta_{i-1,j-1;ij}}(\Delta_{i-1,j-1;ij}).
\end{multline*}
We finally show that this is precisely the same class in $K$-theory as 
\[T_{A_{C,\mathbf{n}_\lambda}}|_{C^{[\mathbf{n}_\lambda]}}-\CE|_{C^{[\mathbf{n}_\lambda]}}\in K_0(C^{[\mathbf{n}_\lambda]}), \]
where $\CE\to A_{C,\mathbf{n}_\lambda} $ is the vector bundle constructed in Theorem \ref{thm: zero locus}. In fact, in the notation of Section \ref{sec: pot on points on curves},
we have
\begin{align*}
    \CX_{ij}|_{C^{[\mathbf{n}_\lambda]}\times C}&=\Delta_{i-1,j;ij}, \quad i\geq 1,\\
     \CY_{ij}|_{C^{[\mathbf{n}_\lambda]}\times C}&=\Delta_{i,j-1;ij}, \quad j\geq 1,\\
      \Gamma^1_{ij}|_{C^{[\mathbf{n}_\lambda]}\times C}&= \Gamma^2_{ij}|_{C^{[\mathbf{n}_\lambda]}\times C}=\Delta_{i-1,j-1;ij}, \quad i,j\geq 1.
\end{align*}
Moreover, the explicit description of the tangent bundle of the Hilbert scheme of points on a smooth curve in terms of its universal subscheme \cite[Lemma IV.2.3]{ACGH_Volume_1} yields
\begin{align*}
    T_{A_{C,\mathbf{n}_\lambda}}=\pi_*\oO_{\CZ_{00}}(\CZ_{00})\oplus\bigoplus_{\substack{(i,j)\in \lambda\\ i\geq 1}}\pi_*\oO_{\CX_{ij}}(\CX_{ij})\oplus\bigoplus_{\substack{(i,j)\in \lambda\\ j\geq 1}}\pi_*\oO_{\CY_{ij}}(\CY_{ij}),
\end{align*}
by which we conclude that
\[     (\BE^\vee)^{\fix}=T_{A_{C,\mathbf{n}_\lambda}}|_{C^{[\mathbf{n}_\lambda]}}-\CE|_{C^{[\mathbf{n}_\lambda]}}\in K_0(C^{[\mathbf{n}_\lambda]}).\]
\end{proof}
In virtue of Proposition \ref{thm: equality virtual classes}, we  denote now on by $T_{C^{[\mathbf{n}_\lambda]}}^{\vir}$ the class in $K$-theory of the dual of the perfect obstruction theory \eqref{eqn: pot coming from zero locus double nested} and the one  induced by the fixed part of  \eqref{eqn: obstruction theory}, which we showed to agree.\\
In order to  compute  stable pair invariants \eqref{eqn: localized PT invariants} one needs to express the virtual normal bundle \eqref{eqn: virtual normal bundle} in terms of $K$-theoretic classes which are easier to handle. For instance, we could express $N^{\vir}$ in terms of pullbacks of the line bundles $L_1,L_2, K_C$ and the universal divisors $\Delta_{ij;lk}$, but that would lead to cumbersome expressions difficult to manipulate. 

\begin{example}\label{ex: weight space t1t2}
As a concrete example, we compute the weight space of $\BE^\vee$ relative to the character  $\tf_1\tf_2$, which we denote by $\BE_{\tf_1\tf_2}^\vee $. 

An application of Grothendieck duality and projection formula on $\pi: C\times C^{[\mathbf{n}_\lambda]}\to C^{[\mathbf{n}_\lambda]} $ yields
\begin{align}\label{eqn: verdier duality on divisors and line bundles}
   \RR\pi_*L(\Delta)= -(\RR\pi_*(K_C\otimes L^{-1}(-\Delta)))^\vee\in K_0(C^{[\mathbf{n}_\lambda]} ), 
\end{align}
for any  divisor $\Delta\subset C\times C^{[\mathbf{n}_\lambda]}$ and any line bundle $L$ on $C$. As in the proof of Proposition \ref{thm: equality virtual classes}, combining \eqref{eqn:betti numbers of graph}, \eqref{eqn: universal exact sequences}, \eqref{eqn: verdier duality on divisors and line bundles}, the identity $K_X=K_C\otimes L_1^{-1}\otimes L^{-1}_2$ and some vanishing of higher direct images yields
\begin{multline*}
    \BE_{\tf_1\tf_2}^\vee=-(\pi_*(K_X\otimes\oO_{\CZ_{00}}(\CZ_{00})))^\vee-\sum_{\substack{(i,j)\in \lambda\\ i\geq 1}}(\pi_*(K_X\otimes\oO_{\Delta_{i-1,j;ij}}(\Delta_{i-1,j;ij})))^\vee\\
    -\sum_{\substack{(i,j)\in \lambda\\ j\geq 1}}(\pi_*(K_X\otimes \oO_{\Delta_{i,j-1;ij}}(\Delta_{i,j-1;ij})))^\vee +\sum_{\substack{(i,j)\in \lambda\\ i,j\geq 1}}(\pi_*(K_X\otimes \oO_{\Delta_{i-1,j-1;ij}}(\Delta_{i-1,j-1;ij})))^\vee.
\end{multline*}
\end{example}
The situation notably simplifies if we impose $X$ to be Calabi-Yau.
 \begin{prop}\label{prop: virtual normal bundle for CY}
 If $X$ is Calabi-Yau, we have an identity  
 \begin{align*}
     N^{\vir}= -T_{C^{[\mathbf{n}_\lambda]}}^{\vir,\vee}\otimes \tf_1\tf_2+ \Omega-\Omega^\vee\otimes \tf_1\tf_2\in K^0_\TT(C^{[\mathbf{n}_\lambda]}),
 \end{align*}
  where $\Omega,\Omega^\vee\in K^0_\TT( C^{[\mathbf{n}_\lambda]})$ 
  have no weight spaces corresponding to the characters $(\tf_1\tf_2)^0, \tf_1\tf_2$.
 \end{prop}
 \begin{proof}
 If $X$ is Calabi-Yau, the perfect obstruction theory \eqref{eqn: obstruction theory} satisfies
 \begin{align}\label{eqn: equivariant serre duality}
\BE^\vee=-\BE\otimes \tf_1\tf_2\in K^0_\TT( C^{[\mathbf{n}_\lambda]})
 \end{align}
      by $\TT$-equivariant Serre duality. Setting $\BE^\vee=W_+-W_-$, where $W_+, W_-\in K^0_\TT( C^{[\mathbf{n}_\lambda]})$ are classes of $\TT$-equivariant vector bundles,  \eqref{eqn: equivariant serre duality} implies that
      \begin{align*}
          W_-=W^\vee_+\otimes \tf_1\tf_2.
      \end{align*}
      We have that $(\BE^\vee)^{\fix}=T^{\vir}_{C^{[\mathbf{n}_\lambda]}}$, therefore
      \begin{align*}
          \BE^\vee=T^{\vir}_{C^{[\mathbf{n}_\lambda]}}-T^{\vir,\vee}_{C^{[\mathbf{n}_\lambda]}}\otimes \tf_1\tf_2+ \Omega-\Omega^\vee\otimes \tf_1\tf_2,
      \end{align*}
      which concludes the argument. 
 \end{proof}
 \begin{remark}\label{rem: remark on Omega explicit}
 A simple computation shows that we could take $\Omega$ to be of the form
 \begin{multline*}
     \Omega=\sum_{\substack{(i,j)\in \lambda\\ (i,j)\neq (0,0)}}\RR\pi_*\left(\oO_{C\times C^{[\mathbf{n}_\lambda]}}(\CZ_{ij}) \otimes L_1^{-i} L_2^{-j}\right) \tf_1^{-i}\tf_2^{-j}\\
     -\sum_{\substack{(i,j),(l,k)\in \lambda\\(i,j)\neq (l,k)\\(i,j)\neq (l+1,k+1)}}\RR\pi_*\left(\oO(\Delta_{ij;lk})\otimes L_1^{i-l}L_2^{j-k} \right)\tf_1^{i-l}\tf_2^{j-k}\\   +\sum_{\substack{(i,j),(l,k)\in \lambda\\(i,j)\neq (l-1,k)\\(i,j)\neq (l,k+1)}}\RR\pi_* \left( \oO(\Delta_{ij;lk})\otimes L_1^{i-l+1}L_2^{j-k} \right)\tf_1^{i-l+1}\tf_2^{j-k}.
 \end{multline*}
 All other choices $\Tilde{\Omega}$ must be of the form
 \begin{align*}
     \Tilde{\Omega}=\Omega+ A+A^\vee\otimes \tf_1\tf_2,
 \end{align*}
 for any $A\in  K^0_\TT( C^{[\mathbf{n}_\lambda]})$  having no weight spaces corresponding to the characters $(\tf_1\tf_2)^0, \tf_1\tf_2$. In particular, this implies that the parity of $\rk \Omega$ is independent by the  choice of $\Omega$.
 \end{remark}

\section{Universality}\label{sec: universal}
\subsection{Universal expression}
In the previous sections, given a triple $(C,L_1,L_2)$ with $C$ an irreducible smooth projective curve and $L_1, L_2$ line bundles on $C$, we reduced stable pair invariants (with no insertions) of $\Tot_C(L_1\oplus L_2)$ to the computation of
\begin{align}\label{eqn: integral on double nested}
    \int_{[C^{[\mathbf{n}_\lambda]}]^{\vir}}e^{\TT}(-N_{C,L_1,L_2}^{\vir})\in \BQ(s_1,s_2),
\end{align}
where $\mathbf{n}_\lambda$ is a reversed plane partition and the virtual normal bundle $N_{C,L_1,L_2}^{\vir}$ is the $\TT$-moving part of the class in $K$-theory \eqref{eqn: dual pot}.\\
We state our main results, describing the generating series of \eqref{eqn: integral on double nested} in terms of three universal functions exploiting the universality techniques used in \cite[Thm. 4.2]{EGL_cobordism} for surfaces. Furthermore, we find  explicit expressions for these universal series under the anti-diagonal restriction $s_1+s_2=0$.
\begin{theorem}\label{thm: universal series}
Let $C$ be a genus $g$ smooth irreducible projective curve and $L_1, L_2$  line bundles over $C$. We have an identity
\begin{align*}
    \sum_{\mathbf{n}_\lambda}q^{|\mathbf{n}_\lambda|} \int_{[C^{[\mathbf{n}_\lambda]}]^{\vir}}e^{\TT}(-N_{C,L_1,L_2}^{\vir})=A_{\lambda}^{g-1}\cdot B_{\lambda}^{\deg L_1}\cdot C_{\lambda}^{\deg L_2}\in \BQ(s_1,s_2)\llbracket q \rrbracket,
\end{align*}
where $A_{\lambda},B_{\lambda},C_{\lambda}\in \BQ(s_1,s_2)\llbracket q \rrbracket$ are fixed universal series for $i=1,2,3$ which only depend on $ \lambda$. Moreover
\begin{align*}
 A_{\lambda}(s_1,s_2)&=A_{\overline{\lambda}}(s_2,s_1),\\
    B_{\lambda}(s_1,s_2)&=C_{\overline{\lambda}}(s_2,s_1).
\end{align*}
\end{theorem}
\begin{proof}
The proof is similar to \cite[Thm. 4.2]{EGL_cobordism}. Consider the map
\[
Z:\CK:=\set{(C,L_1,L_2): C \mbox{ curve}, L_1, L_2 \mbox{ line bundles}}\to \BQ(s_1,s_2)\llbracket q \rrbracket
\]
given by 
\[Z(C,L_1, L_2)=  C_{g,\deg L_1\deg L_2}^{-1}\cdot\sum_{\mathbf{n}_\lambda}q^{|\mathbf{n}_\lambda|} \int_{[C^{[\mathbf{n}_\lambda]}]^{\vir}}e^{\TT}(-N_{C,L_1,L_2}^{\vir}), \]
where $ C_{g,\deg L_1\deg L_2}$ is the leading term of the generating series of the integrals \eqref{eqn: integral on double nested}. 

By Proposition \ref{prop: multiplicativity} the integral \eqref{eqn: integral on double nested} is \emph{multiplicative} and by Corollary \ref{cor: chern dependence integrals}  it is a polynomial  on $g, \deg L_1, \deg L_2$. This implies that $Z$ factors through
\[\CK\xrightarrow{\gamma} \BQ^3\xrightarrow{Z'}  \BQ(s_1,s_2)\llbracket q \rrbracket,\]
where $\gamma(C,L_1, L_2)=(g-1, \deg L_1, \deg L_2)$ and $Z'$ is a linear map.

 A basis of $\BQ^3$ is given by the images
 \begin{align*}
     e_1=\gamma(\BP^1, \oO,\oO),\quad e_2=\gamma(\BP^1, \oO(1),\oO),\quad  e_3=\gamma(\BP^1, \oO,\oO(1)),
 \end{align*}
 and the image of a  generic triple $(C,L_1,L_2) $ can be written as
 \begin{align*}
     \gamma(C,L_1,L_2)=(1-g-\deg L_1-\deg L_2)\cdot e_1+\deg L_1\cdot  e_2+\deg L_2 \cdot e_3.
 \end{align*}
 We conclude that
 \begin{align*}
     Z'(C,L_1,L_2)=Z'(e_1)^{1-g}\cdot (Z'(e_1)^{-1}Z'(e_2))^{\deg L_1}\cdot(Z'(e_1)^{-1}Z'(e_3))^{\deg L_2},
 \end{align*}
 which gives the universal series we were looking for.
 The second claim just follows by interchanging the role of $L_1$ and $L_2$.
\end{proof}
We devote the remainder of Section \ref{sec: universal} to prove the multiplicativity and polynomiality of \eqref{eqn: integral on double nested}. In Section \ref{sec: leading term} we compute the leading term of the generating series of \eqref{eqn: integral on double nested}, while in Section \ref{sec: toric computations} we explicitly compute the integral \eqref{eqn: integral on double nested} in the toric case under the anti-diagonal restriction. These computations will lead to the proof of the second part of the main Theorem \ref{thm: main result intro} (see Theorem \ref{thm: generating series with antidiagonal restriction}).
\subsection{Multiplicativity}
We show now that the integral \eqref{eqn: integral on double nested} is multiplicative. First of all, notice that if $C=C'\sqcup C''$ is a smooth projective curve with two connected components, the construction of Theorem \ref{thm: zero locus} does not directly work and we need to adjust it to define a virtual fundamental class. \\
For any  reversed plane partition $\mathbf{n}_\lambda$  there is an induced stratification
\begin{align*}
    C^{[\mathbf{n}_\lambda]}=\bigsqcup_{\mathbf{n}'_\lambda+\mathbf{n}''_\lambda=\mathbf{n}_\lambda}C'^{[\mathbf{n}'_\lambda]}\times C''^{[\mathbf{n}''_\lambda]}.
\end{align*}
We set 
\begin{align*}
    A_{C,\mathbf{n}_\lambda}:= \bigsqcup_{\mathbf{n}'_\lambda+\mathbf{n}''_\lambda=\mathbf{n}_\lambda}A_{C',\mathbf{n}'_\lambda}\times A_{C'',\mathbf{n}''_\lambda},
\end{align*}
which is a smooth projective variety. Let $\CE_{C',\mathbf{n}'_\lambda},\CE_{C'',\mathbf{n}''_\lambda}$ denote the vector bundles over  $A_{C',\mathbf{n}'_\lambda},A_{C'',\mathbf{n}''_\lambda} $ of Theorem \ref{thm: zero locus}. We define a vector bundle $\CE_{C,\mathbf{n}_\lambda}$ over $A_{C,\mathbf{n}_\lambda}$ by declaring its restriction to any connected component  $A_{C',\mathbf{n}'_\lambda}\times A_{C'',\mathbf{n}''_\lambda} $ to be
\begin{align*}
    \CE_{C,\mathbf{n}_\lambda}|_{A_{C',\mathbf{n}'_\lambda}\times A_{C'',\mathbf{n}''_\lambda}}=\CE_{C',\mathbf{n}'_\lambda}\boxplus \CE_{C'',\mathbf{n}''_\lambda}.
\end{align*}
By Theorem \ref{thm: zero locus}, there exists a section $s$ of  $\CE_{C,\mathbf{n}_\lambda}$ such that 
\begin{align*}
    C^{[\mathbf{n}_\lambda]}\cong Z(s)\hookrightarrow  A_{C,\mathbf{n}_\lambda},
\end{align*}
and therefore an induced virtual fundamental class $[C^{[\mathbf{n}_\lambda]}]^{\vir}$ satisfying
\begin{align}\label{eqn: virt fund disjoint union multipl}
    [C^{[\mathbf{n}_\lambda]}]|^{\vir}_{C'^{[\mathbf{n}'_\lambda]}\times C''^{[\mathbf{n}''_\lambda]}}=[C'^{[\mathbf{n}'_\lambda]}]^{\vir}\boxtimes [C''^{[\mathbf{n}''_\lambda]}]^{\vir}.
\end{align}
By iterating this construction, there exists a natural virtual fundamental class on $C^{[\mathbf{n}_\lambda]}$ for any smooth projective curve $C$ (with any number of connected components).
\begin{prop}\label{prop: multiplicativity}
Let $(C,L_1,L_2)$ be a triple where $C=C'\sqcup C''$ and  $L_i=L'_i\oplus L''_i$ for $i=1,2$, where $L'_i$ are line bundles on $C'$ and $L''_i$ are line bundles on $C''$. Then 
\begin{multline*}
     \sum_{\mathbf{n}_\lambda}q^{|\mathbf{n}_\lambda|} \int_{[C^{[\mathbf{n}_\lambda]}]^{\vir}}e^{\TT}(-N_{C,L_1,L_2}^{\vir})\\= \sum_{\mathbf{n}_\lambda}q^{|\mathbf{n}_\lambda|} \int_{[C'^{[\mathbf{n}_\lambda]}]^{\vir}}e^{\TT}(-N_{C',L'_1,L'_2}^{\vir})\cdot  \sum_{\mathbf{n}_\lambda}q^{|\mathbf{n}_\lambda|} \int_{[C''^{[\mathbf{n}_\lambda]}]^{\vir}}e^{\TT}(-N_{C'',L''_1,L''_2}^{\vir}).
\end{multline*}
\end{prop}
\begin{proof}
Let $\mathbf{n}_\lambda$ be a fixed reversed plane partition. We claim  that the restriction of the virtual normal bundle to the connected component $C'^{[\mathbf{n}'_\lambda]}\times C''^{[\mathbf{n}''_\lambda]}\subset C^{[\mathbf{n}_\lambda]}$ decomposes as
\begin{align}\label{eqn:claim Nvir of disjoint union}
  N^{\vir}_{C,L_1,L_2}|_{C'^{[\mathbf{n}'_\lambda]}\times C''^{[\mathbf{n}''_\lambda]}}= N^{\vir}_{C',L'_1,L'_2}\boxplus N^{\vir}_{C'',L''_1,L''_2}.  
\end{align}
In fact, $N^{\vir}_{C,L_1,L_2}$ is a linear combination of $K$-theoretic classes of the form
\begin{align*}
    \RR\pi_*(\oO_{C\times C^{[\mathbf{n}_\lambda]}}(\Delta)\otimes L_1^aL_2^b)\otimes \tf^\mu\in K^0_\TT(C^{[\mathbf{n}_\lambda]}),
\end{align*}
where $\Delta$ is  a $\BZ$-linear combination of the universal divisors $\CZ_{ij}$  on $C\times C^{[\mathbf{n}_\lambda]}$, $a,b\in \BZ$ and $\tf^\mu$ is a $\TT$-character and notice that
\begin{align*}
    L_1^a L_2^b={L'_1}^a {L'_2}^b\oplus {L''_1}^a {L''_2}^b\in \Pic(C'\sqcup C'').
\end{align*}
Consider the induced  stratification 
\begin{align*}
    C\times C^{[\mathbf{n}_\lambda]}=\bigsqcup_{\mathbf{n}'_\lambda+\mathbf{n}''_\lambda=\mathbf{n}_\lambda}(C'\sqcup C'') \times C'^{[\mathbf{n}'_\lambda]}\times C''^{[\mathbf{n}''_\lambda]}.
\end{align*}
Denote by $\Delta', \Delta''$ the corresponding universal divisor on $C'\times C'^{[\mathbf{n}'_\lambda]}, C''\times C''^{[\mathbf{n}''_\lambda]}$ and consider the projection maps
\begin{align*}
   q_1: C' \times C'^{[\mathbf{n}'_\lambda]}\times C''^{[\mathbf{n}''_\lambda]}&\to C' \times C'^{[\mathbf{n}'_\lambda]},\\
       q_2: C'' \times C'^{[\mathbf{n}'_\lambda]}\times C''^{[\mathbf{n}''_\lambda]}&\to C'' \times C''^{[\mathbf{n}''_\lambda]}.
\end{align*}
On every  component $(C'\sqcup C'') \times C'^{[\mathbf{n}'_\lambda]}\times C''^{[\mathbf{n}''_\lambda]}$ we have
\begin{align*}
    \oO_{C\times C^{[\mathbf{n}_\lambda]}}(\Delta)|_{(C'\sqcup C'') \times C'^{[\mathbf{n}'_\lambda]}\times C''^{[\mathbf{n}''_\lambda]}}=q_1^*\oO_{C'\times C'^{[\mathbf{n}'_\lambda]}}(\Delta')\oplus q_2^*\oO_{C''\times C''^{[\mathbf{n}''_\lambda]}}(\Delta''),
\end{align*}
and similarly
\begin{multline*}
      \oO_{C\times C^{[\mathbf{n}_\lambda]}}(\Delta)\otimes L_1^a L_2^b|_{(C'\sqcup C'') \times C'^{[\mathbf{n}'_\lambda]}\times C''^{[\mathbf{n}''_\lambda]}}=\\
      q_1^*\left(\oO_{C'\times C'^{[\mathbf{n}'_\lambda]}}(\Delta')\otimes {L'_1}^a {L'_2}^b\right)\oplus q_2^*\left(\oO_{C''\times C''^{[\mathbf{n}''_\lambda]}}(\Delta'')\otimes {L''_1}^a {L''_2}^b\right).
\end{multline*}
Consider the cartesian diagram given by the natural projections
\[
\begin{tikzcd}
  C'\times C'^{[\mathbf{n}'_\lambda]}\times C''^{[\mathbf{n}''_\lambda]}\arrow[r, "q_1"]\arrow[d, "\pi"]&C'\times C'^{[\mathbf{n}'_\lambda]}\arrow[d, "\pi_1"]\\
C'^{[\mathbf{n}'_\lambda]}\times C''^{[\mathbf{n}''_\lambda]}\arrow[r, "\Tilde{q}_1"]& C'^{[\mathbf{n}'_\lambda]}.
\end{tikzcd}\]
Flat base change yields
\begin{align*}
     \RR \pi_*q_1^*\left(\oO_{C'\times C'^{[\mathbf{n}'_\lambda]}}(\Delta')\otimes {L'_1}^a {L'_2}^b\right)=\Tilde{q}_1^*\RR {\pi_1}_*\left(\oO_{C'\times C'^{[\mathbf{n}'_\lambda]}}(\Delta')\otimes {L'_1}^a {L'_2}^b\right),
\end{align*}
and analogously for $C''$, which implies
\begin{multline*}
    \RR\pi_*(\oO_{C\times C^{[\mathbf{n}_\lambda]}}(\Delta)\otimes L_1^aL_2^b)\otimes \tf^\mu|_{ C'^{[\mathbf{n}'_\lambda]}\times C''^{[\mathbf{n}''_\lambda]}}\\
    =  \RR{\pi_1}_*(\oO_{C'\times C'^{[\mathbf{n}'_\lambda]}}(\Delta')\otimes {L'_1}^a{L'_2}^b)\otimes \tf^\mu\boxplus   \RR{\pi_2}_*(\oO_{C''\times C''^{[\mathbf{n}''_\lambda]}}(\Delta'')\otimes {L''_1}^a {L''_2}^b)\otimes \tf^\mu,
\end{multline*}
and proves the claim \eqref{eqn:claim Nvir of disjoint union}. Combining \eqref{eqn: virt fund disjoint union multipl} and \eqref{eqn:claim Nvir of disjoint union} yields
\begin{align*}
    \int_{[C^{[\mathbf{n}_\lambda]}]^{\vir}}e^{\TT}(-N_{C,L_1,L_2}^{\vir})&
   &= \sum_{\mathbf{n}'_\lambda+\mathbf{n}''_\lambda=\mathbf{n}_\lambda}\int_{[C'^{[\mathbf{n}'_\lambda]}]^{\vir}}e^{\TT}(-N_{C',L'_1,L'_2}^{\vir})\cdot  \int_{[C''^{[\mathbf{n}''_\lambda]}]^{\vir}}e^{\TT}(-N_{C'',L''_1,L''_2}^{\vir}),
\end{align*}
which concludes the proof.
\end{proof}
\subsection{Chern numbers dependence}
We now show that the integral \ref{eqn: integral on double nested} is a polynomial in the Chern numbers of the triple $(C, L_1, L_2)$. Our strategy is to express the integral on a product of Picard varieties $\Pic^{n}(C)$ --- via the Abel-Jacobi map --- where the integrand is a polynomial expression on  tautological classes. Through this section, we follow the notation as in \cite[Sec. VIII.2]{ACGH_Volume_1} and \cite[Sec. 9, 10.1]{KT_vertical}.
\subsubsection{Tautological integrals on $\Pic^n(C)$}\label{sec: Pic ample}
 Let $C$ be a smooth curve of genus $g$. If $n> 2g-2$,  the Abel-Jacobi map
 \begin{align*}
     \AJ:C^{(n)}&\to \Pic^n(C)\\
     Z&\to [\oO_C(Z)]
 \end{align*}
 is a projective bundle. In fact, consider the diagram
 \[
 \begin{tikzcd}
   C^{(n)}\times C \arrow[d, "\pi"]\arrow[r, "\AJ\times \id"]& \Pic^n(C)\times C\arrow[d, "\overline{\pi}"]\\
    C^{(n)}\arrow[r, "\AJ"] & \Pic^n(C)
 \end{tikzcd}
 \]
 and the Poincaré line bundle $\CP$ on $\Pic^n(C)\times C$, normalized by fixing
 \begin{align*}
     \CP|_{\Pic^n(C)\times \{c\}}\cong \oO_{\Pic^n(C)},
 \end{align*}
 for a certain $c\in C$. Then 
 \begin{align*}
     C^{(n)}\cong \BP(\overline{\pi}_*\CP).
 \end{align*}
 The  universal divisor  $\CZ\subset C^{(n)}\times C$ satisfies \cite[eqn. (67)]{KT_vertical}
 \begin{align*}
      \oO(\CZ)\cong (\AJ\times \id)^*\CP\otimes \pi^*\oO(1) \mbox{ and } \oO(1)\cong \oO(\CZ)|_{C^{(n)}\times \{c\}},
 \end{align*}
and we denote the first Chern class of the latter by
 \[
 \omega:=c_1(\oO(1))\in H^2(C^{(n)}, \BZ).
 \]
 Consider now the product of Abel-Jacobi maps
 \[
 \AJ:C^{(n_1)}\times \dots\times C^{(n_s)}\to \Pic^{n_1}(C)\times \dots \times \Pic^{n_s}(C),
 \]
 where each $n_i> 2g-2$. We denote by $\CP_i$ (the pullback of) the Poincaré line bundle on $\Pic^{n_i}(C)\times C$, each normalized at a point $c_i\in C$, and by $\omega_i$ the first  Chern classes of the tautological bundles on $C^{(n_i)}$. Finally we denote by $\CZ_i\subset C^{(n_i)}\times C$ the universal divisors and by $\CI_i$ (the pullback of) their ideal sheaves, which in this case are line bundles.\\
 We are interested in studying integrals of the form
 \begin{align}\label{eqn: int with f to be exrpessed as universal}
     \int_{C^{(n_1)}\times \dots\times C^{(n_s)}}f,
 \end{align}
 where $f$ is a polynomial in the Chern classes of the $K$-theoretic classes
 \begin{align*}
     \RR\pi_*\RR\hom\left(\bigotimes_{i\in I}\CI_i,\bigotimes_{j\in J}\CI_j\otimes L_k\right),
 \end{align*}
  $L_k$ are line bundles on $C$ and $I,J$ are families of indices (possibly with repetitions).\\
 We assume now that $n_i\gg0$ for all $i=1,\dots, s$. Applying the projection formula and flat base change yields
  \begin{align*}
     &\RR\pi_*\RR\hom\left(\bigotimes_{i\in I}\CI_i,\bigotimes_{j\in J}\CI_j\otimes L_k\right)\\&=  \RR\pi_*\RR\hom\left(\bigotimes_{i\in I}(\AJ\times \id)^*\CP_i^*\otimes \pi^*\oO_{C^{(n_i)}}(-1),\bigotimes_{j\in J}(\AJ\times \id)^*\CP_j^*\otimes \pi^*\oO_{C^{(n_j)}}(-1)\otimes L_k\right)\\
     &= \CF\otimes\RR\pi_*\left(\bigotimes_{i\in I}(\AJ\times \id)^*\CP_i\otimes\bigotimes_{j\in J}(\AJ\times \id)^*\CP_j^*\otimes  L_k \right)\\
     &= \CF\otimes \AJ^*\RR\overline{\pi}_*\left(\bigotimes_{i\in I}\CP_i\otimes\bigotimes_{j\in J}\CP_j^*\otimes  L_k \right),
 \end{align*}
 where $\CF=\bigotimes_{i\in I}\oO_{C^{(n_i)}}(1)\otimes\bigotimes_{j\in J} \oO_{C^{(n_j)}}(-1) $. The Chern classes of the last expression are a linear combination of 
 \[\prod_{i=1}^s\omega_i^{m_i}\cdot \AJ^*\ch_l\left(\RR\overline{\pi}_*\left(\bigotimes_{i\in I}\CP_i\otimes\bigotimes_{j\in J}\CP_j^*\otimes  L_k \right)\right),\]
 where $\ch$ denotes the Chern character for certain $m_i\in\BZ_{\geq 0} $. Integrating this class yields
\begin{multline*}
     \int_{C^{(n_1)}\times \dots\times C^{(n_s)}} \prod_{i=1}^s\omega_i^{m_i}\cdot \AJ^*\ch_l\left(\RR\overline{\pi}_*\left(\bigotimes_{i\in I}\CP_i\otimes\bigotimes_{j\in J}\CP_j^*\otimes  L_k \right)\right)\\
     = \int_{\Pic^{n_1}(C)\times \dots\times \Pic^{n_s}(C)}\AJ_*\prod_{i=1}^s\omega_i^{m_i}\cdot \ch_l\left(\RR\overline{\pi}_*\left(\bigotimes_{i\in I}\CP_i\otimes\bigotimes_{j\in J}\CP_j^*\otimes  L_k \right)\right).
\end{multline*}
Using a standard identity \cite[Sec. 3.1]{Fulton_intersection_theory} we can express the pushforward $\AJ_*\omega_i^{m_i}$ in terms of Segre classes (and therefore Chern characters) of $\overline{\pi}_*\CP_i$. These Chern characters appearing in the integral are computed  by Grothendieck-Riemann-Roch
\begin{align*}
    \ch\left(\RR\overline{\pi}_*\left(\bigotimes_{i\in I}\CP_i\otimes\bigotimes_{j\in J}\CP_j^*\otimes  L_k \right)\right)=\overline{\pi}_*\left(\prod_{i\in I}\ch(\CP_i)\prod_{j\in J}\ch(\CP_j^*)\cdot e^{c_1(L_k)}\cdot \td(T_C) \right).
\end{align*}
The Chern character of the Poincaré line bundle is (cf. \cite[pag. 335]{ACGH_Volume_1})
\begin{align*}
    \ch(\CP_i)=1+n_i[c_i]+\gamma_i-\theta_i[c_i]\in H^*(\Pic^{n_i}(C)\times C, \BZ).
\end{align*}
Here, in the decomposition
\begin{align*}
    H^2(\Pic^{n_i}(C)\times C)\cong H^2(\Pic^{n_i}(C))\oplus\left(H^1(\Pic^{n_i}(C))\otimes H^1(C) \right)\oplus H^2(C),
\end{align*}
we have
\begin{align*}
    &[c_i]\in  H^2(C),\\
    &\theta_i\in H^2(\Pic^{n_i}(C)),\\
    &\gamma_i=[\Delta]^{1,1}\in H^1(C)\otimes H^1(C)\cong H^1(\Pic^{n_i}(C))\otimes H^1(C),
\end{align*}
where $\Delta\subset C\times C$ is the diagonal and $\theta_i$ is the theta divisor. All of this results in 
\begin{align*}
    \int_{C^{(n_1)}\times \dots\times C^{(n_s)}}f=
      \int_{\Pic^{n_1}(C)\times \dots\times \Pic^{n_s}(C)\times C}\Tilde{f},
\end{align*}
where $\Tilde{f}$ is a polynomial in the classes
\begin{align*}
    \theta_i,\, [c_i],\, \gamma_i,\, c_1(T_C),\, c_1(L_k).
\end{align*}
\subsubsection{Extension to all $n$}\label{sec: Pic any n}
In the previous section we assumed that all $n_i> 2g-2$; we explain now how to remove this assumption, following closely \cite[Sec. 10.1]{KT_vertical}.\\

Let $n\in \BZ_{\geq 0}$ and $N>2g-2$. Then $C^{(N)}\cong \BP(\overline{\pi}_*\CP)$, where $\CP$ is the Poincaré line bundle on $\Pic^N(C)\times C$ normalized at $c\in C$.  We can embed  $C^{(n)}$ in $C^{(N)}$ as the zero section of a vector bundle
\[
\begin{tikzcd}
  & \pi_*\oO(\CW)|_{C^{(N)}\times (N-n)c}\arrow[d]\\
  C^{(n)}\cong Z(s)\arrow[r, hook,"\iota"]&C^{(N)},\arrow[u, bend right, swap,  "s"]
\end{tikzcd}
\]
by sending $Z\mapsto Z+(N-n)c$, where $\CW\subset C^{(N)}\times C$ is the universal divisor and $ (N-n)c\subset C$ is an Artinian thickened point. Moreover, the section is regular, therefore
\begin{align*}
    \iota_*[C^{(n)}]=e(\CG)\cap [C^{(N)}]\in H_*(C^{(N)}),
\end{align*}
where $\CG= \pi_*\oO(\CW)|_{C^{(N)}\times (N-n)c}$. Finally, if we denote by $\CZ\subset C^{(n)}\times C$ the universal divisor, we have that $(\iota\times \id)^*\CW=\CZ((N-n)c)$.\\

Recall that we are interested in the integrals \eqref{eqn: int with f to be exrpessed as universal}. Choose $N_i> 2g-2$, denote by $\CI'_i$ the universal ideal sheaves of the universal divisors $\CW_i$ on $C^{(N_1)}\times\dots \times C^{(N_s)}\times C$ and by $\pi'$ the projection map. By base change we can write
\begin{multline*}
     \RR\pi_*\RR\hom\left(\bigotimes_{i\in I}\CI_i,\bigotimes_{j\in J}\CI_j\otimes L_k\right)=\\ \mathbf{L}\iota^*\RR\pi'_*\RR\hom\left(\bigotimes_{i\in I}\CI'_i(-(N_i-n_i)c_i),\bigotimes_{j\in J}\CI'_j(-(N_j-n_j)c_j)\otimes L_k\right),
 \end{multline*}
 therefore
 \begin{align*}
      \int_{C^{(n_1)}\times \dots\times C^{(n_s)}}f= \int_{C^{(N_1)}\times \dots\times C^{(N_s)}} f'\cdot \prod_{i=1}^s e(\CG_i),
 \end{align*}
 where each $\CG_i=\pi_*\oO(\CW_i)|_{C^{(N_i)}\times (N_i-n_i)c_i}$ and $f'$ is a polynomial in the Chern classes of \[\RR\pi'_*\RR\hom\left(\bigotimes_{i\in I}\CI'_i,\bigotimes_{j\in J}\CI'_j\otimes \bigotimes_{i\in I} \oO((N_i-n_i)c_i) \otimes\bigotimes_{j\in J}\oO(-(N_j-n_j)c_j)  \otimes L_k\right) .\]
 The exact sequence
 \begin{align*}
     0\to \oO(\CW_i-(N_i-n_i)c_i)\to \oO(\CW_i)\to \oO(\CW_i)|_{C^{(N_i)}\times (N_i-n_i)c_i}\to 0
 \end{align*}
 yields the identity in $K$-theory
 \begin{align*}
     \CG_i=\RR\pi'_*\CI'^*_i-\RR\pi'_*\left(\CI'^*_i\otimes \oO((N_i-n_i)c_i)\right),
 \end{align*}
 by which we conclude that we can apply the construction in Section \ref{sec: Pic ample} to express
 \begin{align*}
    \int_{C^{(n_1)}\times \dots\times C^{(n_s)}}f=
      \int_{\Pic^{N_1}(C)\times \dots\times \Pic^{N_s}(C)  \times C}\Tilde{f},
\end{align*}
where $\Tilde{f}$ is a polynomial in the classes
\begin{align}\label{eqn: taut classes Pic}
   \theta_i,\, [c_i],\, \gamma_i,\, c_1(T_C),\, c_1(L_k), \, c_1(\oO((N_i-n_i)c_i) ).
\end{align}
Smooth projective curves of the same genus are diffeomorphic to each other, therefore $\Pic^n(C)$ is diffeomorphic to a $g$-dimensional complex torus. By the intersection theory on $\Pic^n(C)$ developed in \cite[Sec. VIII.2]{ACGH_Volume_1} we immediately obtain the following result.
\begin{prop}\label{prop: taut dependence cobordism class with Pic}
Let $f$ be a polynomial in the Chern classes of the $K$-theory classes
 \begin{align*}
     \RR\pi_*\RR\hom\left(\bigotimes_{i\in I}\CI_i,\bigotimes_{j\in J}\CI_j\otimes L_k\right),
 \end{align*}
where  $L_k$ are line bundles on $C$ and $I,J$ are families of indices (possibly with repetitions). Then
 \begin{align*}
    \int_{C^{(n_1)}\times \dots\times C^{(n_s)}}f= \int_{\Pic^{n_1}(C)\times \dots\times \Pic^{n_s}(C)\times C}\Tilde{f},
\end{align*}
where $\Tilde{f}$ is a polynomial in the classes \eqref{eqn: taut classes Pic}. In particular, the integral is a polynomial in the genus $g=g(C)$ and the degrees of the line bundles $\{L_k\}_k$.
\end{prop}
As a corollary, we obtain that the localized contributions on $C^{[\mathbf{n}_\lambda]}$ only depend on the Chern numbers of $(C, L_1, L_2)$.
\begin{corollary}\label{cor: chern dependence integrals}
Let $C$ be a genus $g$  irreducible smooth projective curve and $L_1, L_2$ line bundles on $C$. Then the intersection numbers \eqref{eqn: integral on double nested} are polynomials in $g, \deg L_1, \deg L_2$.
\end{corollary}
\begin{proof}
Let $i:C^{[\mathbf{n}_\lambda]}\hookrightarrow A_{C, \mathbf{n}_\lambda}$ be the embedding of Theorem \ref{thm: zero locus} and to ease the notation set $A_{C, \mathbf{n}_\lambda}=C^{(n_1)}\times\dots \times  C^{(n_s)} $. We claim that
\begin{align*}
N^{\vir}_{C, L_1, L_2}=i^* \Tilde{N}^{\vir}_{C, L_1, L_2}, 
\end{align*}
for a certain class $\Tilde{N}^{\vir}_{C, L_1, L_2}\in K^0_\TT(A_{C, \mathbf{n}_\lambda})$. In fact, $ N^{\vir}_{C, L_1, L_2}$ is a linear combination of classes in $K$-theory of the form 
\begin{align*}
    \RR\pi_*(\oO_{C\times C^{[\mathbf{n}_\lambda]}}(\Delta)\otimes L_1^aL_2^b)\otimes \tf^\mu\in K^0_\TT(C^{[\mathbf{n}_\lambda]}),
\end{align*}
 where $\Delta$ is  a $\BZ$-linear combination of the universal divisors $\CZ_{ij}$  on $C\times C^{[\mathbf{n}_\lambda]}$, $a,b\in \BZ$ and $\tf^\mu$ is a $\TT$-character. Each of such universal divisors $\Delta$ can be expressed, in the Picard group of $ C^{[\mathbf{n}_\lambda]}$, as a linear combination of the divisors $  i^*\CZ_{00}, i^*\CX_{ij}, i^*\CY_{ij}$, with notation as in Section \ref{sec: pot on points on curves}; applying base change  proves the claim. Therefore the integral \eqref{eqn: integral on double nested} can be expressed as an intersection number on the product of symmetric powers of curves $A_{C, \mathbf{n}_\lambda}$
\begin{align*}
    \int_{[C^{[\mathbf{n}_\lambda]}]^{\vir}}e^{\TT}(-N_{C,L_1,L_2}^{\vir})=\int_{A_{C, \mathbf{n}_\lambda}}e^{\TT}(\CE-\Tilde{N}_{C,L_1,L_2}^{\vir}),
\end{align*}
where $\CE$ is the vector bundle of Theorem \ref{thm: zero locus}. In particular, the $K$-theory class of $\CE-\Tilde{N}_{C,L_1,L_2}^{\vir}$ is a linear combination of classes of the form 
 \begin{align*}
     \RR\pi_*\RR\hom\left(\bigotimes_{i\in I}\CI_i,\bigotimes_{j\in J}\CI_j\otimes L_k\right)\otimes \tf^{\mu},
 \end{align*}
  where $L_k$ are line bundles on $C$ and $I,J$ are families of indices (possibly with repetitions). By Proposition \ref{prop: taut dependence cobordism class with Pic} this integral is a polynomial in $g$ and the degrees of $L_k$. We conclude the proof by noticing that all line bundles $L_k$ possibly occuring are a linear combination of $L_1, L_2, K_C$.
\end{proof}
\section{The leading term}\label{sec: leading term}
We compute the leading term of the generating series of the integrals  \eqref{eqn: integral on double nested}, which is essential for the computation of the full generating series in Theorem \ref{thm: generating series with antidiagonal restriction}.
\begin{prop}\label{prop: leading term}
Let $C$ be a smooth projective curve of genus $g$ and $L_1, L_2$ line bundles on $C$. Then under the anti-diagonal restriction $s_1+s_2=0$ we have
\begin{align*}
    \int_{C^{[\mathbf{0}_{\lambda}]}}e^{\TT}(-N_{C,L_1,L_2}^{\vir})=(-1)^{|\lambda|(g-1+\deg L_2)+n(\lambda)\deg L_1+n(\overline{\lambda})\deg L_2}\left(s_1^{|\lambda|}\cdot\prod_{\Box\in \lambda}h(\Box)\right)^{2g-2-\deg L_1-\deg L_2},
\end{align*}
where  $\mathbf{0}_{\lambda}$ is the unique reversed plane partition of size $0$ and underlying Young diagram $\lambda$.
\end{prop}
\begin{proof}
We have that $C^{[\mathbf{0}_{\lambda}]}\cong \pt $ and $[C^{[\mathbf{0}_{\lambda}]}]^{\vir}=[\pt]\in A_*(\pt)$, therefore $\BE^\vee$ is completely $\TT$-movable and
\begin{align*}
    \int_{C^{[\mathbf{0}_{\lambda}]}}e^{\TT}(-N_{C,L_1,L_2}^{\vir})=e^\TT(-N_{C,L_1,L_2}^{\vir})\in \BQ(s_1, s_2),
\end{align*}
where by \eqref{eqn: dual pot}, \eqref{eqn: verdier duality on divisors and line bundles} we express the class in $K$-theory of the virtual normal bundle as
 \begin{multline*}
     N_{C,L_1,L_2}^{\vir}=\sum_{(i,j)\in \lambda}\left(\RR\Gamma(L_1^{-i}L_2^{-j}) \tf_1^{-i}\tf_2^{-j}-\left(\RR\Gamma( K_C L_1^{-i-1}L_2^{-j-1}) \tf_1^{-i}\tf_2^{-j}\right)^\vee \otimes \tf_1\tf_2\right)\\
     -\sum_{(i,j),(l,k)\in \lambda}\left(\RR\Gamma((\oO_C-L_1\tf_1)L_1^{i-l}L_2^{j-k})\tf_1^{i-l}\tf_2^{j-k}-\left(\RR\Gamma((\oO_C-L_1\tf_1)K_CL_1^{i-l-1}L_2^{j-k-1})\tf_1^{i-l}\tf_2^{j-k}\right)^\vee\otimes \tf_1\tf_2\right).
 \end{multline*}
Applying Riemann-Roch,  every line bundle $L$ on $C$ satisfies
\begin{align*}
    \RR\Gamma(L)=&\BC^{\deg L +1-g}\in K^0(\pt)\cong \BZ,
\end{align*}
therefore we can write the virtual normal bundle as
\begin{multline}\label{eqn: Nvir con mu nu}
    N_{C,L_1,L_2}^{\vir}=\sum_\mu(\BC^{m_{\mu}+1-g}\tf^\mu-\BC^{m_{\mu}+ g-1-\deg L_1-\deg L_2}\tf^{-\mu}\tf_1\tf_2)\\-\sum_\nu(\BC^{m_{\nu}+1-g}\tf^\nu-\BC^{m_{\nu}+ g-1-\deg L_1-\deg L_2}\tf^{-\nu}\tf_1\tf_2),
\end{multline}
where the weights $\mu$ range among
\begin{align}\label{eqn: mu}
\begin{split}
    &(-i,-j)\in\BZ^2 \mbox{ such that } (i,j)\in \lambda, (i,j)\neq (0,0), \\
    &(i-l+1,j-k)\in \BZ^2 \mbox{ such that } ((i-j),(l-k))\in \lambda, (i,j)\neq ((l-1,k)(l,k+1)),
\end{split}
\end{align}
  the weights $\nu$ range among
\begin{align}\label{eqn: nu}
    (i-l,j-k)\in   \BZ^2 \mbox{ such that } ((i-j),(l-k))\in \lambda, (i,j)\neq ((l,k),(l+1,k+1)),
\end{align}
and, for a weight $\mu=(\mu_1,\mu_2)$,  we set  $m_{\mu}=\mu_1 \deg L_1+\mu_2 \deg L_2$.
In particular, the weights $1,\tf_1\tf_2$ do not appear in \eqref{eqn: Nvir con mu nu}, as the virtual tangent bundle of the $\TT$-fixed locus has rank 0 and by the explicit description of the weight space of $\tf_1\tf_2 $ in Example \ref{ex: weight space t1t2}. For every weight $\mu$, we compute
\begin{align*}
    e^\TT(\BC^{m_{\mu}+1-g}\tf^\mu-\BC^{m_{\mu}+ g-1-\deg L_1-\deg L_2}\tf^{-\mu}\tf_1\tf_2)=\frac{(\mu\cdot s)^{m_{\mu}+1-g}}{(-\mu\cdot s+s_1+s_2)^{m_{\mu}+ g-1-\deg L_1-\deg L_2}},
\end{align*}
where $s=(s_1,s_2)$ and $\mu\cdot s$ is the standard inner product. With a simple manipulation we end up with
\begin{multline*}
    e^\TT(-N_{C,L_1,L_2}^{\vir})=(-1)^{\sum_{\nu\neq (a,a)}m_{\nu}+\sum_{\mu\neq (b,b)}m_{\mu}+(1-g+\deg L_1+\deg L_2)(\#\set{\nu\neq (a,a)}-\#\set{\mu\neq (b,b)})}
    \\\cdot\frac{\prod_{\nu\neq (a,a)}(\nu\cdot s)^{m_\nu+ 1-g}(\nu\cdot s-s_1-s_2)^{-m_{\nu}+1-g+\deg L_1+\deg L_2}}{\prod_{\mu\neq (b,b)}(\mu\cdot s)^{m_\mu+ 1-g}(\mu\cdot s-s_1-s_2)^{-m_{\mu}+1-g+\deg L_1+\deg L_2}}\\
    \cdot\frac{\prod_{\nu= (a,a)} a^{a(\deg L_1+\deg L_2)+1-g}(1-a)^{(1-a)(\deg L_1+\deg L_2)+1-g}(s_1+s_2)^{2-2g+\deg L_1+\deg L_2}}{\prod_{\mu= (b,b)} b^{b(\deg L_1+\deg L_2)+1-g}(1-b)^{(1-b)(\deg L_1+\deg L_2)+1-g}(s_1+s_2)^{2-2g+\deg L_1+\deg L_2}}.
\end{multline*}
Following the proof of  Lemma \ref{lemma: induction on poly for no powers CY}, we see that for any weight $\nu=(a,a)$ (with $a\in \BZ$) there is either a weight $\mu=(a,a)$ or $\mu=(1-a,1-a)$ (and viceversa). This implies that 
\begin{align*}
    \#\set{\nu= (a,a)}=\#\set{\mu= (b,b)}
\end{align*}
and
\begin{align*}
    \frac{\prod_{\nu= (a,a)} a^{a(\deg L_1+\deg L_2)+1-g}(1-a)^{(1-a)(\deg L_1+\deg L_2)+1-g}(s_1+s_2)^{2-2g+\deg L_1+\deg L_2}}{\prod_{\mu= (b,b)} b^{b(\deg L_1+\deg L_2)+1-g}(1-b)^{(1-b)(\deg L_1+\deg L_2)+1-g}(s_1+s_2)^{2-2g+\deg L_1+\deg L_2}}=1.
\end{align*}
In particular, the anti-diagonal restriction $s_1+s_2=0$ is well-defined; we get
\begin{multline*}
   \left. \frac{\prod_{\nu\neq (a,a)}(\nu\cdot s)^{m_\nu+ 1-g}(\nu\cdot s-s_1-s_2)^{-m_{\nu}+1-g+\deg L_1+\deg L_2}}{\prod_{\mu\neq (b,b)}(\mu\cdot s)^{m_\mu+ 1-g}(\mu\cdot s-s_1-s_2)^{-m_{\mu}+1-g+\deg L_1+\deg L_2}}\right|_{s_1+s_2=0}\\
   =\left(s_1^{\#\nu-\#\mu}\cdot\frac{\prod_{\nu= (\nu_1,\nu_2)}(\nu_1-\nu_2)}{\prod_{\mu= (\mu_1,\mu_2)}(\mu_1-\mu_2)}\right)^{2-2g+\deg L_1+\deg L_2}
\end{multline*}
where the product is over all $\mu,\nu\neq (a,a)$. Moreover, it is immediate to check that $  \#\mu-\#\nu=|\lambda|$. By Lemma \ref{lemma: hook length in size 0}, we conclude that the last expression equals to
\begin{align*}
    \sigma(\lambda)^{\deg L_1+\deg L_2}\cdot\left(s_1^{-|\lambda|}\cdot\prod_{\Box\in \lambda}h(\Box)^{-1}\right)^{2-2g+\deg L_1+\deg L_2},
\end{align*}
where $\sigma(\lambda)$ is defined in Lemma \ref{lemma: hook length in size 0}. We are only left with a sign computation; we conclude by Lemma \ref{lemma: sign for leading term}
\begin{align*}
    (-1)^{\sum_{\nu\neq (a,a)}m_{\nu}+\sum_{\mu\neq (b,b)}m_{\mu}+(1-g+\deg L_1+\deg L_2)|\lambda|}\sigma(\lambda)^{\deg L_1+\deg L_2}&= (-1)^{\rho(\lambda)+ |\lambda|(k_1+k_2+ g-1)}\sigma(\lambda)^{k_1+k_2}\\
    &=(-1)^{|\lambda|(k_2+g-1)+ n(\lambda)k_1+n(\overline{\lambda})k_2},
\end{align*}
where $\rho(\lambda)$ was defined is Lemma \ref{lemma: sign for leading term} and we set $k_i=\deg L_i$.
\end{proof}

\section{Toric computations}\label{sec: toric computations}
\subsection{Torus action}
Let $U_0, U_\infty $ be the standard open cover of $\BP^1$. The torus $\BC^*$ acts on the coordinate functions of $\BP^1$ as $t\cdot x=tx$ (resp. $t\cdot x=t^{-1}x $) in the chart $U_0$ (resp. $U_\infty$). The $\BC^*$-representation of the tangent space at the $\BC^*$-fixed points of $\BP^1$ is
\begin{align*}
    T_{\BP^1,0}&=\BC\otimes \tf^{-1},\\
    T_{\BP^1,\infty}&=\BC\otimes \tf.
\end{align*}
We prove some identities of $\BC^*$-representations that will  be useful later in this section.
\begin{lemma}\label{lemma: identities representations}
Let $Z=Z_0\sqcup Z_\infty\subset \BP^1$ be a closed subscheme, where $Z_0$ (resp. $Z_\infty$) is a closed subscheme of length $n_0$ (resp. $n_\infty$) supported on $0$ (resp. $\infty$). For any $a\in \BZ$, we have the following identities in $K^0_{\BC^*}(\pt)$
\begin{align*}
  &\RR\Gamma(K_{\BP^1}^a)=\begin{cases}\sum_{i=a}^{-a}\tf^i & a\leq 0,\\
     -\sum_{i=-a+1}^{a-1}\tf^i & a\geq  1,
     \end{cases}\\
    &\RR\Gamma(\oO_Z(Z))= \sum_{i=1}^{n_0}\tf^{-i}+\sum_{i=1}^{n_\infty}\tf^{i},\\
    &\RR\Gamma(\oO_Z\otimes K_{\BP^1}^a )= \sum_{i=0}^{n_0-1}\tf^{a+i}+\sum_{i=0}^{n_\infty-1}\tf^{-a-i},\\
    &\RR\Hom(\oO_Z,K_{\BP^1}^a )=-\sum_{i=0}^{n_0-1}\tf^{a-i-1}-\sum_{i=0}^{n_\infty-1}\tf^{-a+i+1}.
\end{align*}
\end{lemma}
\begin{proof}
We have that 
\begin{align*}
     \RR\Gamma(K_{\BP^1}^a)&=\chi(\BP^1,K_{\BP^1}^a )\\
     &=\frac{\tf^a}{1-\tf}-\frac{\tf^{-a}}{1-\tf^{-1}}\\
     &=\begin{cases}\sum_{i=a}^{-a}\tf^i & a\leq 0,\\
     -\sum_{i=-a+1}^{a-1}\tf^i & a\geq  1,
     \end{cases}
\end{align*}
where in the second line we applied the classical $K$-theoretic localization \cite{Tho:formule_Lefschetz} on $\BP^1$.\\
Secondly, we have
\begin{align*}
   \RR\Gamma(\oO_Z(Z))&=H^0(\BP^1,\oO_Z(Z))\\
    &= H^0(U_0,\oO_{Z_0}(Z_0)|_{U_0})+H^0(U_\infty,\oO_{Z_{\infty}}(Z_\infty)|_{U_\infty})\\
    &= H^0(U_0, \oO_{Z_0}\otimes K_{\BP^1}|_{U_0} )^*+H^0(U_\infty, \oO_{Z_\infty}\otimes K_{\BP^1}|_{U_\infty} )^*\\
    &= \sum_{i=1}^{n_0}\tf^{-i}+\sum_{i=1}^{n_\infty}\tf^{i},
\end{align*}
where in the second line we used  \v{C}ech cohomology and in the third line we used \cite[Ex. 3.4.5]{Okounk_Lectures_K_theory}.\\
Thirdly, by applying \v{C}ech cohomology as before we have
\begin{align*}
    \RR\Gamma(\oO_Z\otimes K_{\BP^1}^a)&= H^0(U_0, \oO_{Z_0}\otimes K_{\BP^1}^a|_{U_0} )+H^0(U_\infty, \oO_{Z_\infty}\otimes K_{\BP^1}^a|_{U_\infty} )\\
    &= \sum_{i=0}^{n_0-1}\tf^{a+i}+\sum_{i=0}^{n_\infty-1}\tf^{-a-i}.
\end{align*}
Finally, combining Serre duality and the previous result  yields
\begin{align*}
     \RR\Hom(\oO_Z, K_{\BP^1}^a)&=-\RR\Gamma(\oO_Z\otimes K_{\BP^1}^{1-a})^*\\
     &=-\sum_{i=0}^{n_0-1}\tf^{a-i-1}-\sum_{i=0}^{n_\infty-1}\tf^{-a+i+1}.
\end{align*}
\end{proof}
The $\BC^*$-action on $\BP^1$ naturally lifts to a $\BC^*$-action on the Hilbert scheme of point ${\BP^1}^{[n]}$, whose $\BC^*$-fixed locus  consists of length $n$ closed subschemes $Z\subset \BP^1$ supported on $0, \infty$. Therefore, there is an induced  $\BC^*$-action on $A_{\BP^1, \mathbf{n}_\lambda}$, whose $\BC^*$-fixed locus is 0-dimensional and reduced. This $\BC^*$-action restricts to a  $\BC^*$-action on $ {\BP^{1}}^{[\mathbf{n}_\lambda]}$, whose $\BC^*$-fixed locus is necessarily  0-dimensional and reduced. Moreover, the perfect obstruction theory \eqref{eqn: pot coming from zero locus double nested} is naturally $\BC^*$-equivariant, as all the ingredients of Theorem \ref{thm: zero locus} are.
\begin{prop}\label{prop: C^* pot is movable}
Let $\underline{Z}\in  {\BP^{1}}^{[\mathbf{n}_\lambda], \BC^*} $ be a $\BC^*$-fixed point. Then $ T_{{\BP^{1}}^{[\mathbf{n}_\lambda]}, \underline{Z}}^{\vir}$ is completely $\BC^*$-movable. In particular,  the induced perfect obstruction theory on $ {\BP^{1}}^{[\mathbf{n}_\lambda], \BC^*}$ is trivial.
\end{prop}
\begin{proof}
We need to show that the class in $K$-theory
\begin{align*}
    T_{A_{C,\mathbf{n}_\lambda}}|_{\underline{Z}}-\CE|_{\underline{Z}}\in K^0_{\BC^*}(\pt)
\end{align*}
does not have $\BC^*$-fixed part. Recall that we have an identity in $K^0_{\BC^*}(\pt)$
\begin{align*}
    T_{A_{C,\mathbf{n}_\lambda}}|_{\underline{Z}}=\RR\Gamma(\BP^1,\oO_{Z_{00}}(Z_{00}))+\sum_{\substack{(i,j)\in \lambda\\ i\geq 1}}\RR\Gamma(\BP^1,\oO_{X_{ij}}(X_{ij}))+\sum_{\substack{(i,j)\in \lambda\\ j\geq 1}}\RR\Gamma(\BP^1,\oO_{Y_{ij}}(Y_{ij})),
\end{align*}
where for simplicity we denoted by $X_{ij}= Z_{ij}-Z_{i-1,j}$ and by $Y_{ij}=Z_{ij}-Z_{i,j-1}$. Moreover we have
\begin{align*}
    \CE|_{\underline{Z}}=\sum_{\substack{(i,j)\in \lambda\\ i,j\geq 1}}\RR\Gamma(\BP^1,\oO_{W_{ij}  }(W_{ij}))\in K^0_{\BC^*}(\pt),
\end{align*}
where for simplicity we denoted by $W_{ij}= Z_{ij}-Z_{i-1,j-1}$. Therefore the virtual tangent space is a sum of classes of the form $ \RR\Gamma(\BP^1,\oO_Z(Z))$, with $Z\subset \BP^1$ a closed subscheme, which is entirely $\BC^*$-movable by the description in Lemma \ref{lemma: identities representations}.
\end{proof}
\subsection{Case I: Calabi-Yau}
We compute  the integral \eqref{eqn: integral on double nested} for $C=\BP^1$ in the case of $L_1\otimes L_2=K_{\BP^1} $, showing that it  coincides (up to a sign) with the topological Euler characteristic $ e\left({\BP^{1}}^{[\mathbf{n}_\lambda]}\right)$. 
\begin{theorem}\label{thm: g=0 case}
Let $L_1,L_2$ be line bundles on $\BP^1$ such that $L_1\otimes L_2=K_{\BP^{1}} $. For any reversed plane partition $\mathbf{n}_\lambda$, we have
\begin{align*}
   \left.\left( \int_{\left[{\BP^{1}}^{[\mathbf{n}_\lambda]}\right]^{\vir}}e^\TT(-N_{\BP^{1},L_1,L_2}^{\vir})
   \right)\right|_{s_1+s_2=0}=(-1)^{\deg L_1(c_\lambda+|\lambda|)+ |\lambda|+|\mathbf{n}_\lambda|}\cdot e\left({\BP^{1}}^{[\mathbf{n}_\lambda]}\right).
\end{align*}
\end{theorem}
\begin{proof}
By Graber-Pandharipande   \cite{GP_virtual_localization}, there is an induced perfect obstruction theory and virtual fundamental class on the $\BC^*$-fixed locus ${\BP^{1}}^{[\mathbf{n}_\lambda], \BC^*}$, both trivial by Proposition \ref{prop: C^* pot is movable}. By Proposition \ref{prop: virtual normal bundle for CY}
\begin{align*}
      N_{\BP^{1},L_1,L_2}^{\vir}= -T_{{\BP^1}^{[\mathbf{n}_\lambda]}}^{\vir,\vee}\otimes \tf_1\tf_2+ \Omega-\Omega^\vee\otimes \tf_1\tf_2\in K^0_\TT\left({\BP^1}^{[\mathbf{n}_\lambda]}\right),
\end{align*}
and applying the virtual localization formula with respect to the $\BC^*$-action yields 
\begin{align*}
    \int_{\left[{\BP^{1}}^{[\mathbf{n}_\lambda]}\right]^{\vir}}
    e^\TT(-N_{\BP^{1},L_1,L_2}^{\vir})
     &=\left.\left(\sum_{\underline{Z}\in {\BP^{1}}^{[\mathbf{n}_\lambda],\BC^*}}\frac{e^{\TT\times \BC^*}(T_{{\BP^{1}}^{[\mathbf{n}_\lambda]}, \underline{Z}}^{\vir,\vee}\otimes \tf_1\tf_2)}{e^{\TT\times \BC^*}(T_{{\BP^{1}}^{[\mathbf{n}_\lambda]}, \underline{Z}}^{\vir})}\cdot\frac{e^{\TT\times \BC^*}(\Omega|_{\underline{Z}}^\vee\otimes \tf_1\tf_2)}{e^{\TT\times \BC^*}(\Omega|_{\underline{Z}})}\right)\right|_{s_3=0}\\
    &=\left.\left(\sum_{\underline{Z}\in {\BP^{1}}^{[\mathbf{n}_\lambda],\BC^*}}\frac{e^{\TT\times \BC^*}(\mathsf{V}_{\underline{Z}}^\vee\otimes\tf_1\tf_2) }{e^{\TT\times \BC^*}(\mathsf{V}_{\underline{Z}})}\right)\right|_{s_3=0},
\end{align*}
where $s_3$ is the generator of the equivariant cohomology $H^*_{\BC^*}(\pt)$ and  we denoted by 
\begin{align*}
    \mathsf{V}_{\underline{Z}}= T_{{\BP^{1}}^{[\mathbf{n}_\lambda]}, \underline{Z}}^{\vir}+ \Omega|_{\underline{Z}}\in K^0_{\TT\times \BC^*}(\pt)
\end{align*}
the $\BC^*$-equivariant lift of the $K$-theoretic class $T_{{\BP^{1}}^{[\mathbf{n}_\lambda]}+ \underline{Z}}^{\vir}, \Omega|_{\underline{Z}}\in K^0_{\TT}(\pt) $. Under the Calabi-Yau restriction $s_1+s_2=0$, we have by Lemma \ref{lemma: no weight spaces on Omega}
\begin{align*}
    \left.\left(\frac{e^{\TT\times \BC^*}(\mathsf{V}_{\underline{Z}}^\vee\otimes\tf_1\tf_2) }{e^{\TT\times \BC^*}(\mathsf{V}_{\underline{Z}})}\right)\right|_{s_1+s_2=0}=(-1)^{\rk \mathsf{V}_{\underline{Z}}}.
\end{align*}
 Moreover, by Lemma \ref{lemma:rank Omega}
\begin{align*}
   \rk \mathsf{V}_{\underline{Z}}=\deg L_1(c_\lambda+|\lambda|)+ |\lambda|+|\mathbf{n}_\lambda| \mod 2.
\end{align*}
Therefore, we conclude that
\begin{align*}
        \left.\left(\sum_{\underline{Z}\in {\BP^{1}}^{[\mathbf{n}_\lambda],\BC^*}}\frac{e^{\TT\times \BC^*}(\mathsf{V}_{\underline{Z}}^\vee\otimes\tf_1\tf_2) }{e^{\TT\times \BC^*}(\mathsf{V}_{\underline{Z}})}\right)\right|_{s_+s_2=0}&=\sum_{\underline{Z}\in {\BP^{1}}^{[\mathbf{n}_\lambda],\BC^*}} (-1)^{\rk \mathsf{V}_{\underline{Z}}}\\
        &= (-1)^{\deg L_1(c_\lambda+|\lambda|)+ |\lambda|+|\mathbf{n}_\lambda|}\cdot e\left( {\BP^{1}}^{[\mathbf{n}_\lambda]}\right),
\end{align*}
as the Euler characteristic of a  $\BC^*$-scheme coincides with the Euler characteristic of its $\BC^*$-fixed locus (in our case the number of fixed points).
\end{proof}
Exploiting the close formula for the generating series of the topological Euler characteristic proved in  Theorem \ref{thm: euler char double nested}, we derive the following close expression in the Calabi-Yau case.
\begin{corollary}\label{cor: gen series CY}
Let $L_1,L_2$ be line bundles on $\BP^1$ such that $L_1\otimes L_2=K_{\BP^{1}} $ and $ \lambda$ be a Young diagram. We have
\begin{align*}
   \sum_{\mathbf{n}_\lambda}q^{|\mathbf{n}_\lambda|}\cdot \left.\left( \int_{\left[{\BP^{1}}^{[\mathbf{n}_\lambda]}\right]^{\vir}}e^{\TT}(-N_{C, L_1, L_2}^{\vir})\right)\right|_{s_1+s_2=0}=(-1)^{\deg L_1(c_\lambda+|\lambda|)+ |\lambda|}\cdot\prod_{\Box\in \lambda}(1-(-q)^{h(\Box)})^{-2},
\end{align*}
where the sum is over all reversed plane partition $\mathbf{n}_\lambda$.
\end{corollary}
We devote the remainder of this section to prove the technical lemmas we used in Theorem \ref{thm: g=0 case}.

\begin{lemma}\label{lemma:rank Omega}
Let $C$ be a smooth projective curve of genus $g$,  $L_1,L_2$ be line bundles on $C$ such that $L_1\otimes L_2=K_{C} $ and let $\Omega\in K^0_\TT( {C}^{[\mathbf{n}_\lambda]})$ be the $K$-theory class of Remark \ref{rem: remark on Omega explicit}. Then
\begin{align*}
   \rk \left( T_{{C}^{[\mathbf{n}_\lambda]}}^{\vir}+ \Omega\right)=\deg L_1(c_\lambda+|\lambda|)+ (1-g)|\lambda|+|\mathbf{n}_\lambda| \mod 2.
\end{align*}
\end{lemma}
\begin{proof}
Let $\underline{Z}\in {C}^{[\mathbf{n}_\lambda]}$. Then
\begin{align*}
     \rk T_{{C}^{[\mathbf{n}_\lambda]}}^{\vir}&= \rk T_{C^{[\mathbf{n}_\lambda]}}^{\vir}|_{\underline{Z}}\\
     &= n_{00}+ \sum_{\substack{(i,j)\in \lambda\\ i\geq 1}}(n_{ij}-n_{i-1,j})+\sum_{\substack{(i,j)\in \lambda\\ j\geq 1}}(n_{ij}-n_{i,j-1})-\sum_{\substack{(i,j)\in \lambda\\ i,j\geq 1}}(n_{ij}-n_{i-1,j-1})
     \end{align*}
     For any line bundle $L$ on $C$ and   $\Delta$  a $\BZ$-linear combination of the universal divisors $\CZ_{ij}$  on $C\times C^{[\mathbf{n}_\lambda]}$, Riemann-Roch yields
     \begin{align*}
    \rk \RR\pi_*(\oO(\Delta)\otimes L)
    &=\chi(C, L\otimes \oO(\Delta|_{\underline{Z}}) )\\
    &= \deg L+ \deg \Delta|_{\underline{Z}}+ 1-g.
\end{align*}
Since  $\Omega$ is a sum of $K$-theoretic classes of the form $\RR\pi_*(\oO(\Delta)\otimes L)$, for suitable $L, \Delta$, we have
\begin{multline*}
     \rk\Omega=\sum_{\substack{(i,j)\in \lambda\\ (i,j)\neq (0,0)}}( n_{ij}-i\cdot\deg L_1-j\cdot(2g-2-\deg L_1)+1-g)\\
     -\sum_{\substack{(i,j),(l,k)\in \lambda\\(i,j)\neq (l,k)\\(i,j)\neq (l+1,k+1)}}(n_{lk}-n_{ij}-(i-l)\deg L_1-(j-k)(2g-2-\deg L_1)+1-g )\\   +\sum_{\substack{(i,j),(l,k)\in \lambda\\(i,j)\neq (l-1,k)\\(i,j)\neq (l,k+1)}}(n_{lk}-n_{ij}-(i-l+1)\deg L_1-(j-k)(2g-2-\deg L_1)+1-g ).
 \end{multline*}
 Denote by $\equiv$ the congruence modulo 2. We have
 \begin{align*}
     n_{00}-\sum_{\substack{(i,j)\in \lambda\\ (i,j)\neq (0,0)}} n_{ij}\equiv|\mathbf{n}_\lambda|.
 \end{align*}
 Denote by $V, E,Q$ respectively the number of vertices, edges and squares of the graph associated to $\lambda$ as in Lemma \ref{lemma: vertices and edges of graph}. We have
 \begin{align*}
     \sum_{\substack{(i,j)\in \lambda\\ (i,j)\neq (0,0)}}1- \sum_{\substack{(i,j),(l,k)\in \lambda\\(i,j)\neq (l,k)\\(i,j)\neq (l+1,k+1)}}1+ \sum_{\substack{(i,j),(l,k)\in \lambda\\(i,j)\neq (l-1,k)\\(i,j)\neq (l,k+1)}}1&= V-1-(|\lambda|^2-V-Q)+(|\lambda|^2-E)\\
     &=|\lambda|,
 \end{align*}
 where in the last line we used Lemma \ref{lemma: vertices and edges of graph}. The coefficient of $\deg L_1$ in $\rk \Omega$ is
 \begin{align*}
      &\sum_{\substack{(i,j)\in \lambda\\ (i,j)\neq (0,0)}}(i-j)- \sum_{\substack{(i,j),(l,k)\in \lambda\\(i,j)\neq (l,k)\\(i,j)\neq (l+1,k+1)}}(i-j-l+k)+ \sum_{\substack{(i,j),(l,k)\in \lambda\\(i,j)\neq (l-1,k)\\(i,j)\neq (l,k+1)}}(i-j+l-k+1)\\
      &\equiv \sum_{\substack{(i,j)\in \lambda}}(i+j)+ \sum_{\substack{(i,j),(l,k)\in \lambda\\(i,j)= (l,k)\\(i,j)= (l+1,k+1)}}(i+j+l+k)+ \sum_{\substack{(i,j),(l,k)\in \lambda\\(i,j)= (l-1,k)\\(i,j)= (l,k+1)}}(i+j+l+k)+\sum_{\substack{(i,j),(l,k)\in \lambda\\(i,j)\neq (l-1,k)\\(i,j)\neq (l,k+1)}} 1\\
      &\equiv c_\lambda+ 2E+ |\lambda|\\
      &\equiv c_\lambda+|\lambda|.
 \end{align*}
 Finally
 \begin{multline*}
     \sum_{\substack{(i,j),(l,k)\in \lambda\\(i,j)\neq (l,k)\\(i,j)\neq (l+1,k+1)}}(n_{lk}-n_{ij})+\sum_{\substack{(i,j),(l,k)\in \lambda\\(i,j)\neq (l-1,k)\\(i,j)\neq (l,k+1)}}(n_{lk}-n_{ij})\\+\sum_{\substack{(i,j)\in \lambda\\ i\geq 1}}(n_{ij}-n_{i-1,j})+\sum_{\substack{(i,j)\in \lambda\\ j\geq 1}}(n_{ij}-n_{i,j-1})-\sum_{\substack{(i,j)\in \lambda\\ i,j\geq 1}}(n_{ij}-n_{i-1,j-1})\equiv 0.
 \end{multline*}
 Combining all these identities together, we conclude that 
 \begin{align*}
   \rk \left( T_{{C}^{[\mathbf{n}_\lambda]}}^{\vir}+ \Omega\right)=\deg L_1(c_\lambda+|\lambda|)+ |\lambda|(1-g)+|\mathbf{n}_\lambda| \mod 2.
\end{align*}
\end{proof}

\begin{lemma}\label{lemma: no weight spaces on Omega}
Let $\underline{Z}\in  {\BP^{1}}^{[\mathbf{n}_\lambda], \BC^*}$ be a $\BC^*$-fixed point and set 
\begin{align*}
    \mathsf{V}_{\underline{Z}}= T_{{\BP^{1}}^{[\mathbf{n}_\lambda]}}^{\vir}|_{\underline{Z}}+ \Omega|_{\underline{Z}}\in K^0_{\TT\times \BC^*}(\pt),
\end{align*}
where $\Omega\in K^0_\TT( {\BP^{1}}^{[\mathbf{n}_\lambda]})$ is as in Remark \ref{rem: remark on Omega explicit}. Then we have
\begin{align*}
    \left.\left(\frac{e^{\TT\times \BC^*}(\mathsf{V}_{\underline{Z}}^\vee\otimes\tf_1\tf_2) }{e^{\TT\times \BC^*}(\mathsf{V}_{\underline{Z}})}\right)\right|_{s_1+s_2=0}=(-1)^{\rk \mathsf{V}_{\underline{Z}}}.
\end{align*}
\end{lemma}
\begin{proof}
Denote by $\mathsf{V}_{\underline{Z}}^{CY}$ the sub-representation of $\mathsf{V}_{\underline{Z}} $ consisting of weight spaces corresponding to the characters $(\tf_1\tf_2)^{a}$, for all $a\in \BZ$, with respect to the $\TT\times \BC^*$-action. We claim that $\mathsf{V}_{\underline{Z}}^{CY}$ is of the form
\begin{align*}
    \mathsf{V}_{\underline{Z}}^{CY}=A_{\underline{Z}}+A_{\underline{Z}}^\vee\otimes\tf_1\tf_2,
\end{align*}
for a suitable $A_{\underline{Z}}\in K^0_{\TT\times \BC^*}(\pt)$. \\
{\bf Step I:} Assuming the claim, set
\[\Tilde{\mathsf{V}}_{\underline{Z}}= \mathsf{V}_{\underline{Z}}- \mathsf{V}_{\underline{Z}}^{CY} \]
and  $ \Tilde{\mathsf{V}}_{\underline{Z}}=\sum_{\mu}\tf^{\mu}-\sum_{\nu}\tf^{\nu}\in K^0_{\TT\times \BC^*}(\pt)$, where none of the characters $\tf^\mu, \tf^\nu $ is a power of $\tf_1\tf_2$. Write 
$
e^{\TT\times \BC^*}(\tf^\mu)=\mu\cdot s,
$
where $s=(s_1,s_2,s_3)$. Then we conclude that
\begin{align*}
     \left.\left(\frac{e^{\TT\times \BC^*}(\mathsf{V}_{\underline{Z}}^\vee\otimes\tf_1\tf_2) }{e^{\TT\times \BC^*}(\mathsf{V}_{\underline{Z}})}\right)\right|_{s_1+s_2=0}&= \left.\left(\frac{e^{\TT\times \BC^*}(\Tilde{\mathsf{V}}_{\underline{Z}}^\vee\otimes\tf_1\tf_2) }{e^{\TT\times \BC^*}(\Tilde{\mathsf{V}}_{\underline{Z}})}\right)\right|_{s_1+s_2=0}\\
     &=\left.\left(\prod_{\mu}\frac{-\mu\cdot s+s_1+s_2}{\mu\cdot s}\cdot \prod_{\nu}\frac{\nu\cdot s}{-\nu\cdot s+s_1+s_2}\right)\right|_{s_1+s_2=0}\\
     &=(-1)^{\rk \mathsf{V}_{\underline{Z}}},
\end{align*}
where we used that no $\mu\cdot s, \nu\cdot s$ is a multiple of $ s_1+s_2$ and that $\rk \mathsf{V}_{\underline{Z}}=\rk \Tilde{\mathsf{V}}_{\underline{Z}} \mod 2$.\\
{\bf Step II:} We prove now our claim on $ \mathsf{V}_{\underline{Z}}^{CY}$. Firstly, by Proposition \ref{prop: C^* pot is movable} $  T_{{\BP^{1}}^{[\mathbf{n}_\lambda]}, \underline{Z}}^{\vir}$ is $\BC^*$-movable, which implies that there are no weight spaces corresponding to a power of $\tf_1\tf_2$.\\
It is clear that the $\TT\times \BC^*$-weight spaces of $\Omega|_{\underline{Z}}$ relative to the characters $(\tf_1\tf_2)^a$ are given by $\Omega'|_{\underline{Z}}$, where 
 \begin{multline*}
     \Omega'=\sum_{\substack{(a,a)\in \lambda\\ a\neq 0}}\RR\pi_*\left(\oO_{{\BP^1}\times {\BP^1}^{[\mathbf{n}_\lambda]}}(\CZ_{aa}) \otimes K_{\BP^1}^{-a}\right) (\tf_1\tf_2)^{-a}
     -\sum_{\substack{(i,j),(l,k)\in \lambda\\(i,j)=(l+a,k+a)\\a\neq 0,1}}\RR\pi_*\left(\oO(\Delta_{ij;lk})\otimes K_{\BP^1}^{a} \right)(\tf_1\tf_2)^{a}\\   +\sum_{\substack{(i,j),(l,k)\in \lambda\\(i,j)= (l+a-1,k+a)\\a\neq 0,1}}\RR\pi_* \left( \oO(\Delta_{ij;lk})\otimes K_{\BP^1}^a \right)(\tf_1\tf_2)^{a},
 \end{multline*}
 as we  just considered the weight spaces of $(\tf_1\tf_2)^a$ in $\Omega$ of Remark \ref{rem: remark on Omega explicit}. We notice that $\Omega'|_{\underline{Z}}$ is a sum of $\TT\times\BC^*$-representations of the form
 \begin{align*}
     &\RR\Gamma(\oO_{\BP^1}(Z_a)\otimes K_{\BP^1}^{-a})\otimes (\tf_1\tf_2)^{-a}, \mbox{ if } a\geq 1, \mbox{ or}\\
     &\RR\Gamma(\oO_{\BP^1}(-Z_a)\otimes K_{\BP^1}^{a})\otimes (\tf_1\tf_2)^a, \mbox{ if } a\geq 2,
 \end{align*}
 where $Z_a\subset \BP^1$ are effective divisors. We have the following identities of $\BC^*$-representations
 \begin{align*}
     \RR\Gamma(\oO_{\BP^1}(Z_a)\otimes K_{\BP^1}^{-a})&=\RR\Gamma(K_{\BP^1}^{-a})-\RR\Hom(\oO_{Z_a}, K_{\BP^1}^{-a}), \quad a\geq 1,\\
      \RR\Gamma(\oO_{\BP^1}(-Z_a)\otimes K_{\BP^1}^{a})&=\RR\Gamma(K_{\BP^1}^{a})- \RR\Gamma(\oO_{Z_a}\otimes K_{\BP^1}^{a}), \quad a\geq 2.
 \end{align*}
 By Lemma \ref{lemma: identities representations}, their $\BC^*$-fixed part is
 \begin{align*}
     (\RR\Gamma(K_{\BP^1}^{a}))^{\fix}&=\begin{cases}1 & a\leq -1,\\ -1 &a \geq 2,\end{cases}\\
     (\RR\Hom(\oO_{Z_a}, K_{\BP^1}^{-a}))^{\fix}&=0, \quad a\geq 1,\\
     (\RR\Gamma(\oO_{Z_a}\otimes K_{\BP^1}^{a}))^{\fix}&=0, \quad a\geq 2.
 \end{align*}
 Combining everything together, we conclude that
 \begin{align*}
     \mathsf{V}_{\underline{Z}}^{CY}=\sum_{\substack{(a,a)\in \lambda\\ a\neq 0}}(\tf_1\tf_2)^{-a}
     +\sum_{\substack{(i,j),(l,k)\in \lambda\\(i,j)=(l+a,k+a)\\a\neq 0,1}}\sgn(a)(\tf_1\tf_2)^{a}  - \sum_{\substack{(i,j),(l,k)\in \lambda\\(i,j)= (l+a-1,k+a)\\a\neq 0,1}}\sgn(a) (\tf_1\tf_2)^{a},
 \end{align*}
 where $\sgn$ is the usual sign function. Our claim follows by Lemma \ref{lemma: induction on poly for no powers CY}.
\end{proof}

\subsection{Case II: trivial   vector bundle}
If $L_1=L_2=\oO_{\BP^1}$, the integral  \eqref{eqn: integral on double nested} amounts to the leading term computation of Section \ref{sec: leading term} and a vanishing result, which relies on a combinatorial fact about the topological vertex formalism for stable pairs developed in  \cite{PT_vertex}.
\begin{prop}\label{prop: trivial case full with vanishing}
Under the anti-diagonal restriction $s_1+s_2=0$ we have an identity
\begin{align*}
    \sum_{\mathbf{n}_\lambda}q^{|\mathbf{n}_\lambda|} \left. \int_{[{\BP^1}^{[\mathbf{n}_\lambda]}]^{\vir}}e^{\TT}(-N_{\BP^1,\oO,\oO}^{\vir})\right|_{s_1+s_2=0}=(-s_1^2)^{-|\lambda|}\cdot\prod_{\Box\in \lambda}h(\Box)^{-2}. 
\end{align*}
\end{prop}
\begin{proof}
The leading term is computed in Proposition \ref{prop: leading term}, therefore we just need to prove that the integral vanishes for $ |\mathbf{n}_{\lambda}|>0$. We apply Graber-Pandharipande virtual localization \cite{GP_virtual_localization} with respect to the $\BC^*$-action on $ {\BP^1}^{[\mathbf{n}_\lambda]}$
\begin{align*}
    \int_{[{\BP^1}^{[\mathbf{n}_\lambda]}]^{\vir}}e^{\TT}(-N_{\BP^1,\oO,\oO}^{\vir})=\left.\left(\sum_{\underline{Z}\in {\BP^1}^{[\mathbf{n}_\lambda],\BC^*} } e^{\TT\times \BC^*}(-T^{\vir}_{{\BP^1}^{[\mathbf{n}_\lambda]},\underline{Z}}-N_{\BP^1,\oO,\oO, \underline{Z}}^{\vir})\right)\right|_{s_3=0},
\end{align*}
where $s_3$ is the generator of $H^*_{\BC^*}(\pt)$. The $\BC^*$-action on $ {\BP^1}^{[\mathbf{n}_\lambda]}$ is just the restriction of the natural $ \TT\times \BC^*$-action on  $P_n(\BP^1\times \BC^2, d[\BP^1])$, where $d=|\lambda|$ and $n=|\lambda|+|\mathbf{n}_\lambda|$. This means that, as $\TT\times \BC^*$-representation,
\begin{align*}
    \BE^\vee|_{\underline{Z}}=T^{\vir}_{{\BP^1}^{[\mathbf{n}_\lambda]},\underline{Z}}+N_{\BP^1,\oO,\oO, \underline{Z}}^{\vir}\in K^0_{\TT\times \BC^*}(\pt)
\end{align*}
can be described via the topological vertex formalism of Pandharipande-Thomas \cite[Thm. 3]{PT_vertex}, which states
\begin{align*}
    \BE^\vee|_{\underline{Z}}=\mathsf{V}_{0,\underline{Z}}+\mathsf{V}_{\infty,\underline{Z}}+\mathsf{E_{\underline{Z}}}\in K^0_{\TT\times \BC^*}(\pt),
\end{align*}
where $\mathsf{V}_{0,\underline{Z}},\mathsf{V}_{\infty,\underline{Z}}$ are the vertex terms corresponding to the two toric charts of $\BP^1\times \BC^2$ and $ \mathsf{E_{\underline{Z}}}$ is the edge term. By \cite[Lemma 22]{MPT_curves_K3} we have that  $e^{\TT\times \BC^*}(-\mathsf{V}_{0,\underline{Z}}-\mathsf{V}_{\infty,\underline{Z}})$ is divisible by $(s_1+s_2)$ if $|\mathbf{n}_\lambda|>0$, while $e^{\TT\times \BC^*}(-\mathsf{E_{\underline{Z}}})$ is easily seen to be coprime with $(s_1+s_2)$.\footnote{In fact, $\mathsf{E_{\underline{Z}}}$ is the $\TT$-representation of the tangent space at a $\TT$-fixed point of $\Hilb^{|\lambda|}(\BC^2)$. The claim follows by noticing that the $ \TT$-fixed locus coincides with the fixed locus of the subtorus $\set{t_1t_2=1}\subset \TT$ preserving the Calabi-Yau form of $\BC^2$, which is 0-dimensional and reduced.} This implies that the anti-diagonal restriction $s_1+s_2=0$ is well-defined on every localized term and satisfies
\begin{align*}
        \left.e^{\TT\times\BC^*}(\BE^\vee|_{\underline{Z}})\right|_{s_1+s_2=0}=0,
\end{align*}
by which we conclude the required vanishing.
\end{proof}
\begin{remark}
We could prove Proposition \ref{prop: trivial case full with vanishing} without relying on the vertex formalism for stable pairs, by simply carrying out a detailed (but probably longer) analysis of the weight space of $\tf_1\tf_2$ as in Lemma \ref{lemma: no weight spaces on Omega}.
\end{remark}

\section{Summing  the theory up}
\subsection{Anti-diagonal restriction}
We combine the computations in Section \ref{sec: leading term}, \ref{sec: toric computations} to prove the second part of Theorem \ref{thm: main result intro} from the Introduction.
\begin{theorem}\label{thm: generating series with antidiagonal restriction}
 Under the anti-diagonal restriction $s_1+s_2=0$ the three universal series are
\begin{align*}
   A_{\lambda}(q,s_1,-s_1)&=(-s_1^2)^{|\lambda|}\cdot\prod_{\Box\in \lambda}h(\Box)^2,\\
    B_{\lambda}(-q,s_1,-s_1)&=(-1)^{n(\lambda)}\cdot s_1^{-|\lambda|}\cdot\prod_{\Box\in \lambda}h(\Box)^{-1}\cdot\prod_{\Box\in \lambda}(1-q^{h(\Box)}),\\
    C_{\lambda}(-q,s_1,-s_1)&=(-1)^{n(\overline{\lambda})}\cdot (-s_1)^{-|\lambda|}\cdot\prod_{\Box\in \lambda}h(\Box)^{-1}\cdot\prod_{\Box\in \lambda}(1-q^{h(\Box)}).
\end{align*}
\end{theorem}
\begin{proof}
Let $C_{g,\deg L_1\deg L_2}\in \BQ(s_1,s_2)$ be the leading term of the generating series of Theorem \ref{thm: universal series}. We can write
\begin{multline*}
   C_{g,\deg L_1\deg L_2}^{-1}\cdot \sum_{\mathbf{n}_\lambda}q^{|\mathbf{n}_\lambda|} \int_{[C^{[\mathbf{n}_\lambda]}]^{\vir}}e^{\TT}(-N_{C,L_1,L_2}^{\vir})=\Tilde{A}_{\lambda,1}^{g-1}\cdot \Tilde{A}_{\lambda,2}^{\deg L_1}\cdot \Tilde{A}_{\lambda,3}^{\deg L_2}\\
   =\exp\left((g-1)\cdot\log\Tilde{A}_{\lambda,1}+ \deg L_1\cdot \log\Tilde{A}_{\lambda,2}+\deg L_2\cdot \log\Tilde{A}_{\lambda,3}\right),
\end{multline*}
for suitable $\Tilde{A}_{\lambda,i}\in 1+\BQ(s_1,s_2)\llbracket q \rrbracket$, where $\log \Tilde{A_i} $ are  well-defined  as $\Tilde{A_i}$ are power series starting with 1. The claim is therefore reduced to the computation of the leading term (under the anti-diagonal restriction $s_1+s_2=0$) and to the solution of a linear system. The leading term is computed in Proposition \ref{prop: leading term}, which also shows that it is invertible in $\BQ(s_1,s_2)$.
The linear system  is solved by computing the generating series of the integrals  \eqref{eqn: integral on double nested}, under the anti-diagonal restriction, on a basis of $\BQ^3$. The classes
 \begin{align*}
     \gamma(\BP^1, \oO, \oO), \quad \gamma(\BP^1, \oO(-1), \oO(-1)), \quad
\gamma(\BP^1, \oO, \oO(-2))
\end{align*}
are linearly independent and we computed their generating series in Proposition \ref{prop: trivial case full with vanishing} and Corollary \ref{cor: gen series CY}.
\end{proof}
\subsection{Degree 1}\label{sec: full equivariant}
In degree 1, we consider the  Young diagram consisting of a  single box and the corresponding double nested Hilbert scheme is $C^{[n]}\cong C^{(n)}$, the symmetric power of a smooth projective curve $C$, with universal subscheme $\CZ\subset C\times C^{(n)}$. 
Given line bundles $L_1, L_2$ on $C$,  the class \eqref{eqn: dual pot} in $K$-theory of the virtual normal bundle is
\begin{align*}
    N^{\vir}_{C, L_1, L_2}&=-\RR\pi_*(L_1L_2\otimes \oO_\CZ)\otimes \tf_1\tf_2+ \RR\pi_*L_1\otimes \tf_1+\RR\pi_*L_2\otimes \tf_2\\
    &=-(L_1L_2)^{[n]}\otimes \tf_1\tf_2+ \oO_{C^{(n)}}^{\deg L_1+ 1-g}\otimes \tf_1+\oO_{C^{(n)}}^{\deg L_2+ 1-g}\otimes \tf_2,
\end{align*}
 where $(L_1L_2)^{[n]} $ is the tautological bundle with fibers $(L_1L_2)^{[n]}|_{Z}=H^0(C,L_1L_2\otimes \oO_Z)$.  This yields
\begin{align*}
    \int_{[C^{[\mathbf{n}_\lambda]}]^{\vir}}e^\TT(- N^{\vir}_{C, L_1, L_2})
    &= s_1^{g-1-\deg L_1} s_2^{g-1-\deg L_2}\int_{C^{(n)}}e((L_1L_2)^{[n]}).
\end{align*}
By the universal structure of Theorem \ref{thm: universal series}, we just need to compute the (generating series of the) last integral   for $L_1=L_2=\oO_{C}$ and $L_1\otimes L_2\cong K_{C}$, which yields the explicit universal series
\begin{align*}
   A_{\Box}(q,s_1,s_2)&=s_1 s_2,\\
    B_{\Box}(q,s_1,s_2)&= s_1^{-1}(1+q),\\
    C_{\Box}(q,s_1,s_2)&=s_2^{-1}(1+q).
\end{align*}
\subsection{GW/PT correspondence}
Let $X=\Tot_C(L_1\oplus L_2)$ be a local curve. Combining Theorem \ref{thm: universal series}, Theorem \ref{thm: generating series with antidiagonal restriction} and the description of the $\TT$-fixed locus of $P_n(X, d[C])$ in Proposition \ref{prop: iso of schemes fixed locus}, we obtain our main result.
\begin{theorem}\label{thm: full PT invariants}
The generating series of stable pair invariants satisfies
\begin{align*}
\PT_d(X;q)=\sum_{\lambda\vdash d}\left(q^{-|\lambda|}A_{\lambda}(q)\right)^{g-1}\cdot\left(q^{-n(\lambda)} B_{\lambda}(q)\right)^{\deg L_1}\cdot\left(q^{-n(\overline{\lambda})} C_{\lambda}(q)\right)^{\deg L_2},
\end{align*}
where $A_{\lambda,i}$ are the universal series of Theorem \ref{thm: universal series}. Moreover, under the anti-diagonal restriction $s_1+s_2=0$
\begin{multline*}
      \left.  \PT_d(X;-q)\right|_{s_1+s_2=0}=\\
        (-1)^{d\deg L_2}\sum_{\lambda\vdash d}q^{d(1-g)-\deg L_1 n(\lambda)-\deg L_2 n(\overline{\lambda})}\prod_{\Box\in \lambda}(s_1h(\Box))^{2g-2-\deg L_1- \deg L_2}(1-q^{h(\Box)})^{\deg L_1+\deg L_2}.
\end{multline*}
\end{theorem}
Comparing with Bryan-Pandharipande's results --- cf. Theorem \ref{thm: GW series} and  \cite[Sec. 8]{BP_local_GW_curves} for the fully equivariant result in degree 1 --- we obtain a proof of the Gromov-Witten/stable pairs correspondence for local curves.
\begin{corollary}\label{cor: GW/PT correspondence}
Let $X$ be a local curve. Under the anti-diagonal restriction $s_1+s_2=0$ the GW/stable pair correspondence holds
\begin{align*}
    (-i)^{d(2-2g+\deg L_1+\deg L_2)}\cdot \GW_d(g|\deg L_1, \deg L_2;u)=(-q)^{-\frac{1}{2}\cdot d(2-2g+\deg L_1+\deg L_2)}\PT_d(X, q),
\end{align*}
after the change of variable $q=-e^{iu}$. Moreover it holds fully equivariantly in degree 1.
\end{corollary}
\subsection{Resolved conifold}
In some cases, the moduli space $P_n(X, d[C])$ happens to be proper, for example whenever $H^0(C, L_i)=0$ for $i=1,2$; see \cite[Prop. 3.1]{CKM_Stable_Pairs} for a similar setting for local surfaces. Under this assumption, the invariants are computed as
\begin{align*}
    \PT_d(X;q)=\left.  \PT_d(X;q)\right|_{s_1=s_2=0}
\end{align*}
and can be deduced from the anti-diagonal restriction. The resulting invariants are interesting - that is, non-zero - only in the Calabi-Yau case; in this case, they coincide with the virtual Euler characteristic and Behrend's weighted Euler characteristic of $P_n(X, d[C])$ \cite{Beh_DT_via_microlocal}. Applying Riemann-Roch, this situation may appear only when $H^l(C,L_i)=0$ for $l=0,1$ and $\deg L_i=g-1$ for  $i=1,2$. \\

An interesting example is the resolved conifold $X=\Tot_{\BP^1}(\oO(-1)\oplus \oO(-1))$; in this case, the invariants can be further packaged into a generating series
\begin{align*}
   1+ \sum_{d\geq 1}Q^d\cdot\PT_d(X;-q)&=\sum_{d\geq 0}(-Qq)^{d}\sum_{\lambda\vdash d}q^{ n(\lambda)+ n(\overline{\lambda})}\prod_{\Box\in \lambda}(1-q^{h(\Box)})^{-2}\\
   &= \sum_{\lambda}(-Qq)^{|\lambda|}\frac{q^{n(\lambda)}}{\prod_{\Box\in \lambda}(1-q^{h(\Box)})}\frac{q^{n(\overline{\lambda})}}{\prod_{\Box\in \overline{\lambda}}(1-q^{h(\Box)})}\\
   &= \sum_{\lambda}(-Qq)^{|\lambda|} s_{\lambda}(q)s_{\overline{\lambda}}(q)\\
   &=\prod_{ n\geq 1}(1-Qq^n)^n,
\end{align*}
where we used some identities involving the Schur function $s_\lambda$ (see e.g. \cite{MacD_symmetric_functions_Hall}). This last generating series agrees with the expression of the unrefined limit of the topological vertex of Iqbal-Kozçaz-Vafa \cite{IKV_topological_vertex} and can be seen as a specialization both of the motivic invariants of Morrison-Mozgovoy-Nagao-Szendrői \cite{MMNS_motivic_DT_conifold} and of the $K$-theoretic invariants of Kononov-Okounkov-Osinenko \cite{KOO_2_legDT}.

\section{$K$-theoretic refinement}\label{sec: K-th}
\subsection{$K$-theoretic invariants}
Let $X=\Tot_C(L_1\oplus L_2)$ be a local curve. The perfect obstruction theory on $P_X:=P_n(X, d[C])$ induces a ($\TT$-equivariant) virtual structure sheaf $\oO^{\vir}_{P_X}\in K^\TT_0(P_X)$ \cite{FG_riemann_roch} which depends only on the $K$-theory class of the perfect obstruction theory \cite[Cor. 4.5]{Tho_K-theo_Fulton}. $K$-theoretic PT invariants are defined by  virtual $K$-theoretic localization \cite{FG_riemann_roch} for any $V\in K_0^\TT(P_X)$
\begin{align*}
      \PT^{K}_{d,n}(X, V)&:=\chi(P_n(X, d[C]), \oO_{P_X}^{\vir}\otimes V)\\
      &:=\chi\left(P_n(X, d[C])^{\TT}, \frac{\oO_{P_X^\TT}^{\vir}\otimes V|_{P_X^\TT}}{\Lambda^\bullet N^{\vir, *}} \right)\in \BQ(\tf_1, \tf_2),
\end{align*}
where $\Lambda^\bullet(V):=\sum_{i=0}^{\rk V}(-1)^i \Lambda^i V$ is defined for every locally free sheaf $V$ and then extended  by linearity to any class in $K$-theory.
\begin{remark}
Differently than the case of equivariant cohomology,  $\chi(M, \CF)\in K_0^\TT(\pt)_{\mathrm{loc}}$  is well-defined for any $\TT$-equivariant coherent sheaf $\CF$ on a non-proper scheme $M$, as long as the weight spaces of $\CF$
are finite-dimensional; in this case, the virtual localization formula is an actual theorem, rather than an ad-hoc definition of the invariants.
\end{remark}

Using the description of the $\TT$-fixed locus $(P_n(X, d[C))^\TT$ of Proposition \ref{prop: iso of schemes fixed locus}, $K$-theoretic stable pair invariants on $X$ (with no insertions) are reduced to intersection numbers on $C^{[\mathbf{n}_\lambda]}$. The same techniques of Section \ref{sec: universal} can be applied in this setting, yielding the following result.
\begin{prop}\label{prop: universality K-theory}
Let $C$ be a genus $g$ smooth irreducible projective curve and $L_1, L_2$  line bundles over $C$. We have an identity
\begin{align*}
    \sum_{\mathbf{n}_\lambda}q^{|\mathbf{n}_\lambda|} \chi\left(C^{[\mathbf{n}_\lambda]}, \frac{\oO_{C^{[\mathbf{n}_\lambda]}}^{\vir}}{\Lambda^\bullet N_{C,L_1,L_2}^{\vir, *}} \right)=A_{K,\lambda}^{g-1}\cdot B_{K,\lambda}^{\deg L_1}\cdot C_{K,\lambda}^{\deg L_2}\in \BQ(\tf_1,\tf_2)\llbracket q \rrbracket,
\end{align*}
where $A_{K,\lambda},B_{K,\lambda},C_{K,\lambda}\in \BQ(\tf_1,\tf_2)\llbracket q \rrbracket$ are fixed universal series for $i=1,2,3$, only depending on $\lambda$. Moreover
\begin{align*}
 A_{K,\lambda}(\tf_1,\tf_2)&=A_{K,\overline{\lambda}}(\tf_2,\tf_1),\\
    B_{K,\lambda}(\tf_1,\tf_2)&=C_{K,\overline{\lambda}}(\tf_2,\tf_1).
\end{align*}
\end{prop}
\begin{proof}
The proof follows the same strategy as Theorem \ref{thm: universal series}. We just need to notice that $\Lambda^\bullet(\cdot)$ is multiplicative and that, via virtual Hirzebruch-Riemann-Roch \cite{FG_riemann_roch}, we can express
\begin{align*}
     \chi\left(C^{[\mathbf{n}_\lambda]}, \frac{\oO_{C^{[\mathbf{n}_\lambda]}}^{\vir}}{\Lambda^\bullet N_{C,L_1,L_2}^{\vir, *}} \right)&=\int_{[C^{[\mathbf{n}_\lambda]}]^{\vir}}\ch(-\Lambda^\bullet N_{C,L_1,L_2}^{\vir, *})\cdot \td(T^{\vir}_{C^{[\mathbf{n}_\lambda]}})\\
     &=\int_{A_{C, \mathbf{n}_\lambda}}f,
\end{align*}
where $f$ is  a polynomial expression of    classes of the same form as in Proposition \ref{prop: taut dependence cobordism class with Pic}.
\end{proof}
Denote by $\PT^K_d(X;q)=\sum_{n\in \BZ}q^n\PT^{K}_{d,n}(X) $ the generating series of $K$-theoretic stable pair invariants.
\begin{corollary}
Let $X=\Tot_C(L_1\oplus L_2)$ be a local curve. We have
\begin{align*}
\PT^K_d(X;q)=\sum_{\lambda\vdash d}\left(q^{-|\lambda|}A_{K,\lambda}(q)\right)^{g-1}\cdot\left(q^{-n(\lambda)} B_{K,\lambda}(q)\right)^{\deg L_1}\cdot\left(q^{-n(\overline{\lambda})} C_{K,\lambda}(q)\right)^{\deg L_2}.
\end{align*}
\end{corollary}
\subsection{Nekrasov-Okounkov}
Let $M$ be a scheme with a perfect obstruction theory $\BE$.  Define the  virtual canonical bundle to be $K_{\vir}=\det \BE\in \Pic(M)$. Assume that $ K_{\vir}$ admits a \emph{square root} $K_{\vir}^{1/2}$, that is a line bundle such that $ (K_{\vir}^{1/2})^{\otimes 2}\cong K_{\vir}$. Nekrasov-Okounkov \cite{NO_membranes_and_sheaves} teach us that it is much more natural to consider the \emph{twisted virtual structure sheaf}
\begin{align*}
    \widehat{\oO}_{M}^{\vir}=\oO_{M}^{\vir}\otimes K_{\vir}^{1/2}.
\end{align*}
In \cite{NO_membranes_and_sheaves}, Nekrasov-Okounkov show existence of square roots for $P_X$  (and uniqueness, up to 2-torsion). Nevertheless, even if square roots could not exists on $P_X$ as \emph{line bundles}, they exist as a class $ K_{\vir, P_X}^{1/2}\in K^0(P_X, \BZ[\frac{1}{2}])$ and are unique \cite[Lemma 5.1]{OT_1}. In our setting, we define \emph{ Nekrasov-Okounkov $K$-theoretic stable pair invariants} as\footnote{To take square roots \emph{equivariantly}, we need to replace the torus $\TT$ with the minimal  cover where the characters $\tf_1^{1/2},\tf_2^{1/2}$ are well-defined, see \cite[Sec. 7.2.1]{NO_membranes_and_sheaves}.}
\begin{align*}
      \PT^{\widehat{K}}_{d,n}(X):&= \PT^{K}_{d,n}(X, K_{\vir, P_X}^{1/2})\\
      &=\chi\left(P_n(X, d[C])^{\TT}, \frac{\oO_{P_X^\TT}^{\vir}\otimes K_{\vir, P_X}^{1/2}|_{P_X^\TT}}{\Lambda^\bullet N^{\vir, *}} \right)\in \BQ(\tf_1^{1/2}, \tf_2^{1/2}),
\end{align*}
which are an algebro-geometric analogue of the $\widehat{A}$-genus of a spin manifold. On the $\TT$-fixed locus, we have an identity in $K^0_\TT(P_X^\TT, \BZ[\frac{1}{2}])_{\mathrm{loc}} $
\begin{align*}
     \frac{\oO_{P_X^\TT}^{\vir}\otimes K_{\vir, P_X}^{1/2}|_{P_X^\TT}}{\Lambda^\bullet N^{\vir, *}} &= \frac{\oO_{P_X^\TT}^{\vir}\otimes  K_{\vir, P_X^\TT}^{1/2}\otimes  \left(\det N^{\vir, *}\right)^{1/2}}{{\Lambda^\bullet N^{\vir, *}}}\\
     &=  \frac{\widehat{\oO}_{P_X^\TT}^{\vir}}{\widehat{\Lambda}^\bullet N^{\vir, *}},
\end{align*}
where we define $\widehat{\Lambda}^\bullet(\cdot)=\Lambda^\bullet(\cdot)\otimes \det(\cdot)^{-1/2}$. Again, the same techniques of Section \ref{sec: universal} and Proposition \ref{prop: universality K-theory} and  can be applied in this setting, yielding the following result.
\begin{theorem}
Let $C$ be a genus $g$ smooth irreducible projective curve and $L_1, L_2$  line bundles over $C$. We have an identity
\begin{align*}
    \sum_{\mathbf{n}_\lambda}q^{|\mathbf{n}_\lambda|} \chi\left(C^{[\mathbf{n}_\lambda]},  \frac{\widehat{\oO}_{P_X^\TT}^{\vir}}{\widehat{\Lambda}^\bullet N_{C, L_1, L_2}^{\vir, *}} \right)=A_{\widehat{K},\lambda}^{g-1}\cdot B_{\widehat{K},\lambda}^{\deg L_1}\cdot C_{\widehat{K},\lambda}^{\deg L_2}\in \BQ(\tf_1^{1/2}, \tf_2^{1/2})\llbracket q \rrbracket,
\end{align*}
where $A_{\widehat{K},\lambda},B_{\widehat{K},\lambda},C_{\widehat{K},\lambda}\in \BQ(\tf_1^{1/2}, \tf_2^{1/2})\llbracket q \rrbracket$ are fixed universal series for $i=1,2,3$. Moreover
\begin{align*}
 A_{\widehat{K},\lambda}(\tf_1,\tf_2)&=A_{\widehat{K},\overline{\lambda}}(\tf_2,\tf_1),\\
    B_{\widehat{K},\lambda}(\tf_1,\tf_2)&=C_{\widehat{K},\overline{\lambda}}(\tf_2,\tf_1).
\end{align*}
\end{theorem}
Denote by $\PT^{\widehat{K}}_d(X;q)=\sum_{n\in \BZ}q^n\PT^{\widehat{K}}_{d,n}(X) $ the generating series of Nekrasov-Okounkov $K$-theoretic stable pair invariants.
\begin{corollary}\label{cor: K-theoretic NO universal}
Let $X=\Tot_C(L_1\oplus L_2)$ be a local curve. We have
\begin{align*}
\PT^{\widehat{K}}_d(X;q)=\sum_{\lambda\vdash d}\left(q^{-|\lambda|}A_{\widehat{K},\lambda}(q)\right)^{g-1}\cdot\left(q^{-n(\lambda)} B_{\widehat{K},\lambda}(q)\right)^{\deg L_1}\cdot\left(q^{-n(\overline{\lambda})} C_{\widehat{K},\lambda}(q)\right)^{\deg L_2}.
\end{align*}
\end{corollary}
The techniques of Section \ref{sec: leading term}, \ref{sec: toric computations} can be adapted to compute the generating series of Nekrasov-Okounkov $K$-theoretic invariants under the anti-diagonal restriction $\tf_1\tf_2=1$. In fact, as in the proof of Theorem  \ref{thm: generating series with antidiagonal restriction}, we just need to compute the leading term of the generating series and the cases $g=0$, $L_1\otimes L_2=K_{\BP^1}$ and $L_1=L_2=\oO_{\BP^1}$.  Similarly to the proof of Theorem \ref{thm: g=0 case}, as we work with $C\cong \BP^1$, applying the $K$-theoretic virtual localization formula \cite{FG_riemann_roch} on $(\BP^1)^{[\mathbf{n}_\lambda]}$  yields
\begin{align*}
    \chi\left({\BP^1}^{[\mathbf{n}_\lambda]},  \frac{\widehat{\oO}_{P_X^\TT}^{\vir}}{\widehat{\Lambda}^\bullet N_{C, L_1, L_2}^{\vir, *}} \right)&=\left.\left(\sum_{\underline{Z}\in {\BP^1}^{[\mathbf{n}_\lambda],\BC^*} }\chi\left(\underline{Z}, \frac{1 }{\widehat{\Lambda}^\bullet(T^{\vir, *}_{{\BP^1}^{[\mathbf{n}_\lambda]},\underline{Z}}+N_{\BP^1,L_1,L_2, \underline{Z}}^{\vir, *} )}\right)\right)\right|_{\tf_3=1},
\end{align*}
where $\tf_3$ is the equivariant parameter of the $\BC^*$-action. For a character $\tf^{\mu}\in K_{\TT\times \BC^*}^0(\pt)$, define the operator $[\tf^{\mu}]=\tf^{\frac{\mu}{2}}-\tf^{-\frac{\mu}{2}}$ and extend it by linearity to any  $V\in K_{\TT\times \BC^*}^0(\pt)$. It is proven in \cite[Sec. 6.1]{FMR_higher_rank} that
\begin{align*}
    \chi(\pt, \widehat{\Lambda}^\bullet(V^*))=[V],
\end{align*}
which satisfies $[V^*]=(-1)^{\rk V}[V]$, therefore
\begin{align*}
    \chi\left({\BP^1}^{[\mathbf{n}_\lambda]},  \frac{\widehat{\oO}_{P_X^\TT}^{\vir}}{\widehat{\Lambda}^\bullet N_{C, L_1, L_2}^{\vir, *}} \right)&=\left.\left(\sum_{\underline{Z}\in {\BP^1}^{[\mathbf{n}_\lambda],\BC^*} }[-T^{\vir}_{{\BP^1}^{[\mathbf{n}_\lambda]},\underline{Z}}-N_{\BP^1,L_1,L_2, \underline{Z}}^{\vir}]\right)\right|_{\tf_3=1}.
\end{align*}
Explicit computations yields
\begin{align*}
   A_{\widehat{K},\lambda}(q,\tf_1,\tf^{-1}_1)&=(-1)^{|\lambda|}\cdot F_\lambda^{-2},\\
    B_{\widehat{K}, \lambda}(-q,\tf_1,\tf^{-1}_1)&=(-1)^{n(\lambda)}\cdot F_\lambda\cdot\prod_{\Box\in \lambda}(1-q^{h(\Box)}),\\
    C_{\widehat{K},\lambda}(-q,\tf_1,\tf^{-1}_1)&=(-1)^{|\lambda|+n(\overline{\lambda})}F_\lambda\cdot\prod_{\Box\in \lambda}(1-q^{h(\Box)}),
\end{align*}
where 
\[
F_\lambda=\prod_{\substack{(i,j)\in \lambda\\(i,j)\neq (a,a)}}\frac{1}{\tf_1^{\frac{|j-i|}{2}}-\tf_1^{-\frac{|j-i|}{2}}}\cdot \frac{\prod_{\substack{(i,j), (l,k)\in \lambda\\ (i,j)\neq (l+a,k+a)}} \left(\tf_1^{\frac{|i-j+k-l|}{2}}-\tf_1^{-\frac{|i-j+k-l|}{2}}\right) }{\prod_{\substack{(i,j), (l,k)\in \lambda\\ (i,j)\neq (l+a-1,k+a)}}\left(\tf_1^{\frac{|1+i-j+k-l|}{2}}-\tf_1^{-\frac{|1+i-j+k-l|}{2}}\right)}.
\]
\appendix
\section{The combinatorial identities}
In this appendix we collect the proofs of some technical results on the combinatorics of  Young diagrams we have used.
\begin{lemma}\label{lemma: induction on poly for no powers CY}
Let $\lambda$ be a Young diagram and consider the Laurent polynomials in $\BZ[\tf^{\pm 1}]$ 
\begin{align*}
    g_{\lambda}(\tf)&=\sum_{\substack{(a,a)\in \lambda\\ a\neq 0, a\in \BZ}}\tf^{-a}
     +\sum_{\substack{(i,j),(l,k)\in \lambda\\(i,j)=(l+a,k+a)\\a\neq 0,1, a\in \BZ}}\sgn(a)\tf^{a}  - \sum_{\substack{(i,j),(l,k)\in \lambda\\(i,j)= (l+a-1,k+a)\\a\neq 0,1, a\in \BZ}}\sgn(a)\tf^{a},\\
     h_{\lambda}(\tf)&=\sum_{\substack{(a,a)\in \lambda\\ a\neq 0, a\in \BZ}}\tf^{-a}
     -\sum_{\substack{(i,j),(l,k)\in \lambda\\(i,j)=(l+a,k+a)\\a\neq 0,1, a\in \BZ}}\tf^{a}  + \sum_{\substack{(i,j),(l,k)\in \lambda\\(i,j)= (l+a-1,k+a)\\a\neq 0,1, a\in \BZ}}\tf^{a}.
\end{align*}
     Then we have
     \begin{align*}
             g_{\lambda}(\tf)&=A_\lambda(\tf)+A_\lambda(\tf^{-1})\tf,\\
             h_{\lambda}(\tf)&=B_\lambda(\tf)-B_\lambda(\tf^{-1})\tf,
     \end{align*}
     where $A_{\lambda}(\tf),B_{\lambda}(\tf)\in \BZ[\tf^{\pm 1}]$.
\end{lemma}
\begin{proof}
We prove the first claim by induction on the size of $\lambda$. If $|\lambda|=1 $ this is clear. Suppose now the claim holds for all Young diagrams of size $|\lambda|=n$ and consider a Young diagram of size $|\Tilde{\lambda}|=n+1$ obtained by adding to a Young diagram $\lambda$ a box whose lattice coordinates are $(i,j)\in \BZ^{2}$. 
\begin{itemize}
\item $(i,j)=(i,0)$ or $(i,j)=(0,j)$, with $i,j\neq 0$. We have
\begin{align*}
    g_{\Tilde{\lambda}}&=g_{\lambda}.
\end{align*}
    \item $(i,j)=(i,1), i\geq 1$. We have
    \begin{align*}
        g_{\Tilde{\lambda}}&=g_{\lambda}+\tf^{-1}-\tf^{-1}\\
        &=g_{\lambda}.
    \end{align*}
     \item $(i,j)=(1,j), j\geq 2$. We have
    \begin{align*}
         g_{\Tilde{\lambda}}&=g_{\lambda}-\tf^{-1}-\tf^{2}\\
        &=g_{\lambda}-\tf^{-1}-(\tf^{-1})^{-1}\tf.
    \end{align*}
    \item  $(i,j)=(i,i), i\geq 2$. We have
    \begin{align*}
         g_{\Tilde{\lambda}}&=g_{\lambda}+\tf^{-i}-\sum_{l=1}^i\tf^{-l}+\sum_{l=2}^i \tf^l-\left(\sum_{l=2}^i\tf^{l}-\sum_{l=1}^{i-1}\tf^{-l}\right)\\
        &=g_{\lambda}.
    \end{align*}
    \item  $(i,j), i> j\geq 2$. We have
    \begin{align*}
         g_{\Tilde{\lambda}}&=g_{\lambda}-\sum_{l=1}^j\tf^{-l}+\sum_{l=2}^j \tf^l-\left(\sum_{l=2}^j\tf^{l}-\sum_{l=1}^{j}\tf^{-l}\right)\\
        &=g_{\lambda}.
    \end{align*}
     \item  $(i,j), j> i\geq 2$. We have
      \begin{align*}
         g_{\Tilde{\lambda}}&=g_{\lambda}-\sum_{l=1}^i\tf^{-l}+\sum_{l=2}^i \tf^l-\left(\sum_{l=2}^{i+1}\tf^{l}-\sum_{l=1}^{i-1}\tf^{-l}\right)\\
        &=g_{\lambda}-\tf^{-i}-\tf^{i+1}.
    \end{align*}
\end{itemize}
Therefore the induction step is proven in all possible cases and we conclude the proof.\\
With an analogous analysis one proves the second claim as well. 
\end{proof}
\begin{lemma}\label{lemma: hook length in size 0}
Let $\lambda$ be a Young diagram. Then the following identity holds
\begin{align*}
    \prod_{\substack{(i,j)\in \lambda\\(i,j)\neq (a,a)}}(j-i)\cdot \frac{\prod_{\substack{(i,j), (l,k)\in \lambda\\ (i,j)\neq (l+a-1,k+a)}} (1+i-j+k-l)}{\prod_{\substack{(i,j), (l,k)\in \lambda\\ (i,j)\neq (l+a,k+a)}} (i-j+k-l)}=\sigma(\lambda)\cdot\prod_{\Box\in \lambda}h(\Box),
\end{align*}
where 
\begin{align*}
    \sigma(\lambda)= \prod_{\substack{(i,j)\in \lambda\\(i,j)\neq (a,a)}}\sgn(j-i)\cdot \prod_{\substack{(i,j), (l,k)\in \lambda\\ (i,j)\neq (l+a-1,k+a)}} \sgn(1+i-j+k-l)\cdot\prod_{\substack{(i,j), (l,k)\in \lambda\\ (i,j)\neq (l+a,k+a)}} \sgn(i-j+k-l).
\end{align*}
\end{lemma}
\begin{proof}
The sign $\sigma(\lambda)$ is easily determined, so we just need to compute the absolute value. To ease the notation, we adopt the following convention for the remainder of the proof: we set $|0|=1$, which is merely a formal shortcut to include in the productory  trivial factors we would have otherwise excluded. The claim therefore becomes
\begin{align}\label{eqn: claim hooks with convention}
    \prod_{(i,j)\in \lambda}|i-j|\cdot \prod_{(i,j), (l,k)\in \lambda}\frac{ |1+i-j+k-l|}{|i-j+k-l|}=\prod_{\Box\in \lambda}h(\Box).
\end{align}
Denote the left-hand-side of \eqref{eqn: claim hooks with convention} by $H_\lambda$. We prove this claim on the induction on the size of $\lambda$. If $|\lambda|=1$, the claim is trivially satisfied. Assume it holds for all Young diagrams of size $n$ and consider a Young diagram $\lambda'$ of size $n+1$ obtained by a Young diagram of size $\lambda$ by adding a box with lattice coordinates $(i,j)\in \BZ^2$. We have
\begin{align*}
    H_{\lambda'}&=H_\lambda\cdot |i-j|\cdot \prod_{(l,k)\in \lambda}\frac{ |1+i-j+k-l|}{|i-j+k-l|}\cdot\frac{ |-1+i-j+k-l|}{|i-j+k-l|}.
\end{align*}
To avoid confusion, we denote by $h(\Box)$ (resp. $h'(\Box)$) the hooklength of $\Box\in \lambda $ (resp. $\Box\in \lambda'$).   
The strategy now is to divide the boxes of $\lambda'$ in sub-collections and compute separately each contribution of the product on the right-hand-side.\\
{\bf Step I:} Fix a box $(\Tilde{i}, j)\in \lambda$, with $\Tilde{i}<i$. The contribution of all boxes on the right (on the same row) of $(\Tilde{i},j)$ is
\begin{align*}
    &\prod_{k=j}^{\lambda_{\Tilde{i}}-1}\frac{|i-j+k-\Tilde{i}+1|\cdot|i-j+k-\Tilde{i}-1|}{|i-j+k-\Tilde{i}|\cdot|i-j+k-\Tilde{i}|}\\
    &=\frac{|i-\Tilde{i}+1|\cdot|i-\Tilde{i}-1|}{|i-\Tilde{i}|\cdot|i-\Tilde{i}|}\cdot\frac{|i-\Tilde{i}+2|\cdot|i-\Tilde{i}|}{|i-\Tilde{i}+1|\cdot|i-\Tilde{i}+1|}\cdots\frac{|i-\Tilde{i}+\lambda_{\Tilde{i}}-j|\cdot|i-\Tilde{i}+\lambda_{\Tilde{i}}-j-2|}{|i-\Tilde{i}+\lambda_{\Tilde{i}}-j-1|\cdot|i-\Tilde{i}+\lambda_{\Tilde{i}}-j-1|}\\
    &=\frac{|i-\Tilde{i}-1|\cdot|i-\Tilde{i}+\lambda_{\Tilde{i}}-j|}{|i-\Tilde{i}|\cdot|i-\Tilde{i}+\lambda_{\Tilde{i}}-j-1|}\\
     &=\frac{|i-\Tilde{i}-1|}{|i-\Tilde{i}|}\cdot\frac{h'(\Tilde{i},j)}{h(\Tilde{i},j)}.
\end{align*}
We multiply now the last expression for all boxes $(\Tilde{i}, j)$ with $\Tilde{i}=0,\dots, i-1$
\begin{align*}
    \prod_{\Tilde{i}=0}^{i-1}\frac{|i-\Tilde{i}-1|}{|i-\Tilde{i}|}\cdot\frac{h'(\Tilde{i},j)}{h(\Tilde{i},j)}&= \frac{|i-1|}{|i|}\cdot  \frac{|i-2|}{|i-1|}\cdots  \frac{|0|}{|1|}\cdot\prod_{\Tilde{i}=0}^{i-1}\frac{h'(\Tilde{i},j)}{h(\Tilde{i},j)}\\
    &=\frac{1}{|i|}\cdot\prod_{\Tilde{i}=0}^{i-1}\frac{h'(\Tilde{i},j)}{h(\Tilde{i},j)}.
\end{align*}
This is the contribution of all boxes $(l,k)\in \lambda$ such that $k\geq j$. By symmetry, we get that the contribution of all boxes $(l,k)\in \lambda$ such that $l\geq i$ is given by 
\begin{align*}
    \frac{1}{|j|}\cdot\prod_{\Tilde{j}=0}^{j-1}\frac{h'(i,\Tilde{j})}{h(i,\Tilde{j})}.
\end{align*}
{\bf Step II:}
The contribution of the remaining boxes is given by
\begin{multline*}
   |i-j|\cdot  \prod_{\substack{(l,k)\in \lambda\\ l<i \\ k<j}} \frac{ |1+i-j+k-l|}{|i-j+k-l|}\cdot\frac{ |-1+i-j+k-l|}{|i-j+k-l|} \\
   =|i-j|\cdot \prod_{l=0}^{i-1}\prod_{k=0}^{j-1}\frac{ |1+i-j+k-l|}{|i-j+k-l|}\cdot\frac{ |-1+i-j+k-l|}{|i-j+k-l|} \\
   =|i-j|\cdot \prod_{l=0}^{i-1}\frac{|i-l-j-1|\cdot |i-l|}{|i-l-j|\cdot |i-l-1|}
   =|i-j|\cdot\frac{|i|\cdot |j|}{|i-j|}=|i|\cdot |j|.
\end{multline*}
{\bf Step III:} Using Step I,II and the induction step we have
\begin{align*}
    H_{\lambda'}&=\prod_{\Box\in \lambda}h(\Box)\cdot \prod_{\Tilde{i}=0}^{i-1}\frac{h'(\Tilde{i},j)}{h(\Tilde{i},j)}\cdot\prod_{\Tilde{j}=0}^{j-1}\frac{h'(i,\Tilde{j})}{h(i,\Tilde{j})}\\
    &=\prod_{\Box\in \lambda'}h'(\Box),
\end{align*}
which concludes the proof.
\end{proof}
\begin{lemma}\label{lemma: sign for leading term}
Let $\lambda$ be a Young diagram and $k_1, k_2\in \BZ$. Set 
\begin{multline*}
    \rho(\lambda)=\sum_{\substack{(i,j)\in \lambda\\(i,j)\neq (a,a)}}(i k_1+jk_2)+ \sum_{\substack{(i,j), (l,k)\in \lambda\\ (i,j)\neq (l+a-1,k+a)}} ((i-l+1)k_1+(j-k)k_2)\\+\sum_{\substack{(i,j), (l,k)\in \lambda\\ (i,j)\neq (l+a,k+a)}} ((i-l)k_1+(j-k)k_2).
\end{multline*}
We have
\begin{align*}
    (-1)^{\rho(\lambda)+ |\lambda|(k_1+k_2)}\sigma(\lambda)^{k_1+k_2}=(-1)^{|\lambda|k_2+ n(\lambda)k_1+n(\overline{\lambda})k_2},
\end{align*}
where $\sigma(\lambda)$ was defined in Lemma \ref{lemma: hook length in size 0}.
\end{lemma}
\begin{proof}
Let $\equiv$ denote congruence modulo 2. We have
\begin{multline*}
    \rho(\lambda)+ |\lambda|(k_1+k_2)\equiv\sum_{(i,j)\in \lambda}(i k_1+jk_2)+ (k_1+k_2)\left(\sum_{\substack{(i,j)\in \lambda\\(i,j)= (a,a)}}a+ \sum_{\substack{(i,j), (l,k)\in \lambda\\ (i,j)= (l+a-1,k+a)}} a+\sum_{\substack{(i,j), (l,k)\in \lambda\\ (i,j)= (l+a,k+a)}} a\right)\\
     +k_1\sum_{(i,j), (l.k)\in \lambda}1+ |\lambda|(k_1+k_2)\\
     \equiv n(\lambda)k_1+n(\overline{\lambda})k_2+ |\lambda|k_2
     +(k_1+k_2)\left(\sum_{\substack{(i,j)\in \lambda\\(i,j)= (a,a)}}a+ \sum_{\substack{(i,j), (l,k)\in \lambda\\ (i,j)= (l+a-1,k+a)}} a+\sum_{\substack{(i,j), (l,k)\in \lambda\\ (i,j)= (l+a,k+a)}} a\right),
\end{multline*}
therefore the statement of the lemma reduces to the following claim
\begin{align}\label{eqn:claim in proof signs of leading term}
    (-1)^{\sum_{\substack{(i,j)\in \lambda\\(i,j)= (a,a)}}a+ \sum_{\substack{(i,j), (l,k)\in \lambda\\ (i,j)= (l+a-1,k+a)}} a+\sum_{\substack{(i,j), (l,k)\in \lambda\\ (i,j)= (l+a,k+a)}} a}\cdot \sigma(\lambda)=1.
\end{align}
Given a lattice point $\mu=(\mu_1, \mu_2)\in \BZ^2$ define 
\begin{align*}
    \tau(\mu)=\begin{cases}(-1)^{\mu_1} &\mu_1=\mu_2,\\
    \sgn(\mu_1-\mu_2) &\mu_1\neq \mu_2.
    \end{cases}
\end{align*}
Notice that the left-hand-side of \eqref{eqn:claim in proof signs of leading term} can be rewritten as
\begin{align*}
    \prod_{(i,j)\in \lambda}\tau(-i,-j)\prod_{(i,j), (l,k)\in \lambda}\tau(i-l,j-k)\tau(1+i-l,j-k),
\end{align*}
and denote this last expression by $F_\lambda$. We prove the claim \eqref{eqn:claim in proof signs of leading term} on induction on the size of $ \lambda$. If $|\lambda|=1$, the result is clear. Assume it holds for all Young diagrams of size $n$ and consider a Young diagram $\lambda'$ of size $n+1$ obtained from a Young diagram of size $\lambda$ by adding a box with lattice coordinates $(i,j)\in \BZ^2$; we have
\begin{align*}
    F_{\lambda'}=F_\lambda\cdot \tau(-i,-j)\prod_{(l,k)\in \lambda}\tau(i-l,j-k)\tau(l-i, k-j)\tau(1+i-l,j-k)\tau(1+l-i, k-j).
\end{align*}
We analyze the contribution of every box $(l,k)\in \lambda$ in the product above. We say that a box $(l,k)$ is in the same diagonal as $(i,j)$ if it is of the form
\begin{align*}
    (l,k)=(i+a,j+a), \quad a\in \BZ.
\end{align*}
The contribution of the boxes on the same diagonal of $(i,j)$ is 
\begin{align*}
    \tau(-a,-a)\tau(a,a)\tau(1-a,-a)\tau(1+a,a)=1.
\end{align*}
 We say that a box $(l,k)$ is in the \emph{lower diagonal} of $(i,j)$ if it is of the form 
\begin{align*}
    (l,k)=(i+a,j+a-1), \quad a\in \BZ.
\end{align*}
The contribution of the boxes in the lower diagonal of $(i,j)$ is
\begin{align*}
    \tau(-a,1-a)\tau(a,a-1)\tau(1-a,1-a)\tau(1+a,a-1)=(-1)^a.
    \end{align*}
 We say that a box $(l,k)$ is in the \emph{upper diagonal} of $(i,j)$ if it is of the form 
\begin{align*}
    (l,k)=(i+a-1,j+a), \quad a\in \BZ.
\end{align*}
The contribution of the boxes in the upper diagonal of $(i,j)$ is
\begin{align*}
    \tau(1-a,-a)\tau(a-1,a)\tau(2-a,-a)\tau(a,a)=(-1)^{a+1}.
    \end{align*}
    A completely analogous analysis shows that all other boxes $(i,j)\in \lambda$ do not contribute to the product. Therefore the contribution to the sign in the product is just given by the boxes in the upper or lower  diagonal of $(i,j)$, as displayed in the picture
\begin{figure}[!h]
    \centering
 \young(11111,1\minusone 1\minusone,111\ijbox,11,1)
 \caption{The label $\pm 1$ in every box represents the contribution to the final sign, with $(i,j)=(2,3)$. }
\end{figure}

 If we denote by $\delta_{+}, \delta_{-}$ respectively the number of boxes in the upper diagonal and lower diagonal, we conclude that
    \begin{align*}
       \prod_{(l,k)\in \lambda}\tau(i-l,j-k)\tau(l-i, k-j)\tau(1+i-l,j-k)\tau(1+l-i, k-j)= (-1)^{\lceil\frac{\delta_{+}}{2}\rceil+\lfloor\frac{\delta_{-}}{2}\rfloor},
    \end{align*}
    where $ \lceil\cdot\rceil,\lfloor \cdot\rfloor$ denote the ususal \emph{ceiling} and \emph{floor} functions. One readily proves that 
    \begin{align*}
           (\delta_+, \delta_{-})= \begin{cases}
        (j+1,j)& i>j,\\
        (i,i) & i=j,\\
        (i,i+1) & i<j.
        \end{cases}
    \end{align*}
            With a direct analysis we can show that
        \begin{align*}
           (-1)^{\lceil\frac{\delta_{+}}{2}\rceil+\lfloor\frac{\delta_{-}}{2}\rfloor}\cdot\tau(-i,-j)=1
        \end{align*}
       in all the three cases, by which we conclude the inductive step.
\end{proof}

\bibliographystyle{amsplain-nodash}
\bibliography{The_Bible}
\end{document}

%% file: macros.tex
\usepackage[colorinlistoftodos]{todonotes}

\definecolor{antiquewhite}{rgb}{0.98, 0.92, 0.84}
\definecolor{buff}{rgb}{0.94, 0.86, 0.51}
\definecolor{palecopper}{rgb}{0.85, 0.54, 0.4}
\definecolor{fluorescentyellow}{rgb}{0.8, 1.0, 0.0}

\definecolor{britishracinggreen}{rgb}{0.0, 0.26, 0.15}
\definecolor{cobalt}{rgb}{0.0, 0.28, 0.67}
\DeclareSymbolFont{usualmathcal}{OMS}{cmsy}{m}{n}
\DeclareSymbolFontAlphabet{\mathcal}{usualmathcal}

\newcommand{\TT}{\mathbf{T}}

\newcommand{\BC}{{\mathbb{C}}}

\newcommand{\BE}{{\mathbb{E}}}
\newcommand{\BF}{{\mathbb{F}}}

\newcommand{\BI}{{\mathbb{I}}}

\newcommand{\BL}{{\mathbb{L}}}

\newcommand{\BP}{{\mathbb{P}}}
\newcommand{\BQ}{{\mathbb{Q}}}

\newcommand{\BZ}{{\mathbb{Z}}}

\newcommand{\CE}{{\mathcal E}}
\newcommand{\CF}{{\mathcal F}}
\newcommand{\CG}{{\mathcal G}}

\newcommand{\CI}{{\mathcal I}}

\newcommand{\CK}{{\mathcal K}}

\newcommand{\CP}{{\mathcal P}}

\newcommand{\CW}{{\mathcal W}}
\newcommand{\CX}{{\mathcal X}}
\newcommand{\CY}{{\mathcal Y}}
\newcommand{\CZ}{{\mathcal Z}}

\newcommand{\pt}{{\mathsf{pt}}}
\newcommand{\ch}{{\mathrm{ch}}}
\newcommand{\td}{{\mathrm{td}}}
\newcommand{\fix}{\mathsf{fix}}
\newcommand{\mov}{\mathsf{mov}}

\DeclareMathOperator{\Hilb}{Hilb}
\DeclareMathOperator{\Sets}{Sets}
\DeclareMathOperator{\Sch}{Sch}

\DeclareMathOperator{\Sym}{Sym}

\DeclareMathOperator{\vd}{vd}

\DeclareMathOperator{\vir}{\mathrm{vir}}

\DeclareMathOperator{\Var}{Var}

\DeclareMathOperator{\Db}{D^b}

\DeclareMathOperator{\rk}{rk}
\DeclareMathOperator{\op}{op}
\DeclareMathOperator{\tr}{tr}

\newcommand{\derived}{\mathbf{D}}
\newcommand{\dd}{\mathrm{d}}

\newcommand{\into}{\hookrightarrow}

\DeclareFontFamily{OT1}{rsfs}{}
\DeclareFontShape{OT1}{rsfs}{n}{it}{<-> rsfs10}{}
\DeclareMathAlphabet{\curly}{OT1}{rsfs}{n}{it}
\renewcommand\hom{\curly H\!om}

\newcommand\Hom{\operatorname{Hom}}

\newcommand{\RR}{\mathbf R}

\newcommand\Tot{\operatorname{Tot}}

\newcommand\id{\operatorname{id}}

\newcommand{\Pic}{\mathop{\rm Pic}\nolimits}
\newcommand{\PT}{\mathsf{PT}}

\newcommand{\GW}{\mathsf{GW}}



\usepackage[all]{xy}
\usepackage{tikz}
\usepackage{tikz-cd}
\usepackage{adjustbox}
\usepackage{rotating}
\usepackage{comment}

\usetikzlibrary{matrix,shapes,intersections,arrows,decorations.pathmorphing}
\tikzset{commutative diagrams/arrow style=math font}
\tikzset{commutative diagrams/.cd,
mysymbol/.style={start anchor=center,end anchor=center,draw=none}}

\tikzset{
shift up/.style={
to path={([yshift=#1]\tikztostart.east) -- ([yshift=#1]\tikztotarget.west) \tikztonodes}
}
}

\theoremstyle{definition}

\newtheorem*{lemma*}{Lemma}
\newtheorem*{theorem*}{Theorem}
\newtheorem*{example*}{Example}
\newtheorem*{fact*}{Fact}
\newtheorem*{notation*}{Notation}
\newtheorem*{definition*}{Definition}
\newtheorem*{prop*}{Proposition}
\newtheorem*{remark*}{Remark}
\newtheorem*{corollary*}{Corollary}

\newtheorem*{conventions*}{Conventions}

\newtheorem{definition}{Definition}[section]

\newtheorem{example}[definition]{Example}

\newtheorem{remark}[definition]{Remark}

\newtheoremstyle{thm} 
        {3mm}
        {3mm}
        {\slshape}
        {0mm}
        {\bfseries}
        {.}
        {1mm}
        {}
\theoremstyle{thm}
\newtheorem{theorem}[definition]{Theorem}
\newtheorem{corollary}[definition]{Corollary}
\newtheorem{lemma}[definition]{Lemma}
\newtheorem{prop}[definition]{Proposition}

\newtheoremstyle{ex} 
        {3mm}
        {3mm}
        {}
        {0mm}
        {\scshape}
        {.}
        {1mm}
        {}
\theoremstyle{ex}

\newtheoremstyle{sol} 
        {3mm}
        {3mm}
        {}
        {0mm}
        {\scshape}
        {.}
        {1mm}
        {}
\theoremstyle{sol}

\usepackage{amsmath,float}
\usetikzlibrary{decorations.markings}

\usepackage{amssymb}
\usepackage{enumitem}
\usepackage{comment}
\newcommand{\sgn}{{\mathrm{sgn}}}
\usepackage{bm}

   \DeclareMathOperator{\oO}{\mathcal{O}}
   
   \DeclareMathOperator{\tf}{\mathfrak{t}}

\newcommand\AJ{\operatorname{AJ}}